\RequirePackage{fix-cm}
\documentclass[twocolumn]{svjour3}          
\smartqed  
\usepackage{subfigure}
\usepackage{graphicx}
\usepackage{booktabs}
\usepackage{caption}

\newlength{\tab}
\usepackage{amscd}
\usepackage{amssymb}
\usepackage{amsmath}
\usepackage{mathrsfs}
\newcommand{\R}{\mathbb{R}}

\newcommand{\ul}{\textbf}

\newcommand{\www}{\mbox{\boldmath$\omega$}}

\newcommand{\KKK}{\mbox{\boldmath$\kappa$}}
\newcommand{\TTT}{\mbox{\boldmath$\tau$}}
\newcommand{\OOmega}{\mbox{\boldmath$\Omega$}}
\newcommand{\desda}{\Leftrightarrow}

\DeclareMathAlphabet\gothic{U}{euf}{m}{n}
\twocolumn

\newcommand{\argmin}{\mathop{\mathrm{arg\,min}}\limits}

\begin{document}
\title{Locally Adaptive Frames in the Roto-Translation Group and their Applications
in Medical Imaging}

\titlerunning{Locally Adaptive Frames in the Roto-Translation Group}        

\author{R. Duits$^*$ \and  M.H.J. Janssen$^*$ \and J. Hannink \and G.R. Sanguinetti}

\authorrunning{Duits, Janssen, Hannink, Sanguinetti}
%

\institute{{*\textbf{Joint main authors:}}\\[8pt]
R. Duits \at  CASA,
              Eindhoven University of Technology, \\
              Tel.: +31-40-2472859,
              \email{R.Duits@tue.nl}  \\         
                  \and
           M.H.J. Janssen \at
              CASA,
              Eindhoven University of Technology, \\
              Tel.: +31-40-2475571,        %
              \email{M.H.J.Janssen@tue.nl} \\ \\
              \textbf{co-authors:} \\ \\
              J. Hannink \at Digital Sports Group, \\
              University of Erlangen-N\"urnberg (FAU),\\
              Tel.: +49-9131-85-27830,
              \email{julius.hannink@fau.de}
              \and
              G.R. Sanguinetti \at
              CASA,
              Eindhoven University of Technology, \\
              Tel.: +31-40-2478777,        %
              \email{G.R.Sanguinetti@tue.nl}
               }       


\maketitle

\begin{abstract}
Locally adaptive differential frames (gauge frames) are a well-known effective tool in image analysis, used in differential invariants and PDE-flows. However, at complex structures such as crossings or junctions, these frames are not well-defined. Therefore, we generalize the notion of gauge frames on images to gauge frames on data representations $U:\R^{d} \rtimes S^{d-1} \to \mathbb{R}$ defined on the extended space of positions and orientations, which we relate to data on the roto-translation group $SE(d)$, $d=2,3$. This allows to define multiple frames per position, one per orientation. We compute these frames via exponential curve fits in the extended data representations in $SE(d)$. These curve fits minimize first or second order variational problems which are solved by spectral decomposition of, respectively, a structure tensor or Hessian of data  on $SE(d)$. We include these gauge frames in differential invariants and crossing preserving PDE-flows acting on extended data representation $U$ and 
we show their advantage compared to the standard left-invariant frame on $SE(d)$. Applications include crossing-preserving filtering and improved segmentations of the vascular tree in retinal images, and new 3D extensions of coherence-enhancing diffusion via invertible orientation scores.

\end{abstract}

\section{Introduction \label{ch:intro}}
\label{intro}
\begin{figure}
\includegraphics[width=\hsize]{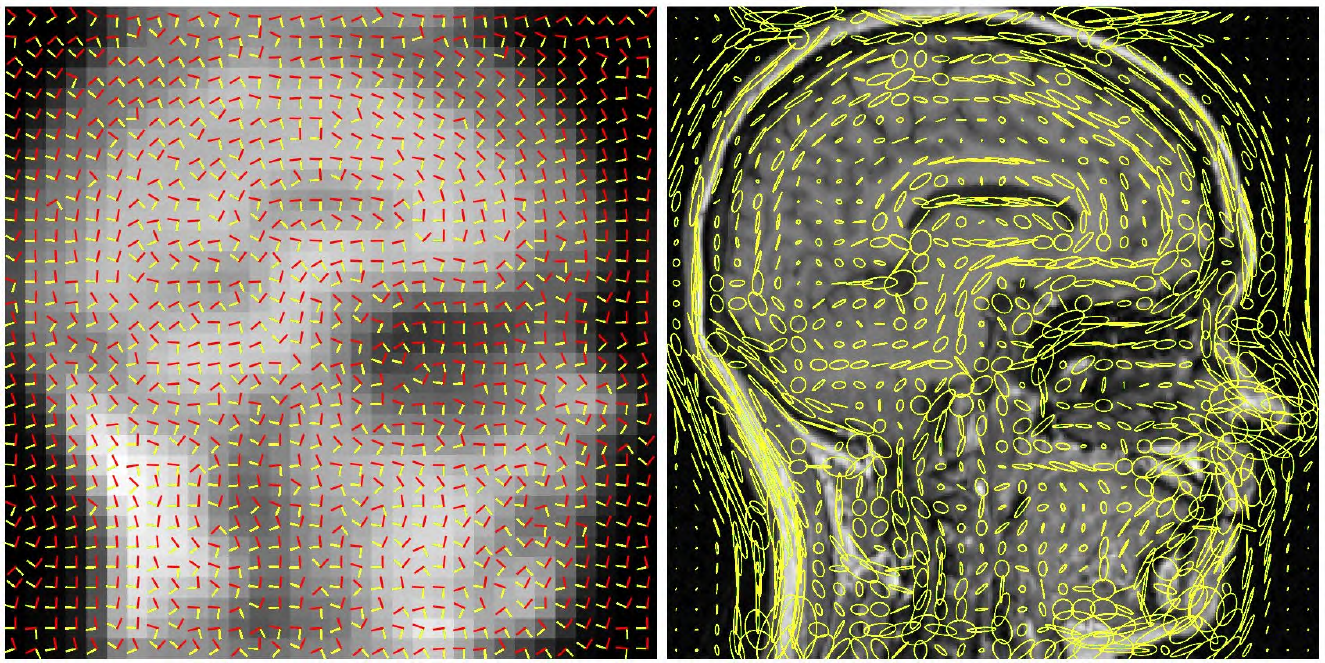}
\caption{Left: Locally adaptive frames (gauge frames) in the image domain computed as the eigenvectors of the Hessian of the image at each location.
Right: Such gauge frames can be used for adaptive anisotropic diffusion and geometric reasoning. However, at complex structures such as blob-structures/crossings, the gauge frames are ill-defined causing fluctuations. \label{fig:1}}
\end{figure}

Many existing image analysis techniques rely on differential frames that are locally adapted to image data. This includes methods based on differential invariants \cite{Sapiro,HaarRomenybook,Florackbook,Lindeberg}, partial differential equations \cite{Sapiro,Weickert,Guichard}, and non-linear and morphological scale spaces \cite{Burgeth,Breuss,Welk}, used in various image processing tasks such as tracking and line detection \cite{Bekkers}, corner detection and edge focussing \cite{HaarRomenybook,Bergholm}, segmentation \cite{Staal}, active contours \cite{Cao,Kimmel}, DTI data processing \cite{KindlmannVIZ,KindlmannTMI}, feature based clustering etc. These local coordinate  frames (also known as `gauge frames' according to \cite{Florackbook,Blom,HaarRomenybook}) provide differential frames directly adapted to the local image structure via a structure tensor or a Hessian of the image. Typically the structure tensor (based on 1st order Gaussian derivatives) is used for adapting to edge-like structures while the Hessian (based on 2nd order Gaussian derivatives) is used for adapting to line-like structures.  The primary benefit of the gauge frames is that they allow to include adaptation for anisotropy and curvature in a rotation and translation invariant way. 
See Fig.~\ref{fig:1}, where we have depicted local adaptive frames based on eigenvector decomposition of the image Hessian at some given scale, of the MR-image in the background.

\begin{figure}
\includegraphics[width=\hsize]{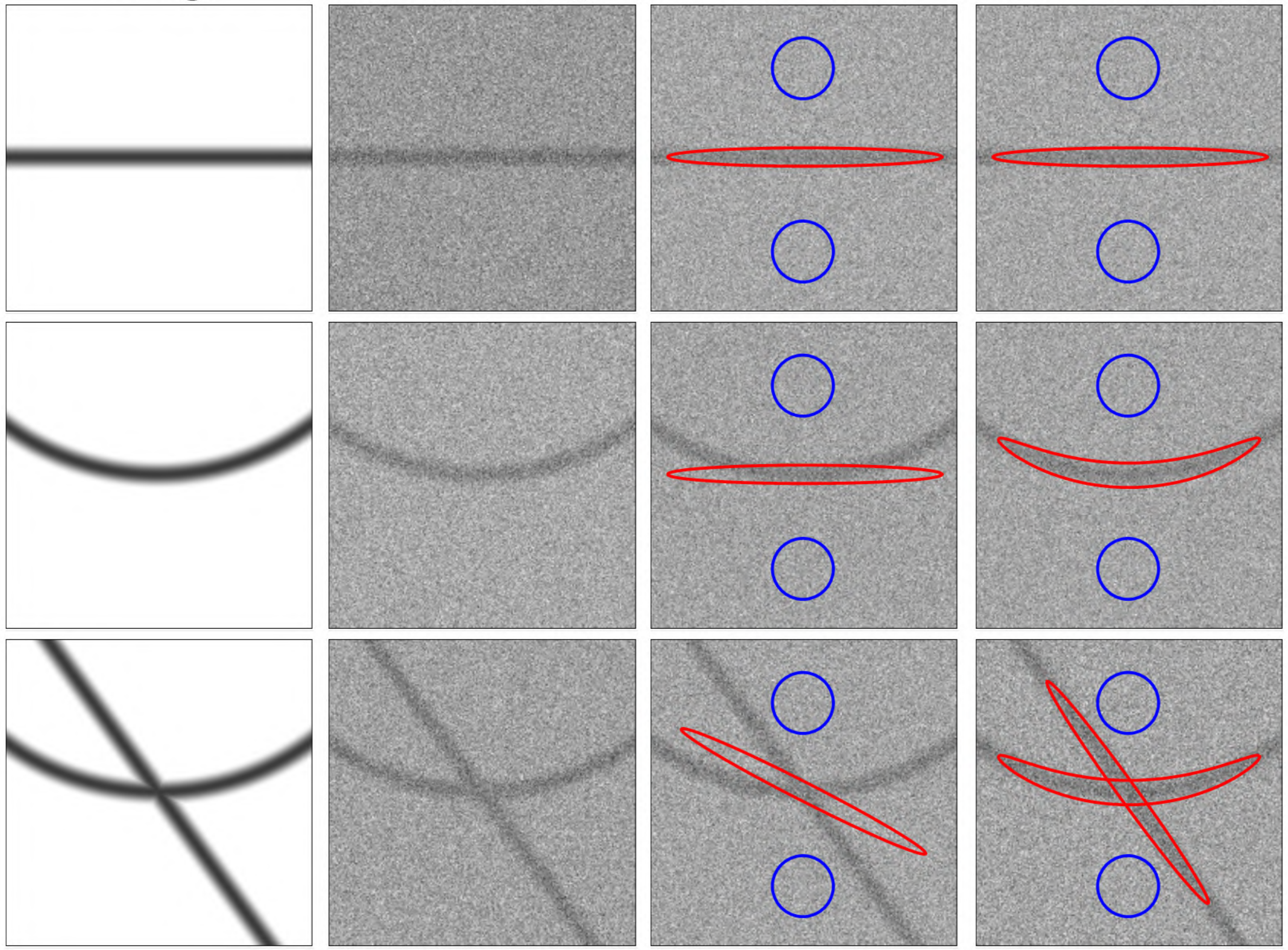}
\caption{We aim for adaptive anisotropic diffusion of images 
which takes into account curvature.
At areas with low orientation confidence (in blue)
isotropic diffusion is required, whereas at areas with high orientation confidence (in red) anisotropic diffusion with curvature adaptation is required. Application of locally adaptive frames in the image domain suffers from interference (3rd column), whereas application of locally adaptive frames in the domain $\R^{d} \rtimes S^{d-1}$
allows for adaptation along all the elongated structures (4th column).
\label{fig:2}}
\end{figure}

It is sometimes problematic that such locally adapted differential frames are directly placed in the image domain $\R^{d}$ $(d=2,3)$, as at the vicinity of complex structures, e.g. crossings, textures, bifurcations, one typically requires
multiple local spatial coordinate frames. To this end, one effective alternative is to extend the image domain to the joint space of positions and orientations $\R^{d} \rtimes S^{d-1}$.  The advantage is that it allows to disentangle oriented structures involved in crossings, and to include curvature, cf.~\!Fig.~\!\ref{fig:2}. Such extended domain techniques rely on various kinds of lifting, such as coherent state transforms (also known as invertible orientation scores) \cite{Ali,QAM1,Bekkers,FrankenPhDThesis}, continuous wavelet transforms \cite{DuitsIJCV,QAM1,Sharma,Bekkers}, orientation lifts \cite{Zweck,Boscain}, or orientation channel representations \cite{Felsberg}. In case one has to deal with more complex diffusion weighted MRI techniques, the data in extended position orientation domain can be obtained after a modelling procedure as in  \cite{CSD,Tuch,Aganj,Sinnaeve}. In this article we will not discuss in detail on how such a new image representation or lift $U:\R^{d} \rtimes S^{d-1} \to \mathbb{R}$ is to be constructed from grey-scale image $f:\R^{d} \to \R$, and we assume it to be a sufficiently smooth given input. Here  $U(\ul{x},\ul{n})$ is to be considered as a probability density of finding a local oriented structure
(i.e. an elongated structure) at position $\ul{x} \in \R^{d}$ with orientation $\ul{n} \in S^{d-1}$. 

When processing data in the extended position orientation domain it is often necessary to equip the domain with a structure that links the data across different orientation channels, in such a way that a notion of alignment between local orientations is taken into account. This is achieved by relating data on positions and orientations to data on the roto-translation group $SE(d)=\R^{d} \rtimes SO(d)$. 
This idea resulted in contextual image analysis methods
\cite{Mumford,Zweck,August,Duits2006,CittiFranceschiello,QAM1,Sharma,vanAlmsickthesis,DuitsIJCV2010,Chirikjian,Siddiqi} and appears in models of low level visual perception and their relation with the functional architecture of the visual cortex \cite{Citti,Petitot,Sanguinetti2010,Boscain,Barbieri,SanguinettiPhd,Zucker}.
Following the conventions in \cite{DuitsIJCV2010} we denote functions on the coupled space of positions and orientations by $U:\mathbb{R}^d\rtimes S^{d-1}\to \mathbb{R}$. Then, its extension $\tilde{U}:SE(d)\to\mathbb{R}$ is given by:
\begin{equation} \label{tildeU}
\tilde{U}(\ul{x},\ul{R}):=U(\ul{x},\ul{R}\ul{a})
\end{equation}
for all $\ul{x}\in \mathbb{R}^d$ and all rotations $\ul{R}\in SO(d)$, and given
reference axis $\ul{a} \in S^{d-1}$. Throughout this article $\ul{a}$ is chosen as follows:
\begin{equation} \label{conventie}
d=2 \Rightarrow \ul{a}=(1,0)^T, \ d=3 \Rightarrow \ul{a}=(0,0,1)^T.
\end{equation}
Then, we can identify the joint space of positions and orientations $\mathbb{R}^d\rtimes S^{d-1}$ by:
\begin{equation} \label{PositionsAndOrientations}
\R^{d} \rtimes S^{d-1}:= SE(d)/(\{\ul{0}\} \times SO(d-1)),
\end{equation}
where this quotient structure is due to (\ref{tildeU}), and where $SO(d-1)$ is identified with all rotations on $\R^d$ that map reference axis $\ul{a}$ onto itself.
Note that in Eq. (\ref{tildeU}) the tilde indicates we consider data on the group instead of data on the quotient.
If $d=2$ the tildes can be ignored as $\mathbb{R}^2\rtimes S^1 =SE(2)$. However, for $d\geq 3$ this distinction is crucial and necessary details on (\ref{PositionsAndOrientations}) will follow in the beginning of Section~\ref{ch:SE3}.

In this article, our quest is to find locally optimal differential frames in $SE(d)$  relying on similar Hessian- and/or structure-tensor type of techniques for gauge frames on images, recall Fig.~\!\ref{fig:1}. Then, the frames can be used to construct crossing-preserving differential invariants and adaptive diffusions of data in $SE(d)$. In order to find these optimal frames our main tool is the theory of curve fits. Early works on curve fits have been presented in \cite{Parent1989} where the notion of \emph{curvature consistency} is applied to inferring local curve orientations, based on neighbourhood co-circularity continuation criteria. This approach was extended to 2D texture flow inference in \cite{BenShahar2003}, by lifting images in position and orientation domain and inferring multiple Cartan frames at each point. 
Our work is embedded in a Lie group framework where we consider the notion of \emph{exponential curve} fits via formal variational methods. Exponential curves in the  $SE(d)$-curved geometry  are the equivalents of \emph{straight}\footnote{Exponential curves are auto-parallels w.r.t.~`-'Cartan~connection, see Appendix A, Eq. (\ref{AutoParallel}).} lines in the Euclidean geometry.  If $d=2$, the spatial projection of these exponential curves are osculating circles, which are used: for constructing the curvature consistency
in \cite{Parent1989},  for defining the tensor voting fields in \cite{medioni},
and for local modeling association fields in \cite{Citti}. If $d=3$, the spatial projection of exponential curves are spirals with constant curvature and torsion. Based on co-helicity principles, similar spirals have been used in neuroimaging applications \cite{SavadjievMEDIA2006} or for modelling heart fibers \cite{SavadjievPNAS2012}. In these works curve fits are obtained via efficient discrete optimization techniques, which are
beyond the scope of this article.

\begin{figure*}
\centering
\includegraphics[width=.9\hsize]{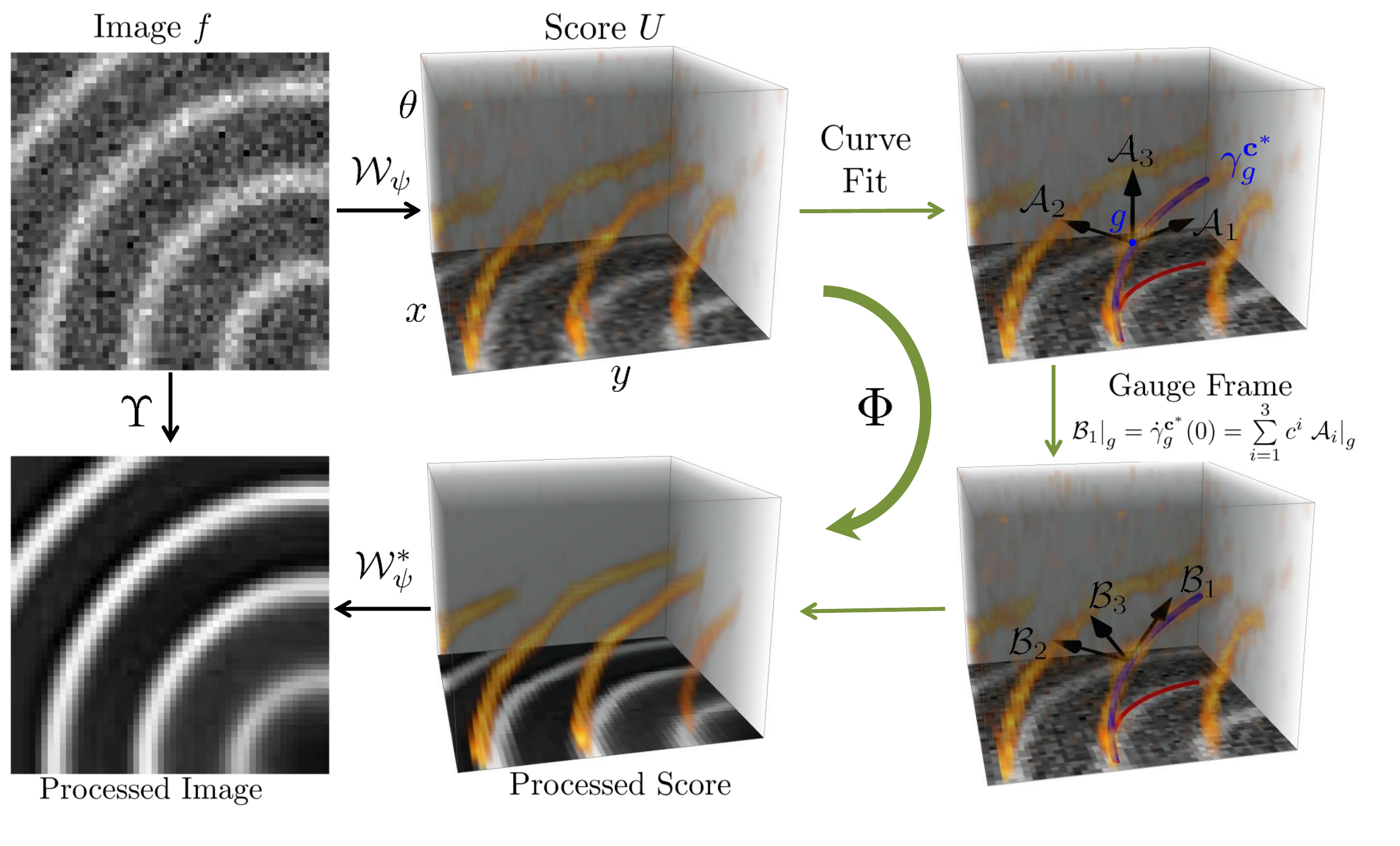}
\caption{The overall pipeline of image processing $f \mapsto \Upsilon f$ via left-invariant operators $\Phi$.  In this pipeline we construct an invertible orientation score $W_\psi f$ (Section \ref{sec:ios}), we fit an exponential curve (Section \ref{ch:SE2},\ref{ch:SE3}), we obtain the gauge frame (Section \ref{ch:togauge} and App. \ref{app:a}), we construct a non-linear diffusion, and finally we apply reconstruction (Section \ref{sec:ios}). The main focus of this paper is curve fitting, where we compute per element $g=(x,y,\theta)$
an exponential curve fit $\gamma^{\ul{c}^*}_{g}(t)$ (in blue, with spatial projection in red) with tangent
$\dot{\gamma}_g^{\ul{c}^{*}}(0)=\ul{c}^{*}(g)=(c^1,c^2,c^3)^T$ at $g$. Based on this fit we construct for each $g$ a local frame $\{\left.\mathcal{B}_{1}\right|_{g},
\left.\mathcal{B}_{2}\right|_{g},\left.\mathcal{B}_{3}\right|_{g}\}$ which are used in our operators $\Phi$ on the lift (here $\Phi$ is a non-linear diffusion operator).
\label{Fig:Intro}}
\end{figure*}

In Fig.\!~\ref{Fig:Intro}, we present an example for $d=2$ of the overall pipeline of including locally adaptive frames in a suitable diffusion operators $\Phi$ acting in the lifted domain $\R^{2}\rtimes S^1$. 
For $d>2$ the same pipeline applies. Here, an exponential curve fit  $\gamma^{\ul{c}^*}_{g}(t)$  (in blue, with spatial projection in red) at a group element $g\in SE(d)$ is characterized by $(g, \ul{c}^{*}(g))$, i.e. a starting point $g$ and an tangent vector $\ul{c}^{*}(g)$ that should be aligned with the structures of interest. In essence, this paper explains in detail how to compute $\ul{c}^{*}(g)$ as this will be the principal direction the differential frame will be aligned with, and then gives appropriate conditions for fixing the remaining directions in the frame. 

The main contribution of this article is to provide a general theory for finding locally adaptive frames in the roto-translation group $SE(d)$, for $d=2,3$.
Some preliminary work on exponential curve fits of the second order on $SE(2)$ has been presented in \cite{FrankenPhDThesis,FrankenIJCV,Sharma}. In this paper we formalize these previous methods (Theorems \ref{th:0b} and \ref{th:0c})  and we extend them to first-order exponential curve fits (Theorem \ref{th:0}). Furthermore, we generalize both approaches to the case $d=3$ (Theorems \ref{th:2}, \ref{th:3b}, \ref{th:6}, \ref{app:th} and  \ref{th:8}). All theorems contain new results except for Theorems \ref{th:0b} and \ref{th:0c}. The key ingredient is to consider the  fits as formal variational curve optimization problems with exact solutions derived by spectral decomposition of structure tensors and Hessians of the data $\tilde{U}$ on $SE(d)$. In the $SE(3)$-case we show that in order to obtain torsion-free exponential curve fits with well-posed projection on $\R^{3}\rtimes S^{2}$, one must resign to a two-fold optimization algorithm. To show the potential of considering  these locally adaptive frames, we employ them in  medical image analysis applications, in improved differential invariants and improved crossing-preserving diffusions. Here, we provide for the first time coherence enhancing diffusions via 3D invertible orientation scores \cite{Janssen,MichielMaster}, extending previous methods \cite{FrankenPhDThesis,FrankenIJCV,Sharma} to the 3D Euclidean motion group.

\subsection{Structure of the Article \label{ch:structure}}

We start the body of this article reviewing preliminary differential geometry tools in Section \ref{ch:geomTools}. Then, in Section \ref{ch:togauge} we describe how a given exponential curve fit induces the locally adaptive frame.
In Section \ref{ch:gauge} we provide an introduction by reformulating the standard gauge frames construction in images in a group theoretical setting. This gives a roadmap towards $SE(2)$-extensions explained in Section~\ref{ch:SE2}, where we deal with  exponential curve fits of the 1st order in Subsection~\ref{ch:SE2-order1} computed via a structure tensor, and exponential curves fits of 2nd order in Section~\ref{ch:SE2-order2} computed via the Hessian of the data $\tilde{U}$. In the latter case we have 2 options for the curve optimization problem, one solved by the symmetric sum, and one by the symmetric product of the non-symmetric Hessian.
The curve fits in $SE(2)$ in Section~\ref{ch:SE2}, are extended to curve fits in $SE(3)$ in Section~\ref{ch:SE3}. It starts with preliminaries on the quotient (\ref{PositionsAndOrientations}) and then it follows the same structure as the previous section. Here we present the two-fold algorithm for computing the torsion free exponential curve fits. 

In Section~\ref{ch:app} we consider experiments regarding medical imaging applications and feasibility studies. 
 We first recall the theory of invertible orientation scores needed for the applications. In the $SE(2)$-case we present crossing-preserving multi-scale vessel enhancing filters in retinal imaging, and  in the $SE(3)$-case we include a proof of concept of crossing-preserving (coherence enhancing diffusion) steered by gauge frames via invertible 3D orientation scores.

Finally, there are 5 appendices. Appendix \ref{app:a} supplements Section \ref{ch:togauge} by explaining the construction of the frame for $d=2,3$. Appendix \ref{app:B} describes the geometry of neighboring exponential curves needed for formulating the variational problems. Appendix \ref{app:C}  complements the two-fold approach in Section \ref{ch:SE3}. Appendix \ref{app:remcohappy} provides the definition of the Hessian used in the paper. Finally, Appendix \ref{app:TableOfNotations} contains a list of symbols, their explanations and references to the equation in which they are defined. We advise the reader to keep track of this table. Especially, in the more technical sections: Section \ref{ch:SE2} and \ref{ch:SE3}.

\section{Differential Geometrical Tools \label{ch:geomTools}}
Relating our data to data on the Euclidean motion group, via Eq.~\!(\ref{tildeU}), allows us to use tools from Lie group theory and differential geometry.
In this section we explain these tools that are important for our notion of an exponential curve fit to smooth data $\tilde{U}:SE(d) \to \mathbb{R}$.
Often, we consider the case $d=2$ for basic illustration.
Later on, in Section~\ref{ch:SE3}, we consider the case $d=3$ and extra technicalities on the quotient structure will enter. 




\subsection{The Roto-Translation Group}\label{ch:SEd}

The data $\tilde{U}:SE(d) \to \R$ is defined on the group $SE(d)$ of rotations and translations acting on $\R^{d}$. As the concatenation of two rigid body motions is again a rigid body motion, the group $SE(d)$ is equipped with the following group product:
\begin{equation}\label{groupproduct}
\begin{array}{l}
g g'=(\ul{x},\ul{R}) (\ul{x}', \ul{R}')= (\ul{R} \ul{x}'+\ul{x}, \ul{R} \ul{R}'),   \\
\textrm{ with }g=(\ul{x},\ul{R}), \ \ g'=(\ul{x}',\ul{R}') \in SE(d),
\end{array}
\end{equation}
where we recognize the semi-direct product structure $SE(d)=\R^{d} \rtimes SO(d)$, of the translation group $\R^{d}$ with rotation group {\small $SO(d)=\{\ul{R} \in \R^{d\times d}\,|\, \ul{R}^{T}\!=\!\ul{R}^{-1}\!\!,\det \ul{R}\!=\!1\}$}. The groups $SE(d)$ and $SO(d)$ have dimension
\begin{equation}\label{dims}
\begin{array}{ll}
r_d&:=\textrm{dim}(SO(d))=\frac{(d-1)d}{2}, \\
n_d&:=\textrm{dim}(SE(d))=\frac{d(d+1)}{2}=d+r_d.
\end{array}
\end{equation}
Note that $n_2=3$, $n_3=6$.
One may represent elements $g$ from $SE(d)$ by the following matrix representation
\begin{equation}\label{thedw}
\begin{array}{l}
g \equiv M(g)=\left(
\begin{array}{cc}
\ul{R} & \ul{x} \\
\ul{0}^T & 1
\end{array}
\right), \textrm{ which indeed satisfies } \\
M(g \, g')=M(g)\, M(g').
\end{array}
\end{equation}
We will often avoid this embedding into the set of invertible $(d+1) \times (d+1)$ matrices, in order to focus on the geometry rather than the algebra.

\subsection{Left-Invariant Operators}

In image analysis applications operators $\tilde{U} \mapsto \tilde{\Phi}(\tilde{U})$ need to be left-invariant and not right-invariant \cite{DuitsRPHDThesis,FrankenPhDThesis}. Left-invariant operators $\tilde{\Phi}$ in the extended domain correspond to rotation and translation invariant operators $\Upsilon$ in the image domain, which is a desirable property. On the other hand, right-invariance boils down to isotropic operators in the image domain which is an undesirable restriction.
By definition $\tilde{\Phi}$ is left-invariant and not right-invariant if it commutes
with the left-regular representation $\mathcal{L}$ (and not with the right-regular representation $\mathcal{R}$).
Representations $\mathcal{L},\mathcal{R}$ are given by
\begin{equation}\label{eq:leftrightregularrepr}
\begin{array}{l}
(\mathcal{L}_h \tilde{U})(g)=\tilde{U}(h^{-1}g),\ \
(\mathcal{R}_h \tilde{U})(g)=\tilde{U}(gh), \
\end{array}
\end{equation}
for all $h,g \in SE(d)$. So operator $\tilde{\Phi}$ must satisfy $\tilde{\Phi} \circ \mathcal{L}_{g}= \mathcal{L}_{g} \circ \tilde{\Phi}$ and
$\tilde{\Phi} \circ \mathcal{R}_{g}\neq \mathcal{R}_{g} \circ \tilde{\Phi}$
for all $g \in SE(d)$.

\subsection{Left-Invariant Vector Fields and Dual Frame}\label{ch:LeftInvariantDerivatives}

A special case of left-invariant operators are left-invariant derivatives. More precisely (see Remark~\ref{rem:1} below), we need to consider left-invariant vector fields $g \mapsto \mathcal{A}_{g}$, as the left-invariant derivative $\mathcal{A}_{g}$ depends on the location $g$ where it is attached.
Intuitively, the left-invariant vector fields $\{\mathcal{A}_{i}\}_{i=1}^{n_i}$ provide a local moving frame of reference in the tangent bundle $T(SE(d))$, that comes in naturally when including alignment of local orientations in the image processing of $\tilde{U}$.

Formally, the left-invariant vector fields are obtained by taking a basis $\{A_{i}\}_{i=1}^{n_d} \in T_{\gothic{e}}(SE(d))$ in the tangent space at the unity element
$\ \gothic{e}:=(\ul{0},I),\ $
and then one uses
the push-forward $(L_{g})_*$ of the left multiplication
\begin{equation}\label{leftmultiplication}
  L_gh=gh,
\end{equation}
to obtain the corresponding tangent vectors in the tangent space $T_{g}(SE(d))$. Thus one associates to each
$A_{i}$ a left-invariant field $\mathcal{A}_{i}$ given by
\begin{equation} \label{LINVDef}
\left.\mathcal{A}_{i}\right|_{g} = (L_{g})_*A_i, \textrm{ for all }g \in SE(d), \ i=1,\ldots, n_d,
\end{equation}
where we consider each $\mathcal{A}_{i}$ as a differential operator on smooth locally defined functions $\tilde{\phi}$ given by
\[
\left.\mathcal{A}_{i}\right|_{g}\tilde{\phi}= (L_{g})_*A_i\tilde{\phi} :=A_{i}(\tilde{\phi} \circ L_g).
\]
An explicit way to construct and compute the differential operators $\left.\mathcal{A}_{i}\right|_{g}$ from $A_{i}=\left.\mathcal{A}_{i}\right|_{\gothic{e}}$ is via
\begin{equation}\label{LINV}
\left.\mathcal{A}_{i}\right|_{g}\tilde{\phi}=\mathcal{A}_{i}\tilde{\phi}(g)= \lim \limits_{\epsilon \to 0} \frac{\tilde{\phi}(g \, e^{\epsilon A_i})-\tilde{\phi}(g)}{\epsilon},
\end{equation}
where $A \mapsto e^A=\sum \limits_{k=0}^{\infty} \frac{A^k}{k!}$ denotes the matrix exponential from Lie algebra $T_{\gothic{e}}(SE(d))$ to Lie group $SE(d)$.
The differential operators $\{\mathcal{A}_{i}\}_{i=1}^{n_d}$ induce a corresponding dual frame
$\{\omega^{i}\}_{i=1}^{n_d}$, which is a basis for the co-tangent bundle $T^{*}(SE(d))$. This dual frame is given by
\begin{equation} \label{coframe}
\langle \omega^{i}, \mathcal{A}_{j} \rangle =\delta^{i}_{j} \textrm{ for }i,j=1,\ldots n_d,
\end{equation}
where $\delta^{i}_{j}$ denotes the Kronecker delta. Then the
derivative of a differentiable function $\tilde{\phi}: SE(d) \to \mathbb{R}$ is expressed as follows
\begin{equation}\label{ExteriorDerivative}
{\rm d}\tilde{\phi} = \sum \limits_{i=1}^{n_d}\mathcal{A}_{i}\tilde{\phi} \; \omega^{i} \in T^{*}(SE(d)).
\end{equation}
\begin{remark}\label{rem:1}
In differential geometry, there exist two equivalent viewpoints \cite[Ch. 2]{Aubin} on tangent vectors $\mathcal{A}_{g} \in T_{g}(SE(d))$: either one considers them as tangents to locally defined curves; or one considers them as differential operators on locally defined functions.
The connection between these viewpoints is as follows. We identify a
tangent vector $\dot{\tilde{\gamma}}(t) \in T_{\tilde{\gamma}(t)}(SE(d))$
with the differential operator $(\dot{\tilde{\gamma}}(t))(\tilde{\phi}) := \frac{d}{dt} \tilde{\phi}(\tilde{\gamma}(t))$
for all locally defined, differentiable, real-valued functions $\tilde{\phi}$.
\end{remark}
Next we express tangent vectors explicitly in the left-invariant moving frame of reference, by taking a directional derivative:
\begin{equation}\label{ddtcurve}
\!
\boxed{
\frac{d}{dt} \tilde{\phi}(\tilde{\gamma}(t)) =\langle {\rm d}\tilde{\phi}(\tilde{\gamma}(t)),
\dot{\tilde{\gamma}}(t) \rangle
  =
  \sum \limits_{i=1}^{n_d} \dot{\tilde{\gamma}}^{i}(t) \left.\mathcal{A}_{i}\right|_{\tilde{\gamma}(t)}\tilde{\phi}
}
\end{equation}
with $\dot{\tilde{\gamma}}(t)= \sum \limits_{i=1}^{n_d} \dot{\tilde{\gamma}}^{i}(t)\left.\mathcal{A}_{i}\right|_{\tilde{\gamma}(t)}$, and with $\tilde{\phi}$ smooth and defined on an open set around $\tilde{\gamma}(t)$. Eq.~(\ref{ddtcurve}) will play a crucial role in Section~\ref{ch:SE2} (exponential curve fits for $d=2$) and
Section~\ref{ch:SE3} (exponential curve fits for $d=3$).
\begin{example}
For $d=2$ we take $A_{1}=\left.\partial_{x}\right|_{\gothic{e}}$, $A_{2}=\left.\partial_y\right|_{\gothic{e}}$, $A_{3}=\left.\partial_{\theta}\right|_{\gothic{e}}$.
Then we have the left-invariant vector fields
\begin{equation} \label{VFs2}
\begin{array}{l}
\left.\mathcal{A}_{1}\right|_{(x,y,\theta)}:= \cos \theta \left.\frac{\partial}{\partial x}\right|_{(x,y,\theta)} + \sin \theta \left. \frac{\partial}{\partial y}\right|_{(x,y,\theta)}, \\
\left.\mathcal{A}_{2}\right|_{(x,y,\theta)}:= -\sin \theta \left. \frac{\partial}{\partial x}\right|_{(x,y,\theta)} + \cos \theta \left.\frac{\partial}{\partial y}\right|_{(x,y,\theta)}, \\
\left.\mathcal{A}_{3}\right|_{(x,y,\theta)}:= \left.\frac{\partial}{\partial \theta}\right|_{(x,y,\theta)}.
\end{array}
\end{equation}
The dual frame is given by
\begin{equation}\label{dual}
\begin{array}{l}
\omega^{1}=\cos \theta {\rm d}x + \sin \theta {\rm d}y, \\
\omega^{2}=-\sin \theta {\rm d}x
+\cos \theta {\rm d}y, \\
\omega^{3}={\rm d}\theta.
\end{array}
\end{equation}
For explicit formulas for left-invariant vector fields in $SE(3)$ see \cite{Chirikjian1,DuitsIJCV2010}.
\end{example}

\subsection{Exponential Curves in $SE(d)$}\label{ch:ExponentialCurves}

Let $(\ul{c}^{(1)},\ul{c}^{(2)})^T \in \R^{d+r_d}=\R^{n_d}$ be a given column vector, where $\ul{c}^{(1)}=(c^1,\ldots,c^{d}) \in \R^d$ denotes the spatial part and $\ul{c}^{(2)}=(c^{d+1},\ldots, c^{n_d}) \in \R^{r_d}$ denotes the rotational part. The unique exponential curve passing through $g \in SE(d)$
with initial velocity $\ul{c}(g)=\sum \limits_{i=1}^{n_d} c^{i}\left.\mathcal{A}_{i}\right|_{g}$ equals
\begin{equation} \label{expcurvend}
\boxed{
\tilde{\gamma}^{\ul{c}}_{g}(t)= g \; e^{t \sum \limits_{i=1}^{n_d} c^{i}A_{i}}
}
\end{equation}
with $A_{i}=\left.\mathcal{A}_{i}\right|_{\gothic{e}}$ denoting a basis of $T_{\gothic{e}}(SE(d))$.
%
In fact such exponential curves satisfy
\begin{equation}
\label{ec}
\dot{\tilde{\gamma}}(t)= \sum \limits_{i=1}^{n_d} c^{i} \left.\mathcal{A}_{i}\right|_{\tilde{\gamma}(t)}
\end{equation}
and thereby have constant velocity in the moving frame of reference, i.e. $\dot{\tilde{\gamma}}^{i}=c^{i}$ in Eq.~\!(\ref{ddtcurve}). A way to compute the exponentials is via matrix exponentials and (\ref{thedw}).
\begin{example}
\mbox{If $d=2$ we have exponential curves:}
\begin{equation} \label{horexp} \! \!
\begin{array}{rl}
 \tilde{\gamma}_{g_0}^{\ul{c}}(t)&= g_0 \, e^{t(c^{1}\! A_{1}+c^{2}\! A_{2} +c^{3} \! A_{3})} = (x(t),y(t),\theta(t)),
\end{array}
\end{equation}
which are circular spirals with
\begin{equation}  \label{horexp1}
  \begin{array}{rl}
 x(t)& = x_0 + \frac{c^{1}}{c^{3}}(\sin(c^{3}t \!+\!\theta_0)-\sin(\theta_{0})) \\
  &\,\,\, +\frac{c^{2}}{c^{3}}(\cos(c^{3}t \!+\!\theta_0)-\cos(\theta_{0}))\ , \\
 y(t)& = y_0-\frac{c^{1}}{c^{3}}(\cos(c^{3}t \!+\!\theta_0)-\cos(\theta_{0}))  \\
  &\,\,\,+ \frac{c^{2}}{c^{3}}(\sin(c^{3}t \!+\!\theta_0)-\sin(\theta_{0}))\ , \\
 \theta(t) &= \theta_{0} + t c^{3} ,
 \end{array}
\end{equation}
for the case $c^{3}\neq 0$, and all $t\geq 0$ and straight lines with
\begin{equation}  \label{horexp2}
  \begin{array}{rl}
 x(t) &= x_0+ t (c^{1}\cos\theta_0 - c^2 \sin \theta_0), \\
 y(t) &= y_0+t (c^{1} \sin \theta_{0}+c^{2} \cos \theta_0), \\
 \theta(t) &= \theta_{0},
 \end{array}
\end{equation}
for the case $c^{3}=0$, where $g_0=(x_0,y_0,\theta_{0}) \in SE(2)$.
See the left panel in Fig.~\!\ref{fig:geometry}.
\end{example}
\begin{example}
For $d=3$, the formulae for exponential curves in $SE(3)$ are given in for example \cite{Chirikjian1,DuitsIJCV2010}. Their spatial part are circular spirals
with torsion $\TTT(t)=\frac{|\ul{c}^{(1)} \cdot \ul{c}^{(2)}|}{\|\ul{c}^{(1)}\|}\, \KKK(t)$ and curvature \begin{equation} \label{curvature}
\begin{array}{ll}
\KKK(t) &= \frac{1}{\|\ul{c}^{(1)}\|^2}\big( \cos(t \|\ul{c}^{(2)}\|)\, \ul{c}^{(2)}\times \ul{c}^{(1)}  \\
 &+ \frac{\sin(t \|\ul{c}^{(2)}\|)}{\|\ul{c}^{(2)}\|} \ul{c}^{(2)} \times \ul{c}^{(2)} \times \ul{c}^{(1)}\big).
\end{array}
\end{equation}
Note that their magnitudes are constant:
\begin{equation} \label{torsion1}
|\kappa|= \frac{\|\ul{c}^{(1)} \times \ul{c}^{(2)}\|}{\|\ul{c}^{(1)}\|^2}\textrm{ and }
|\tau|= \frac{|\ul{c}^{(1)} \cdot \ul{c}^{(2)}|\cdot |\kappa|}{\|\ul{c}^{(1)}\|}.
\end{equation}
\end{example}

\subsection{Left-Invariant Metric Tensor on $SE(d)$}\label{ch:Metric}

We use the following (left-invariant) metric tensor:
\begin{equation} \label{metrictensor2}
\left.\gothic{G}_{\mu}\right|_{\tilde{\gamma}}(\dot{\tilde{\gamma}},\dot{\tilde{\gamma}})=
\mu^2 \sum \limits_{i=1}^{d} |\dot{\tilde{\gamma}}^{i}|^2
+
\sum \limits_{i=d+1}^{n_d} |\dot{\tilde{\gamma}}^{i}|^2,
\end{equation}
where $\dot{\tilde{\gamma}}=\sum_{i=1}^{n_d}\dot{\tilde{\gamma}}^{i} \left.\mathcal{A}_{i}\right|_{\tilde{\gamma}}$, and with stiffness parameter $\mu$ along any smooth curve $\tilde{\gamma}$ in $SE(d)$.
Now, for the special case of exponential curves, one has $\dot{\tilde{\gamma}}^{i}= c^{i}$ is constant. The metric allows us to normalize the speed along the curves by imposing
a normalization constraint
\begin{equation} \label{normmu}
\begin{array}{rl}
\|\ul{c}\|_{\mu}^2 :=\|\ul{M}_{\mu}\ul{c}\|^2&=
\mu^{2}\sum \limits_{i=1}^{d}|c^{i}|^{2} + \sum \limits_{i=d+1}^{n_d}|c^{i}|^{2}\\
&=\mu^{2} \|\ul{c}^{(1)}\|^{2} +  \|\ul{c}^{(2)}\|^{2}=1, \\[6pt]
   \multicolumn{2}{l}{\textrm{with } \ul{M}_{\mu}:= \begin{pmatrix}
   \mu I_d  && \ul{0} \\
   \ul{0} && I_{r_d} \end{pmatrix} 
  \in \R^{n_d\times n_d}.}
\end{array}
\end{equation}
We will use this constraint in the fitting procedure in order to ensure that our exponential curves (\ref{ec}) are parameterized by Riemannian arclength $t$.

\subsection{Convolution and Haar-measure on $SE(d)$}
In general a convolution of data $\tilde{U}:SE(d)\to \R$ with kernel $\tilde{K}:SE(d)\to \R$ is given by
{\small
\begin{equation} \label{Haarmeasure}
\begin{array}{l}
(\tilde{K} * \tilde{U})(g) = \int \limits_{SE(d)} \tilde{K}(h^{-1}g)\, \tilde{U}(h)\, {\rm d}\overline{\mu}(h)= \\
\int \limits_{\R^d}\!\int \limits_{SO(d)} \tilde{K}((\ul{R}')^{-1}\!(\ul{x}\!-\!\ul{x}'),(\ul{R}')^{-1}\!\ul{R})\, {\rm d}\ul{x}' {\rm d}\mu_{SO(d)}(\ul{R}'), \\
\textrm{with } {\rm d}\overline{\mu}(h)={\rm d}\ul{x}' {\rm d}\mu_{SO(d)}(\ul{R}'),
\end{array}
\end{equation}}
for all $h=(\ul{x}',\ul{R}') \in SE(d)$, where Haar measure $\overline{\mu}$ is the direct product of the usual Lebesgue measure on $\R^d$ with the Haar measure on $SO(d)$.

\subsection{Gaussian Smoothing and Gradient on $SE(d)$}\label{ch:GaussianSmoothing}

We define the regularized data
\begin{equation} \label{VEET}
\tilde{V} := \tilde{G}_{\ul{s}} * \tilde{U},
\end{equation}
where $\ul{s}=(s_p,s_o)$ are the spatial and angular scales respectively of the separable Gaussian smoothing kernel defined by
\begin{equation} \label{Gtilde}
\tilde{G}_{\ul{s}}(\ul{x},\ul{R}):= G_{s_p}^{\R^d}(\ul{x}) \; G_{s_o}^{S^{d-1}}(\ul{R}\ul{a}).
\end{equation}
This smoothing kernel is a product of the heat kernel {\small $G_{s_p}^{\R^d}(\ul{x})=\frac{e^{-\frac{\|\ul{x}\|^2}{4s_p}}}{(4\pi s_p)^{d/2}}$} on $\R^{d}$ centered at $\ul{0}$ with spatial scale $s_p>0$,
and a heat kernel $G_{s_o}^{S^{d-1}}(\ul{R}\ul{a})$ on $S^{d-1}$ centered around $\ul{a} \in S^{d-1}$ with angular scale $s_o>0$.

By definition the gradient $\nabla \tilde{U}$ of image data $\tilde{U}:SE(d) \to \R$ is the Riesz representation vector of the derivative ${\rm d}\tilde{U}$:
\begin{equation} \label{gradientmu2}
\begin{array}{rl}
\nabla \tilde{U} :=& \gothic{G}_{\mu}^{-1}{\rm d}\tilde{U} = \sum \limits_{i=1}^{d}\mu^{-2} \mathcal{A}_{i}\tilde{U} \, \mathcal{A}_{i}
\! +  \!\sum \limits_{j=d+1}^{n_d}\! \mathcal{A}_{j}U \, \mathcal{A}_{j} \\
  \equiv & \ul{M}_{\mu^{-2}} (\mathcal{A}_{1}\tilde{U},\ldots, \mathcal{A}_{n_d}\tilde{U})^T,
\end{array}
\end{equation}
relying on $\ul{M}_{\mu}$ as defined in (\ref{normmu}). Here, following standard conventions in differential geometry, $\gothic{G}_{\mu}^{-1}$ denotes the inverse of the linear map associated to
the metric tensor (\ref{metrictensor2}). Then, the \emph{Gaussian} gradient is defined by
\begin{equation} \label{GGR2}
\nabla^{\ul{s}} \tilde{U} :=\nabla \tilde{V}=\nabla (\tilde{G}_{\ul{s}} * \tilde{U})=\nabla \tilde{G}_{\ul{s}}* \tilde{U}.
\end{equation}

\subsection{Horizontal Exponential Curves in $SE(d)$}
Typically, in the distribution $\tilde{U}$ (e.g. if $\tilde{U}$ is an orientation score of a grey-scale image) the mass is concentrated around so-called horizontal exponential curves in $SE(d)$ (see Fig. \ref{Fig:Intro}). Next we explain this notion of horizontal exponential curves.

A curve $t \mapsto (x(t),y(t)) \in \R^2$  can be lifted to a curve $t \mapsto \tilde{\gamma}(t)=(x(t),y(t),\theta(t))$ in $SE(2)$ via
 \begin{equation} \label{lift}
 \theta(t)=\arg\{\dot{x}(t)+ i \, \dot{y}(t)\}.
 \end{equation}
Generalizing to $d\geq 2$, one can lift a curve $t \mapsto \ul{x}(t) \in \R^d$ towards a curve $t \mapsto \gamma(t)=(\ul{x}(t),\ul{n}(t))$ in $\R^{d}\rtimes S^{d-1}$ by setting
\[\ul{n}(t)=\|\dot{\ul{x}}(t)\|^{-1}\dot{\ul{x}}(t).\]
A curve $t \mapsto \ul{x}(t)$ can be lifted towards a family of lifted curves $t \mapsto \tilde{\gamma}(t)=(\ul{x}(t),\ul{R}_{\ul{n}(t)})$ into the roto-translation group $SE(d)$ by setting $\ul{R}_{\ul{n}(t)} \in SO(d)$ such that it maps reference axis $\ul{a}$ onto $\ul{n}(t)$:
\begin{equation} \label{lifttogroup}
\ul{R}_{\ul{n}(t)}\ul{a}=\ul{n}(t)=\|\dot{\ul{x}}(t)\|^{-1}\dot{\ul{x}}(t). 
\end{equation}
Here we use $\ul{R}_{\ul{n}}$ to denote any rotation that maps reference axis $\ul{a}$ onto $\ul{n}$.
Clearly, the choice of rotation is not unique for $d>2$, e.g.~if $d=3$ then $\ul{R}_{\ul{n}} \ul{R}_{\ul{a},\alpha}\ul{a}=\ul{a}$ regardless the value of $\alpha$, where $\ul{R}_{\ul{a},\alpha}$ denotes the counter-clockwise 3D rotation about axis $\mathbf{a}$ by angle $\alpha$.

Next we study the implication of restriction (\ref{lifttogroup}) on the tangent bundle of $SE(d)$.
\begin{itemize}
\item
For $d=2$, we have restriction $\dot{\ul{x}}(t)=(\dot{x}(t), \dot{y}(t))=\|\dot{\ul{x}}(t)\|(\cos \theta(t),\sin \theta(t))$, i.e.
\begin{equation}\label{distributionSE2}
\begin{array}{ll}
\dot{\tilde{\gamma}} \in \left.\Delta\right|_{\tilde{\gamma}}, \
\textrm{with }\Delta&=\textrm{span}\{\cos \theta \partial_x\!+\!\sin \theta \partial_y, \partial_{\theta}\}\!\!\!\\
&=\textrm{span}\{\mathcal{A}_1,\mathcal{A}_3\},
\end{array}
\end{equation}
where $\Delta$ denotes
the so-called horizontal part of tangent bundle $T(SE(2))$. See Fig.~\!\ref{fig:geometry}.
\begin{figure}
\includegraphics[width=1.05\hsize]{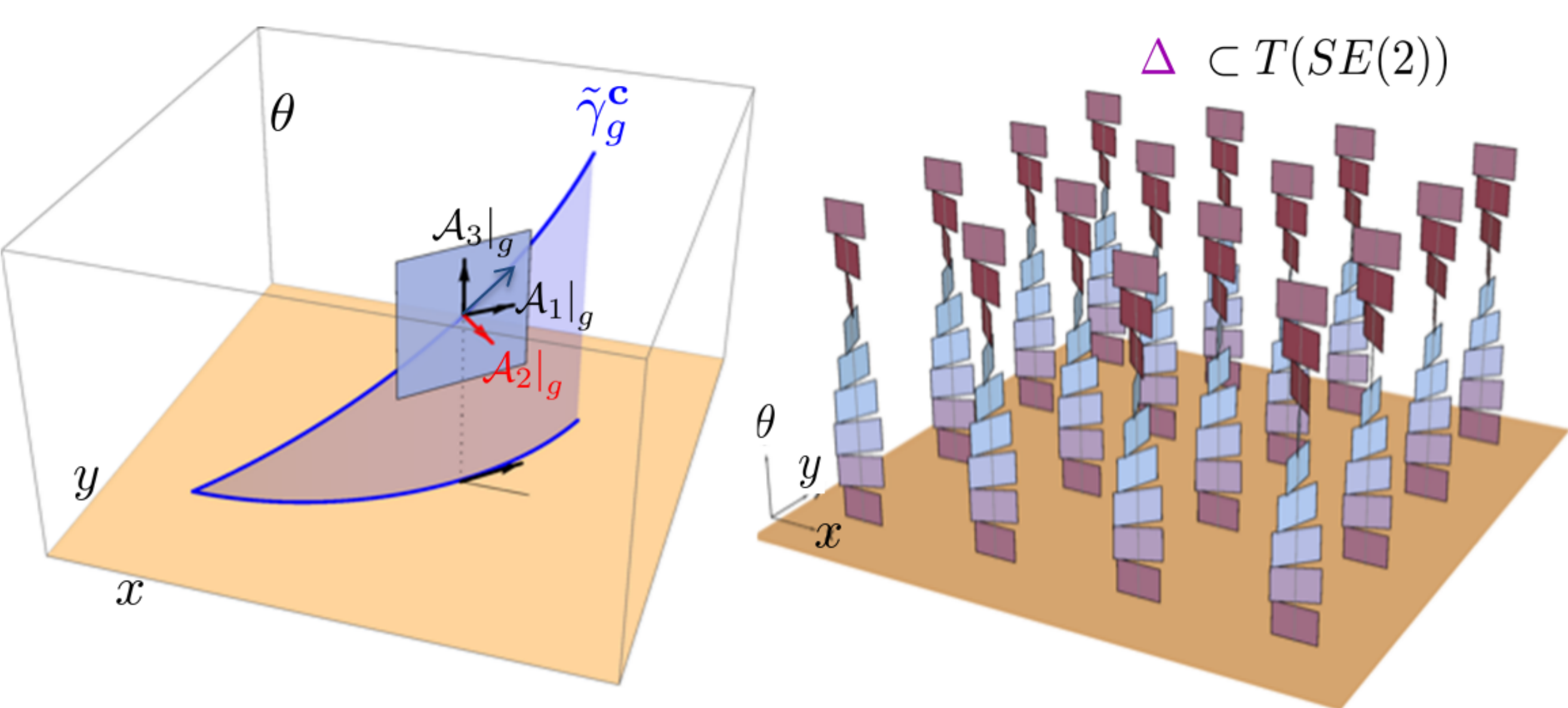}
\caption{Left: horizontal exponential curve $\tilde{\gamma}^{\ul{c}}_g$ in $SE(2)$ with $\ul{c}=(1,0,1)$. Its projection on the ground plane reflects co-circularity, and the curve can be obtained by a lift (\ref{lift}) from its spatial projection. Right: the distribution $\Delta$ of horizontal tangent vector fields as a sub-bundle in the tangent bundle $T(SE(2))$. \label{fig:geometry}}
\end{figure}
\item
For $d=3$, we impose the constraint:
\begin{equation} \label{distributionSE3}
\begin{array}{l}
\dot{\tilde{\gamma}}(t) \in \Delta_{\tilde{\gamma}(t)}, \
\textrm{with }\Delta:= \textrm{span}\{\mathcal{A}_{3},\mathcal{A}_{4}, \mathcal{A}_{5}\},
\end{array}
\end{equation}
where $\mathcal{A}_3= \ul{n} \cdot \nabla_{\R^3}$, since then spatial transport is always along $\ul{n}$ which is required for for (\ref{lifttogroup}).

\end{itemize}
Curves $\tilde{\gamma}(t)$ satisfying the constraint (\ref{distributionSE2}) for $d=2$, and (\ref{distributionSE3}) for $d=3$  
are called \emph{horizontal} curves. Note that $\textrm{dim}(\Delta)=d$.


Next we study how the restriction applies to the particular case of exponential curves on $SE(d)$.
\begin{itemize}
\item
For $d=2$ horizontal exponential curves are obtained from (\ref{horexp}), (\ref{horexp1}), (\ref{horexp2}) by setting $c^{2}=0$.
\item
For $d=3$, we use a different reference axis $\ul{a}$, and horizontal exponential curves are obtained from (\ref{expcurvend}) by setting $c^1=c^2=c^6=0$.
\end{itemize}
If exponential curves are not horizontal, then we indicate how much the local tangent of the exponential curve points outside the spatial part of $\Delta$, by a `deviation from horizontality angle' $\chi$, which is given by:
\begin{equation} \label{dH}
\chi =\arccos \left( \left|\frac{\ul{c}^{(1)}\cdot \ul{a}}{\|\ul{c}^{(1)}\|} \right| \right).
\end{equation}
\begin{example}
In case $d=2$ we have $n_2=3$, $\ul{a}=(1,0)^T$. The horizontal part of the tangent bundle $\Delta$ is given by (\ref{distributionSE2}), and horizontal exponential curves are obtained from (\ref{horexp}) by setting $c^2=0$.
For exponential curves in general, we have  deviation from horizontality angle
\begin{equation}\label{dH2}
\chi=\arccos \left( \left|\frac{c^{1}}{\sqrt{|c^{1}|^2+|c^{2}|^2}}\right|\right).
\end{equation}
An exponential curve in $SE(2)$ is horizontal if and only if $\chi=0$.
See Fig.~\ref{fig:geometry}, where in the left we have depicted a horizontal exponential curve and where in the right we have visualized distribution $\Delta$.
\end{example}
\begin{example}
In case $d=3$, we have $n_3=6$, $\ul{a}=(0,0,1)^T$. The horizontal part of the tangent bundle is given by (\ref{distributionSE3}), and horizontal exponential curves are characterized by $c^{3}, c^{4}, c^{5}$ whereas $c^{1}=c^{2}=c^6=0$. By Eq.~\!(\ref{torsion1}) these curves have zero torsion $|\tau|=0$ and constant curvature  $\frac{\sqrt{(c^4)^2+(c^{5})^2}}{c^3}$ and thus they are planar circles. For exponential curves in general, we have deviation from horizontality angle
\[
\chi = \arccos \left( \left|\frac{c^3}{\sqrt{|c^1|^2+|c^2|^2+|c^3|^2}} \right| \right).
\]
An exponential curve in $SE(3)$ is horizontal if and only if $\chi=0$ and $c^6=0$.
\end{example}

\section{From Exponential Curve Fits to Gauge Frames on $SE(d)$\label{ch:togauge}}

In Section~\ref{ch:SE2} and Section~\ref{ch:SE3} we will discuss techniques to find an exponential curve $\tilde{\gamma}^{\ul{c}}_g(t)$ that fits the data $\tilde{U}:SE(d)\to \R$ locally. 
Let $\ul{c}(g)=(\tilde{\gamma}^{\ul{c}}_g)'(0)$ be its tangent vector at $g$.

In this section we assume that the tangent vector $\ul{c}(g)=(\ul{c}^{(1)}(g),\ul{c}^{(2)}(g))^T \in R^{d+r_d}=\R^{n_d}$ is given. From this vector we will construct a locally adaptive frame $
\{\left.\mathcal{B}_{1}\right|_g,\ldots,\left.\mathcal{B}_{n_d}\right|_{g}\}$, orthonormal w.r.t. $\gothic{G}_{\mu}$-metric in such a way that:
\begin{enumerate}
\item the main spatial generator ($\mathcal{A}_{1}$ for $d=2$ and $\mathcal{A}_{d}$ for $d>2$) is mapped onto $\left.\mathcal{B}_{1}\right|_{g}= \sum \limits_{i=1}^{n_d}c^{i}(g)
\left.\mathcal{A}_{i}\right|_{g}$,
\item the spatial generators $\{\left.\mathcal{B}_{i}\right|_{g}\}^d_{i=2}$ are obtained from the other left-invariant spatial generators $\{\left.\mathcal{A}_{i}\right|_g\}^d_{i=1}$ by a planar rotation of $\ul{a}$ onto $\frac{\ul{c}^{(1)}}{\| \ul{c}^{(1)}\|}$ by angle $\chi$. In particular, if $\chi=0$, the other spatial generators do not change their direction. This allows us to still distinguish spatial generators and angular generators in our adapted frame.
\end{enumerate}
Next we provide for each $g \in SE(d)$ the explicit construction of a rotation matrix $\ul{R}^{\ul{c}(g)}$ and a scaling by $\ul{M}_{\mu^{-1}}$ on $T_g(SE(d))$, which maps frame
$\{\left.\mathcal{A}_{1}\right|_g,\ldots,\left.\mathcal{A}_{n_d}\right|_{g}\}$ onto
$\{\left.\mathcal{B}_{1}\right|_g,\ldots,\left.\mathcal{B}_{n_d}\right|_{g}\}$.

The construction for $d>2$ is technical and provided in Theorem~\ref{app:a} in Appendix~\ref{app:a}. However,
the whole construction of the rotation matrix $\ul{R}^{\ul{c}}$ via a concatenation of two subsequent rotations 
is similar to the case $d=2$ that we will explain next.

\begin{figure}[b]
\includegraphics[width=\hsize]{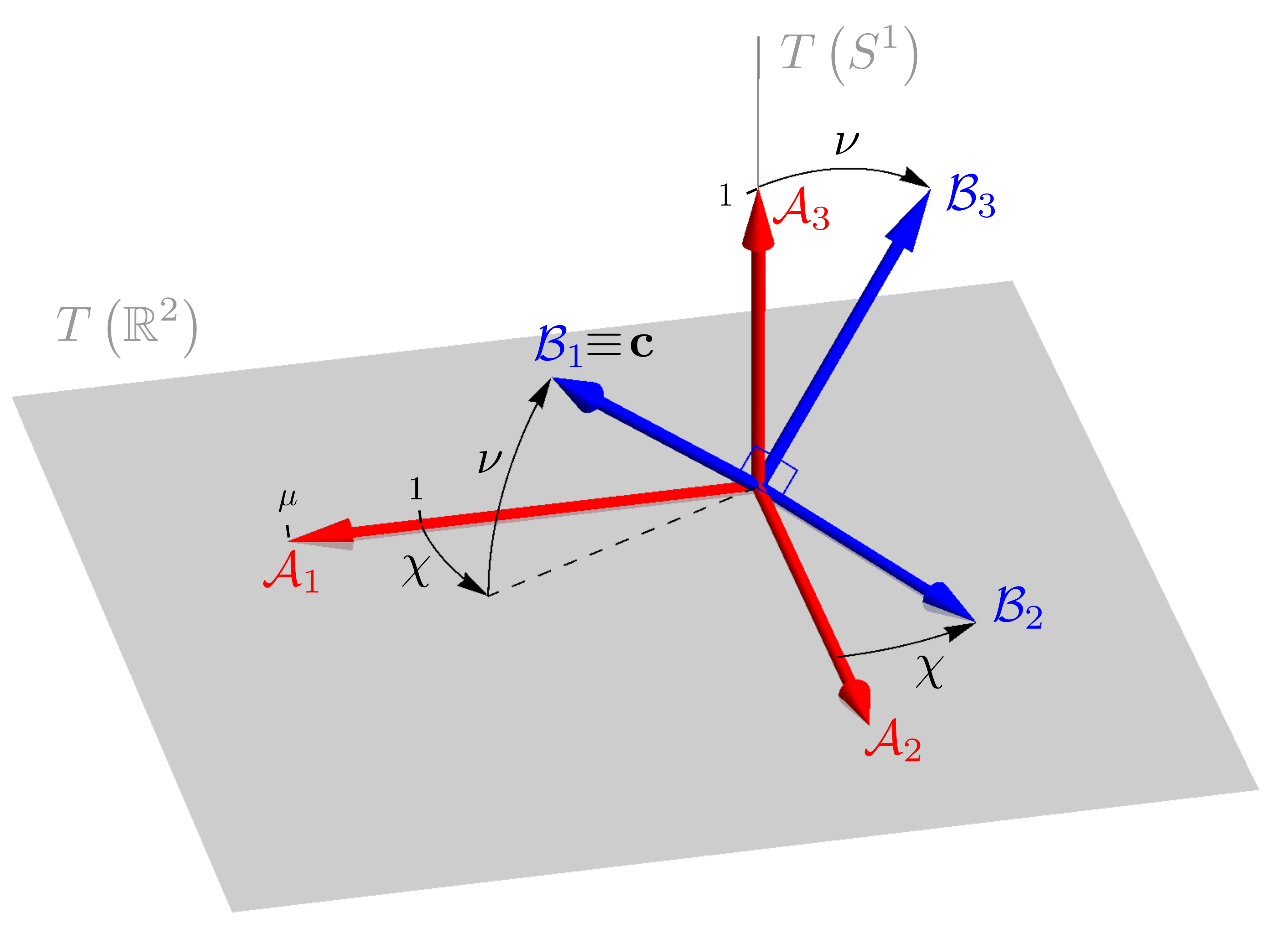}
\caption{Locally adaptive frame $\{\left.\mathcal{B}_{1}\right|_{g}, \left.\mathcal{B}_{2}\right|_{g}, \left.\mathcal{B}_{3}\right|_{g}\}$
(in blue) in $T_{g}(SE(2))$ (with $g$ placed at the origin) is obtained from frame $\{\left.\mathcal{A}_{1}\right|_{g}, \left.\mathcal{A}_{2}\right|_{g}, \left.\mathcal{A}_{3}\right|_{g}\}$ (in red) and $\ul{c}(g)$, via normalization and
two subsequent rotations $\ul{R}^{\ul{c}}=\ul{R}_{2}\ul{R}_{1}$, see Eq.~\!(\ref{gauge}), revealing deviation from horizontality $\chi$ in $R_{1}$, spherical angle $\nu$ in Eq.\!~(\ref{R1R2}). Vector field $\mathcal{A}_{1}$ takes a spatial derivative in direction $\ul{n}$, whereas $\mathcal{B}_{1}$ takes a derivative along the tangent $\ul{c}$ of the local exponential curve fit.
  \label{fig:gauge}}
\end{figure}

Consider $d=2$ where the frames $\{\mathcal{A}_{1},\mathcal{A}_{2},\mathcal{A}_{3}\}$ and $\{\mathcal{B}_{1},\mathcal{B}_{2},\mathcal{B}_{3}\}$ are depicted in Fig.~\!\ref{fig:gauge}
The explicit relation between the normalized gauge frame and the left-invariant vector field frame is given by
\begin{equation} \label{gauge}
\underline{\mathcal{B}}:= (\ul{R}^{\ul{c}})^T \ul{M}_{\mu}^{-1}\underline{\mathcal{A}},
\end{equation}
with $\underline{\mathcal{A}}:=(\mathcal{A}_1,\mathcal{A}_{2},\mathcal{A}_{3})^T$,
$\underline{\mathcal{B}}:=(\mathcal{B}_1,\mathcal{B}_{2},\mathcal{B}_{3})^T$,
and with rotation matrix
{\small
\begin{equation} \label{R1R2}
\begin{array}{l}
\ul{R}^{\ul{c}}=\ul{R}_{2}\ul{R}_{1} \in SO(3),
\textrm{ with } \\
\ul{R}_{2}= \left(
\begin{array}{ccc}
 \cos \chi &\! -\sin \chi & 0 \\
 \sin \chi & \cos \chi & 0\\
0 &  0 & 1
\end{array}
\right)\!,
\ul{R}_{1}=
\left(
\begin{array}{ccc}
 \cos \nu& 0 & \sin \nu \\
0 & 1 & 0 \\
\!-\sin \nu &  0 & \cos \nu
\end{array}
\right),
\end{array}
\end{equation}
}
where the rotation angles are the deviation from horizontality angle $\chi$ and
the spherical angle \\
\[
\nu= \arcsin \left(\frac{c^3}{\|\ul{c}\|_{\mu}} \right) 
\in [-\pi/2,\pi/2].
\]
Recall that $\chi$ is given by (\ref{dH2}).
The multiplication $\ul{M}_{\mu}^{-1}\underline{\mathcal{A}}$ ensures that each of the vector fields in the locally adaptive frame is normalized w.r.t. the $\gothic{G}_{\mu}$-metric, recall (\ref{metrictensor2}).

\begin{remark}
When imposing isotropy (w.r.t. the metric $\gothic{G}_{\mu}$) in the plane orthogonal to $\mathcal{B}_{1}$, 
there is a unique choice $\ul{R}^{\ul{c}}$ mapping $(1,0,0)^T$ onto $(\mu c^{1},\mu c^{2},c^{3})^T$ such that it keeps the other spatial generator in the spatial subspace of $T_g(SE(2))$ (and with $\chi=0 \desda
\mathcal{B}_{2}=\mu^{-1}\mathcal{A}_{2}$). This choice is given by (\ref{R1R2}).

\end{remark}
The generalization to the $d$-dimensional case of the construction of a locally adaptive frame $\{\mathcal{B}_{i}\}_{i=1}^{n_d}$ from
$\{\mathcal{A}_{i}\}_{i=1}^{n_d}$ and the tangent vector $\ul{c}$ of a given exponential curve fit $\tilde{\gamma}^{\ul{c}}_g(\cdot)$ to data
$\tilde{U}:SE(d)\to \mathbb{R}$ is explained in Theorem~\ref{app:th} in Appendix~\ref{app:a}.

\section{Exponential Curve Fits in $\R^{d}$ \label{ch:3-1} \label{ch:gauge}}

In this section we reformulate the classical construction of a locally adaptive frame to image $f$ at location $\ul{x} \in \R^d$, in a group-theoretical way.
This reformulation seems technical at first sight, but helps in understanding the formulation of projected exponential curve fits in the higher dimensional Lie group $SE(d)$.

\subsection{Exponential Curve Fits in $\R^{d}$ of the 1st Order}
We will take the structure tensor approach \cite{Bigun1987,Knustsson1989}, which will be shown to yield first-order exponential curve fits.

The Gaussian gradient
\begin{equation} \label{GGR22}
\nabla^{s} f =\nabla G_{s} * f,
\end{equation}
with Gaussian kernel
\begin{equation}\label{GaussianKernelRd}
  G_{s}(\ul{x})= (4\pi s)^{-d/2}e^{-\frac{\|\ul{x}\|^2}{4s}},
\end{equation}
is used in the definition of the structure matrix:
\begin{equation}
\mathbf{S}^{s,\rho}(f) = G_{\rho} * \nabla^{s} f \, (\nabla^{s} f)^T,
\end{equation}
with $s=\frac{1}{2}\sigma^2_{s}$, and $\rho=\frac{1}{2}\sigma^2_{\rho}$ the scale of regularization typically yielding a non-degenerate and positive definite matrix. In the remainder we use short notation $\mathbf{S}^{s,\rho}:=\mathbf{S}^{s,\rho}(f)$. The structure matrix appears in solving the following optimization problem where for all $\ul{x} \in \R^{d}$ we aim to find optimal tangent vector
\begin{equation} \label{optR2a}
\!
\begin{array}{ll}
\ul{c}^{*}(\ul{x})  &= \argmin_{\scriptsize \begin{array}{c}\ul{c} \in \R^{d},\\
\|\ul{c}\|=1
\end{array}} \; \int \limits_{\R^{d}}\! G_{\rho}(\ul{x}\!-\!\ul{x}')
|\nabla^{s}f(\ul{x}') \cdot \ul{c}|^2{\rm d}\ul{x}'\! \\
 &=\argmin_{\scriptsize \begin{array}{c}\ul{c} \in \R^{d},\\
\|\ul{c}\|=1
\end{array}} \; \ul{c}^{T} \mathbf{S}^{s,\rho}(\ul{x}) \ul{c}.
\end{array}
\end{equation}
In this optimization problem we find the tangent $\ul{c}^{*}(\ul{x})$ which minimizes a (Gaussian) weighted average of the squared directional derivative $|\nabla^{s}f(\ul{x}') \cdot \ul{c}|^2$ in the neighborhood of $\ul{x}$. The second identity in (\ref{optR2a}), which directly follows from the definition of the structure matrix, allows us to solve optimization problem (\ref{optR2a}) via the Euler-Lagrange equation
\begin{equation}
\mathbf{S}^{s,\rho}(\ul{x})\; \ul{c}^{*}(\ul{x}) = \lambda_{1} \ul{c}^{*}(\ul{x}),
\end{equation}
since the minimizer is found as the eigenvector $\ul{c}^{*}(\ul{x})$ with the smallest eigenvalue $\lambda_1$.

Now let us put Eq.\!~(\ref{optR2a}) in group-theoretical form by reformulating it as an exponential curve fitting problem. This is helpful in our subsequent generalizations to $SE(d)$.
On $\R^{d}$ exponential curves are straight lines:
\begin{equation} \label{expcurvesRd}
\gamma_{\ul{x}}^{\ul{c}}(t) =\ul{x} + \exp_{\R^{d}}(t \ul{c}) =\ul{x} + t \ul{c},
\end{equation}
and on $T(\R^{d})$ we impose the standard flat metric tensor
$\gothic{G}(\ul{c},\ul{d})= \sum_{i=1}^{d} c^{i} d^{i}$.
In (\ref{optR2a}) we replace the directional derivative by a time derivative (at $t=0$) when moving over an exponential curve:
\begin{equation} \label{optR2}
\boxed{
\begin{array}{l}
\ul{c}^{*}(\ul{x})  = \argmin_{\ul{c} \in \R^{d}, \|\ul{c}\|=1} \\[8pt]
\int \limits_{\R^{d}} G_{\rho}(\ul{x}-\ul{x}')
\left|\; \left.\frac{d}{dt} (G_{s}*f)(\gamma_{\ul{x}',\ul{x}}^{\ul{c}}(t))\, \right|_{t=0}\; \right|^2
 \,{\rm d}\ul{x}',  \\[8pt]
 \end{array}
 }
\end{equation}
where
\begin{equation} \label{paralleltransportRd}
t \mapsto \gamma_{\ul{x}',\ul{x}}^{\ul{c}}(t)= \gamma_{\ul{x}}^{\ul{c}}(t) -\ul{x} + \ul{x}'= \gamma_{\ul{x}'}^{\ul{c}}(t).
\end{equation}
Because in (\ref{optR2a}) we average over directional derivatives in the neighborhood of $\ul{x}$ we now average the time derivatives over a \emph{family of neighboring exponential curves} $\gamma_{\ul{x}',\ul{x}}^{\ul{c}}(t)$, which are defined to start at neighboring positions $\ul{x}'$ but having the same spatial velocity as $\gamma_{\ul{x}}^{\ul{c}}(t)$. In $\R^{d}$ the distinction between $\gamma_{\ul{x}',\ul{x}}^{\ul{c}}(t)$ and $\gamma_{\ul{x}'}^{\ul{c}}(t)$ is not important but it will be in the $SE(d)$-case.
\begin{definition} \label{def:1}
Let $\ul{c}^{*}(\ul{x}) \in T_{\ul{x}}(\R^{d})$ be the minimizer in (\ref{optR2}).
We say $\gamma_{\ul{x}}(t)= \ul{x}+\exp_{\R^{d}}({t \ul{c}^{*}(\ul{x}))}$ is the first-order exponential curve fit to image data
$f: \R^{d} \to \R$ at location $\ul{x}$.
\end{definition}
\subsection{Exponential Curve Fits in $\R^{d}$ of the 2nd Order}
For second-order exponential curve fits we need the Hessian matrix defined by
\begin{equation}\label{hessianRd}
  (\ul{H}^{s}(f))(\ul{x})= \left[\partial_{x_j}\partial_{x_i} (G_{s} * f)(\ul{x})\right] ,
\end{equation}
with $G_s$ the Gaussian kernel given in Eq.~\!(\ref{GaussianKernelRd}). From now on we use short notation $\ul{H}^{s}:=\ul{H}^{s}(f)$. When using the Hessian matrix for curve fitting we aim to solve
\begin{equation} \label{optR2Hessian}
\!
\begin{array}{ll}
 \ul{c}^{*}(\ul{x})  &= \argmin_{\scriptsize \begin{array}{c}\ul{c} \in \R^{d},\\
\|\ul{c}\|=1
\end{array}} \; | \ul{c}^{T} \mathbf{H}^{s}(\ul{x}) \ul{c} |.
\end{array}
\end{equation}
In this optimization problem we find the tangent $\ul{c}^{*}(\ul{x})$ which minimizes the second-order directional derivative of (Gaussian) regularized data $G_{s} * f$. When all Hessian eigenvalues have the same sign we can solve the optimization problem (\ref{optR2Hessian}) via the Euler-Lagrange equation
\begin{equation}
\mathbf{H}^{s}(\ul{x})\; \ul{c}^{*}(\ul{x}) = \lambda_{1} \ul{c}^{*}(\ul{x}),
\end{equation}
and the minimizer is found as the eigenvector $\ul{c}^{*}(\ul{x})$ with the smallest eigenvalue $\lambda_1$.

Now, we can again put Eq.\!~(\ref{optR2Hessian}) in group-theoretical form by reformulating it as an exponential curve fitting problem. This is helpful in our subsequent generalizations to $SE(d)$.
We again rely on exponential curves as defined in (\ref{expcurvesRd}). In (\ref{optR2Hessian}) we replace the second order directional derivative by a second order time derivative (at $t=0$) when moving over an exponential curve:
\begin{equation} \label{optR2HessianGroup}
\boxed{
\begin{array}{l}
\ul{c}^{*}(\ul{x})  = \!\!\! \argmin_{\ul{c} \in \R^{d}, \|\ul{c}\|=1} \left|\; \left.\frac{d^2}{dt^2} (G_{s}*f)(\gamma_{\ul{x}}^{\ul{c}}(t))\, \right|_{t=0}\; \right|.
 \end{array}
 }
\end{equation}

\begin{remark}
In general the eigenvalues of Hessian matrix $\ul{H}^s$ do not have the same sign. In this case we  still take $\ul{c}^{*}(g)$ as the eigenvector with smallest absolute eigenvalue (representing minimal absolute principal curvature), though this no longer solves (\ref{optR2Hessian}).
\end{remark}

\begin{definition} \label{def:SecondOrderFitRd}
Let $\ul{c}^{*}(\ul{x}) \in T_{\ul{x}}(\R^{d})$ be the minimizer in (\ref{optR2HessianGroup}).
We say $\gamma_{\ul{x}}(t)= \ul{x}+\exp_{\R^{d}}({t \ul{c}^{*}(\ul{x}))}$ is the second-order exponential curve fit to image data
$f: \R^{d} \to \R$ at location $\ul{x}$.
\end{definition}

\begin{remark}
  In order to connect optimization problem (\ref{optR2HessianGroup}) with the first order optimization (\ref{optR2}) we note that (\ref{optR2HessianGroup}) can also be written as an optimization over a family of curves $\gamma_{\ul{x}',\ul{x}}^{\ul{c}}$ defined in (\ref{paralleltransportRd}):
  \begin{equation} \label{optR2HessianGroupFamily}
\begin{array}{l}
\ul{c}^{*}(\ul{x})  =
\argmin_{\scriptsize \begin{array}{c}\ul{c} \in \R^{d},\\\|\ul{c}\|=1\end{array}}
\left|\; \int \limits_{\R^{d}} \!\! G_{s}(\ul{x}-\ul{x}') \! \left.\frac{d^2}{dt^2} (f)(\gamma_{\ul{x}',\ul{x}}^{\ul{c}}(t)) \right|_{t=0} \right|
 \,{\rm d}\ul{x}',
 \end{array}
\end{equation}
because of linearity of the second-order time derivative.
\end{remark}

\section{Exponential Curve Fits in $SE(2)$ \label{ch:SE2}}
As mentioned in the introduction we distinguish between two approaches:
a first order optimization approach based on a structure tensor on $SE(2)$, and a second order
optimization approach based on the Hessian on $SE(2)$.
The first order approach is new while the second order approach formalizes the results in \cite{FrankenPhDThesis,QAM2}. They also serve as an introduction to the new, more technical, $SE(3)$-extensions in Section~\ref{ch:SE3}.

All curve optimization problems are based on the idea that a curve (or a family of curves) fits the data well if a certain quantity is preserved along the curve. This preserved quantity is the data $\tilde{U}(\tilde{\gamma}(t))$ for the first order optimization, and the time derivative
$\frac{d}{dt}\tilde{U}(\tilde{\gamma}(t))$ or the gradient $\nabla\tilde{U}(\tilde{\gamma}(t))$ for the second order optimization.  After introducing a family of curves similar to the ones used in Section \ref{ch:gauge} we will, for all three cases, first pose an optimization problem, and then give its solution in a subsequent theorem.

In this section we rely on the group-theoretical tools explained in Section \ref{ch:geomTools} (only the case d=2), listed in subtables E.1 and E.2 in our table of notations. Furthermore we introduce notations listed in the first part of subtable E.3.
\subsection{Neighboring Exponential Curves in $SE(2)$ \label{ch:FamilyExponentialCurvesSE2}}

\begin{figure}[htbp]
	\centering
	\includegraphics[width=0.95\hsize]{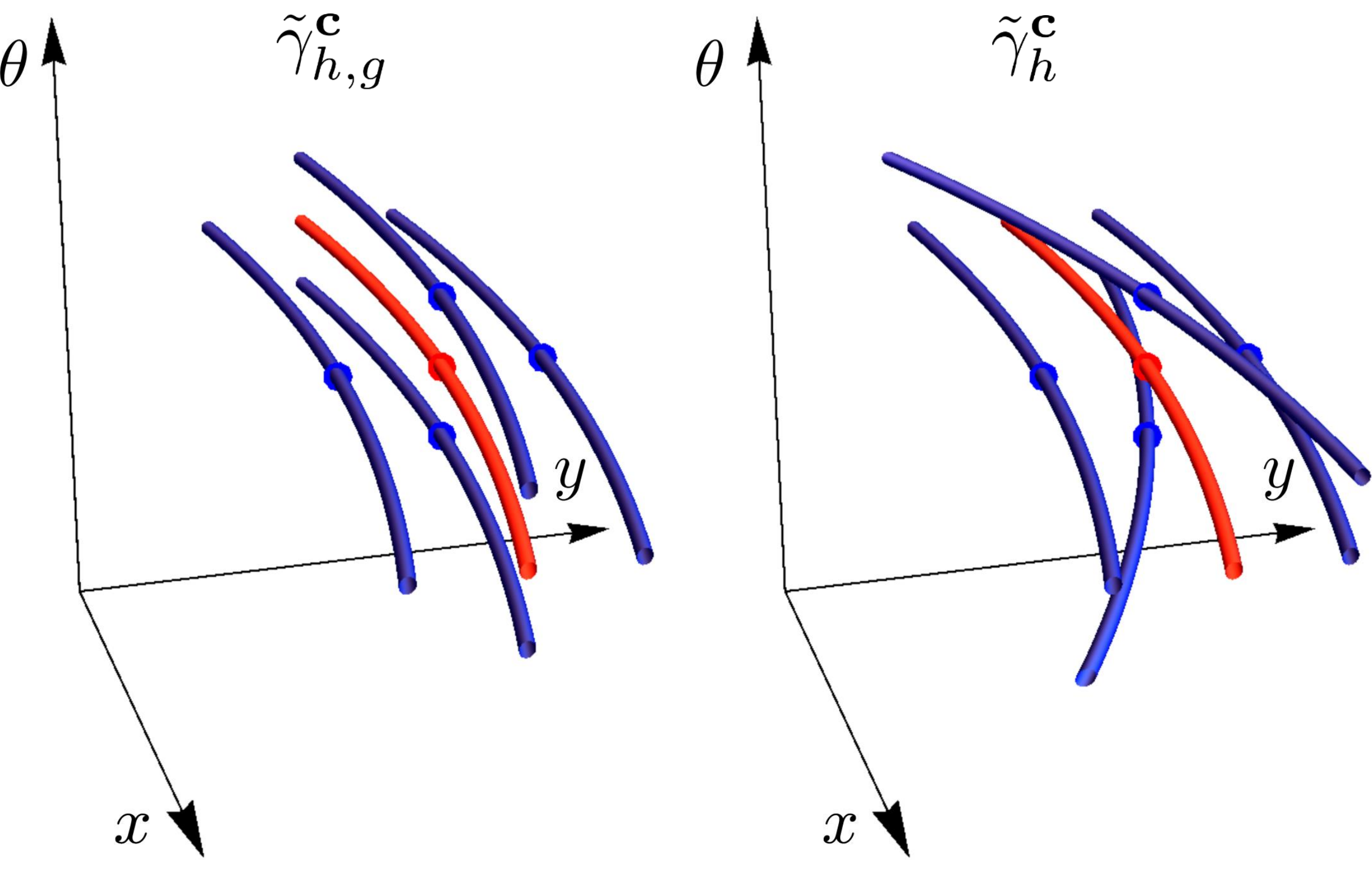}
	\caption{Family of neighboring exponential curves, given a fixed point $g \in SE(2)$ and a fixed tangent vector
$\ul{c}=\ul{c}(g) \in T_{g}(SE(2))$.
Left: Our choice of family of exponential curves $\gamma_{h,g}^{\ul{c}}$ for neighboring $h \in SE(2)$.
Right: Exponential curves  $\gamma^{\ul{c}}_h$ with $\ul{c}=\ul{c}(g)$ are not suited for local averaging in our curve fits.
The red curves start from $g$ (indicated with a dot), the blue curves from $h\neq g$.
}
\label{fig:gonz}
\end{figure}

Akin to (\ref{paralleltransportRd}) we fix reference point $g \in SE(2)$ and velocity components $\ul{c}=\ul{c}(g) \in \R^3$, and we shall rely on a family $\{\tilde{\gamma}_{h,g}^{\ul{c}}\}$ of neighboring exponential curves around $\tilde{\gamma}^{\ul{c}}_g$. As we will show in subsequent Lemma \ref{lemma:SE2} neighboring curve
$\tilde{\gamma}^{\ul{c}}_{h,g}$ departs from $h$ and has the same spatial and rotational velocity as the curve
$\tilde{\gamma}^{\ul{c}}_g$ departing from $g$. This geometric idea is visualized in Fig.~\ref{fig:gonz}, where it is intuitively explained why one needs
the initial velocity vector $\tilde{\ul{R}}_{h^{-1}g}\ul{c}$, instead of $\ul{c}$ in the following definition for the exponential curve departing from a neighboring point $h$ close to $g$.
\begin{definition}\label{def:SE2neighbors}
Let $g \in SE(2)$ and $\ul{c}=\ul{c}(g) \in \R^{3}$ be given. Then we define the family $\{\tilde{\gamma}_{h,g}^{\ul{c}}\}$
of neighboring exponential curves
\begin{equation} \label{paralleltransportSE2}
t \mapsto \tilde{\gamma}_{h,g}^{\ul{c}}(t):=\tilde{\gamma}_{h}^{\tilde{\ul{R}}_{h^{-1}g} \ul{c}}(t),
\end{equation}
with rotation-matrix $\tilde{\ul{R}}_{h^{-1}g} \in SO(3)$ defined by
\begin{equation}\label{RtildeSE2}
              \tilde{\ul{R}}_{h^{-1}g} := \begin{pmatrix}
              (\ul{R}')^T \ul{R} &&& \ul{0} \\
              \ul{0} &&& 1
              \end{pmatrix},
\end{equation}
for all $g=(\ul{x},\ul{R}) \in SE(2)$ and all $h=(\ul{x}',\ul{R}') \in SE(2)$, with $\ul{R},\ul{R}' \in SO(2)$ a counterclockwise rotation by respectively angle $\theta$ and $\theta'$.
\end{definition}
\begin{lemma}\label{lemma:SE2}
Exponential curve $\tilde{\gamma}_{h,g}^{\ul{c}}$ departing from $h\in SE(2)$ given by (\ref{paralleltransportSE2}) has the same spatial and angular velocity as exponential curve $\tilde\gamma^{\ul{c}}_g$ departing from \mbox{$g \in SE(2)$}.

On the Lie algebra level; we have that the initial velocity component vectors of the curves $\tilde\gamma^{\ul{c}}_g$ and $\tilde{\gamma}_{h,g}^{\ul{c}}$ relate via $\ul{c} \mapsto \tilde{\ul{R}}_{h^{-1}g} \ul{c}$.

On the Lie group level; we have that the curves themselves $\tilde\gamma^{\ul{c}}_g(\cdot)=(\ul{x}_{g}(\cdot),\ul{R}_{g}(\cdot))$, $\tilde{\gamma}_{h,g}^{\ul{c}}(\cdot)=
(\ul{x}_{h}(\cdot),\ul{R}_{h}(\cdot))$ relate via
\begin{equation} \label{2star}
\begin{array}{ll}
\ul{x}_{h}(t) &=\ul{x}_{g}(t)-\ul{x}+\ul{x}', \\
\ul{R}_{h}(t) &=\ul{R}_{g}(t) \ul{R}^{-1}\ul{R}' \desda \theta_h(t)=\theta_g(t) -\theta +\theta'.
\end{array}
\end{equation}
\end{lemma}
\textbf{Proof }The proof follows from the proof of a more general theorem on the $SE(3)$ case which follows later (in Lemma~\ref{lemma:SE3}).
\begin{remark}
Eq.~\!(\ref{2star}) is the extension of Eq.~\!(\ref{paralleltransportRd}) on $\R^{2}$ to the $SE(2)$ group.
\end{remark}
Additional geometric background is given in Appendix~\ref{app:B}.

\subsection{Exponential Curve Fits in $SE(2)$ of the 1st Order \label{ch:SE2-order1}}
For first-order exponential curve fits we solve an optimization problem similar to (\ref{optR2}) given by
\begin{equation}\label{minproblemSE2firstorder}
\boxed{
\ul{c}^{*}(g) =
\argmin_{\scriptsize \begin{array}{c} \ul{c}\in \R^{3}, \\ \|\ul{c}\|_{\mu}\!=\!1 \end{array}} \int \limits_{\;\;SE(2)} \!\!\!\!\!\!
\tilde{G}_{\boldsymbol{\rho}}(h^{-1}g)
 \left|\left.\frac{d}{dt} \tilde{V}(\tilde{\gamma}_{h,g}^{\ul{c}}(t)) \right|_{t=0}\right|^{2}
\!{\rm d}\overline{\mu}(h),
}
\end{equation}
with $\tilde{V}=\tilde{G}_{\ul{s}}*\tilde{U}$, $g=(\ul{x},\ul{R})$, $h=(\ul{x}',\ul{R}')$ and ${\rm d}\overline{\mu}(h)={\rm d}\ul{x}' {\rm d}\mu_{SO(2)}(\ul{R}')={\rm d}\ul{x}' {\rm d}\theta'$.
Here we first regularize the data with spatial and angular scale $\ul{s}=(s_{p},s_o)$ and then average over a family of curves where we use spatial and angular scale $\boldsymbol{\rho}=(\rho_{p},\rho_o)$. Here $s_{p},\rho_p>0$ are isotropic scales on $\R^{2}$ and $s_{o}, \rho_o>0$ are scales on $S^{1}$ of separable Gaussian kernels, recall (\ref{Gtilde}).
Recall also (\ref{normmu}) for the definition of the norm $\|\cdot\|_{\mu}$.
When solving this optimization problem the following structure matrix appears
\begin{equation} \label{tensor}
\begin{array}{ll}
(\ul{S}^{\ul{s},\boldsymbol{\rho}}(\tilde{U}))(g) = \int \limits_{SE(2)}  \tilde{G}_{\boldsymbol{\rho}}(h^{-1} g) \cdot \\
\qquad \tilde{\ul{R}}_{h^{-1} g}^T  \nabla \tilde{V}(h)  (\nabla \tilde{V}(h))^{T} \tilde{\ul{R}}_{h^{-1} g} \, {\rm d}\overline{\mu}(h).
\end{array}
\end{equation}
In the remainder we use short notation $\ul{S}^{\ul{s},\boldsymbol{\rho}}:=\ul{S}^{\ul{s},\boldsymbol{\rho}}(\tilde{U})$. We assume that $\tilde{U}$, $\boldsymbol{\rho}, \ul{s}$, and $g$, are chosen such that
$\ul{S}^{\ul{s},\boldsymbol{\rho}}(g)$ is a non-degenerate matrix. The optimization problem is solved in the next theorem.
\begin{theorem}[First Order Fit via Structure \mbox{Tensor}] \label{th:0}
The normalized eigenvector $\ul{M}_{\mu}\ul{c}^{*}(g)$ with smallest eigenvalue of the rescaled structure matrix \\
$\ul{M}_{\mu}\ul{S}^{\ul{s},\boldsymbol{\rho}}(g)\ul{M}_{\mu}$
provides the solution $\ul{c}^{*}(g)$ to optimization problem (\ref{minproblemSE2firstorder}).
\end{theorem}
\textbf{Proof }
We will apply four steps. In the first step we write the time-derivative as a directional derivative, in the second step we express the directional derivative in the gradient. In the third step we put the integrand in matrix-vector form. In the final step we express our optimization functional in the structure tensor and solve the Euler-Lagrange equations.

For the first step we use (\ref{paralleltransportSE2}) and the fundamental property (\ref{ec}) of exponential curves such that via application of (\ref{ddtcurve}):
\begin{equation}
  \begin{array}{rl}
  \left| \left.\frac{d}{dt} \left( \tilde{V}(\tilde{\gamma}^{\ul{c}}_{h,g}(t)) \right)\right|_{t=0} \right|^2  &=
  \left| \langle {\rm d}\tilde{V}|_{
  \tilde{\gamma}^{\ul{c}}_{h,g}(0)}, \dot{\tilde{\gamma}}^{\ul{c}}_{h,g}(0)\rangle \right|^2 \\
    &= \left| \langle {\rm d}\tilde{V}|_{h},\tilde{\mathbf{R}}_{h^{-1}g} \ul{c} \rangle \right|^2,
  \end{array}
\end{equation}
where we use short notation $\tilde{\mathbf{R}}_{h^{\!-\!1}g} \ul{c} \! = \! \sum_{i=1}^{3} (\tilde{\mathbf{R}}_{h^{\!-\!1}g} \ul{c})^i \!\mathcal{A}_i |_h$.

In the second step we use the definition of the gradient (\ref{gradientmu2}) and the metric tensor (\ref{metrictensor2}) to rewrite this expression to
\begin{equation}
  \begin{array}{l}
  \left| \langle {\rm d}\tilde{V}|_{h},\tilde{\mathbf{R}}_{h^{-1}g} \ul{c}) \rangle \right|^2
     = \left| \left.\gothic{G}_{\mu}\right|_{h}(\nabla \tilde{V}(h), \tilde{\mathbf{R}}_{h^{-1}g} \ul{c} ) \right|^{2} .
  \end{array}
\end{equation}

Then, in the third step we write this in vector-matrix form and obtain
\begin{equation}
  \begin{array}{l}
  \left| \left.\gothic{G}_{\mu}\right|_{h}\!(\nabla \tilde{V}(h), \tilde{\mathbf{R}}_{h^{\!-\!1}g} \ul{c} ) \right|^{2} \!\!= \left| \ul{c}^{T} \ul{M}_{\mu^{2}} \tilde{\mathbf{R}}_{h^{\!-\!1}g}^{T} \nabla \tilde{V}(h) \right|^{2}  \\
 =\ul{c}^{T} \ul{M}_{\mu^{2}} \tilde{\mathbf{R}}_{h^{\!-\!1}g}^{T} \nabla \tilde{V}(h) (\nabla \tilde{V}(h))^{T} \tilde{\mathbf{R}}_{h^{\!-\!1}g} \ul{M}_{\mu^{2}} \ul{c},
  \end{array}
\end{equation}
where we used the fact that $\ul{M}_{\mu^{2}}$ and $\tilde{\mathbf{R}}_{h^{\!-\!1}g}^{T}$ commute.

Finally, we use the structure tensor definition (\ref{tensor}) to rewrite the convex optimization functional in (\ref{minproblemSE2firstorder}) as
\begin{equation}
\begin{array}{rl}
\mathcal{E}(\ul{c}):=& \!\! \int \limits_{SE(2)} \!\!
\tilde{G}_{\boldsymbol{\rho}}(h^{-1}g)
 \left|\left.\frac{d}{dt} \tilde{V}(\tilde{\gamma}_{h,g}^{\ul{c}}(t)) \right|_{t=0}\right|^{2}
\!{\rm d}\overline{\mu}(h) \\
 =&\ul{c}^T \ul{M}_{\mu^2} \ul{S}^{\ul{s},\boldsymbol{\rho}} \ul{M}_{\mu^2} \ul{c},
\end{array}
\end{equation}
while the boundary condition $\|\ul{c}\|_{\mu}\!=\!1$ can be written as
\begin{equation}
\varphi(\ul{c}):=\ul{c}^{T} \ul{M}_{\mu^2} \ul{c}-1=0.
\end{equation}
The Euler-Lagrange equation reads $\nabla \mathcal{E}(\ul{c}^{*})=\lambda_{1}\nabla\varphi(\ul{c}^{*})$, with $\lambda_{1}$ the smallest eigenvalue
of $\ul{M}_{\mu}\ul{S}^{\ul{s},\boldsymbol{\rho}}(g)\ul{M}_{\mu}$ and we have
\begin{equation} \label{corr:1B}
\begin{array}{rcl}
\mathbf{M}_{\mu^{2}} \ul{S}^{\ul{s},\boldsymbol{\rho}}(g) \mathbf{M}_{\mu^{2}} \ul{c}^{*}(g)&=& \lambda_{1}\, \mathbf{M}_{\mu^2} \ul{c}^{*}(g) \vspace{-1mm} \\
&\Updownarrow& \vspace{-1mm} \\
\mathbf{M}_{\mu} \ul{S}^{\ul{s},\boldsymbol{\rho}}(g) \mathbf{M}_{\mu} (\mathbf{M}_{\mu} \ul{c}^{*}(g))&=&
\lambda_{1}(\mathbf{M}_{\mu} \ul{c}^{*}(g)),
\end{array}
\end{equation}
from which the result follows. 
$\hfill \Box$
\\
\\
The next remark explains the frequent presence of the $\ul{M}_{\mu}$ matrices in (\ref{corr:1B}).
\begin{remark}
The diagonal $\ul{M}_{\mu}$ matrices enter the functional due to the gradient definition (\ref{gradientmu2}), and they enter the boundary condition via $\|\ul{c}\|_{\mu}^2=\ul{c}^T \ul{M}_{\mu^2} \ul{c}=1$. In both cases they come from the metric tensor (\ref{metrictensor2}). Parameter $\mu$ which controls the stiffness of the exponential curves has physical dimension $[\textrm{Length}]^{-1}$. As a result, the normalized eigenvector $\mathbf{M}_{\mu} \ul{c}^{*}(g)$ is, in contrast to $\ul{c}^{*}(g)$, dimensionless.
\end{remark}

\subsection{Exponential Curve Fits in $SE(2)$ of the 2nd Order\label{ch:SE2-order2}}

We now discuss the second order optimization approach based on the Hessian matrix.
At each $g\in SE(2)$ we define a $3\times 3$ non-symmetric Hessian matrix
\begin{equation}\label{ourhessian}
	(\ul{H}^{\ul{s}}(\tilde{U}))(g)= \left[\mathcal{A}_{j}\mathcal{A}_{i} (\tilde{V})(g)\right] \,, \;\;\; \textrm{with } \tilde{V}=\tilde{G}_{\ul{s}}*\tilde{U},
\end{equation}
and where $i$ denotes the row index and where $j$ denotes the column index,
and with $\tilde{G}_{\ul{s}}$ a Gaussian kernel with isotropic spatial part as described in Eq. (\ref{Gtilde}). In the remainder we write $\ul{H}^{\ul{s}} :=\ul{H}^{\ul{s}}(\tilde{U})$.
\begin{remark}
As the left-invariant vector fields are non-commutative
there are many ways to define the Hessian matrix on $SE(2)$, since the ordering of the left-invariant derivatives matters.
From a differential geometrical point of view our choice (\ref{ourhessian}) is correct, as we motivate
in Appendix~\ref{app:remcohappy}.
\end{remark}
For second-order exponential curve fits we consider 2 different optimization problems. In the first case we minimize the second order derivative along the exponential curve:
\begin{equation}\label{SE2secondOrderOptimization}
\boxed{
\ul{c}^{*}(g) =
\argmin_{\ul{c}\in \R^{3}, \|\ul{c}\|_{\mu}=1} 
\left|\;\left.\frac{d^2}{dt^2} \tilde{V}(\tilde{\gamma}_{g}^{\ul{c}}(t)) \right|_{t=0}\;\right|.
}
\end{equation}
In the second case we minimize the norm of the first order derivative of the gradient
of the neighboring family of exponential curves:
\begin{equation}\label{SE2secondOrderOptimizationGradient}
\boxed{
\begin{array}{l}
\ul{c}^{*}(g) =
\argmin_{\ul{c}\in \R^{3}, \|\ul{c}\|_{\mu}=1} \,\,\, \int \limits_{SE(2)} \tilde{G}_{\boldsymbol{\rho}}(h^{-1}g) \cdot \\
\gothic{G}_{\mu}\left(\left.\frac{d}{dt} \nabla \tilde{V}(\tilde{\gamma}_{h,g}^{\ul{c}}(t))\right|_{t=0},\left.\frac{d}{dt} \nabla \tilde{V}(\tilde{\gamma}_{h,g}^{\ul{c}}(t))\right|_{t=0}\right)\,
{\rm d}\overline{\mu}(h),
\end{array}
}	
\end{equation}
with again $\tilde{V}=\tilde{G}_{\ul{s}}*\tilde{U}$.
\begin{remark}
Optimization problem (\ref{SE2secondOrderOptimization}) can also be written as an optimization problem over the neighboring family of curves, as it is equivalent to problem:
\begin{equation}\label{SE2secondOrderOptimization2}
\ul{c}^{*}(g) = \! \argmin_{\scriptsize \begin{array}{c} \ul{c}\in \R^{3}, \\ \|\ul{c}\|_{\mu}\!=\!1 \end{array} }
\Bigg| \int \limits_{SE(2)} \!\!\!\!\! \tilde{G}_{\ul{s}}(h^{-1}g)
\left.\frac{d^2}{dt^2} \tilde{U}(\tilde{\gamma}_{h,g}^{\ul{c}}(t)) \right|_{t=0}\; \!\!\! {\rm d}\overline{\mu}(h) \Bigg| .
\end{equation}
\end{remark}

In the next two theorems we solve these optimization problems.

\begin{theorem}[Second Order Fit via Symmetric Sum Hessian] \label{th:0b}
Let $g \in SE(2)$ be such that the eigenvalues of the rescaled symmetrized Hessian
\[\frac{1}{2} \ul{M}_{\mu}^{-1}(
\ul{H}^{\ul{s}}(g)+ (\ul{H}^{\ul{s}}(g))^T) \ul{M}_{\mu}^{-1}
\]
have the same sign. Then the normalized eigenvector $\ul{M}_{\mu}\ul{c}^{*}(g)$ with smallest eigenvalue of the rescaled symmetrized Hessian matrix
provides the solution $\ul{c}^{*}(g)$ of optimization problem (\ref{SE2secondOrderOptimization}).
\end{theorem}
\textbf{Proof }Similar to the proof of Theorem~\ref{th:0} we first write the time derivative as a directional derivative using Eq. (\ref{ddtcurve}). Since now we have a second order derivative this step is applied twice:
\begin{equation} \label{Proof2Step1}
\begin{array}{ll}
\bigg|\; \frac{d^2}{dt^2} \tilde{V}(\tilde{\gamma}_{g}^{\ul{c}}(t))\, \Big|_{t=0}\; \bigg| &=
\bigg|\; \frac{d}{dt} \sum \limits_{i=1}^{3} c^{i}\mathcal{A}_{i} \tilde{V}(\tilde{\gamma}_{g}^{\ul{c}}(t)) \Big|_{t=0} \; \bigg|\\
&=\bigg| \sum \limits_{i,j=1}^{3} c^{i}c^j \mathcal{A}_{j}(\mathcal{A}_{i}\tilde{V})(g) \bigg|.
\end{array}
\end{equation}
Then we write the result in matrix-vector form and split the matrix in a symmetric and anti-symmetric part
\begin{equation} \label{Proof2Step2}
\begin{array}{l}
\left| \sum \limits_{i,j=1}^{3} c^{i}c^j \mathcal{A}_{j}(\mathcal{A}_{i}\tilde{V})(g) \right|
=\left|\ul{c}^{T} \ul{H}^{\ul{s}}(g)  \ul{c} \right| \\[5pt]
=\left| \frac{1}{2}\ul{c}^{T} (\ul{H}^{\ul{s}}(g)+ (\ul{H}^{\ul{s}}(g))^T)  \ul{c} \right. \\
\qquad +
\left.
\frac{1}{2} \ul{c}^{T} (\ul{H}^{\ul{s}}(g)- (\ul{H}^{\ul{s}}(g))^T)  \ul{c} \right| \\
=\frac{1}{2}\left|\ul{c}^{T} (\ul{H}^{\ul{s}}(g)+ (\ul{H}^{\ul{s}}(g))^T)  \ul{c} \right|,
\end{array}
\end{equation}
where only the symmetric part remains. Finally, the optimization functional in (\ref{SE2secondOrderOptimization}) (which is convex if the eigenvalues have the same sign) can be written as
\begin{equation}
\begin{array}{rl}
 \mathcal{E}(\ul{c}) :=& \left|\;\left.\frac{d^2}{dt^2} \tilde{V}(\tilde{\gamma}_{g}^{\ul{c}}(t)) \right|_{t=0}\;\right| \\
=& \frac{1}{2}\left|\ul{c}^{T} (\ul{H}^{\ul{s}}(g)+ (\ul{H}^{\ul{s}}(g))^T)  \ul{c} \right|.
\end{array}
\end{equation}
Again we have the boundary condition $\varphi(\ul{c})=\ul{c}^T \ul{M}_{\mu^2} \ul{c}-1=0$. The result follows using the Euler-Lagrange formalism
$\nabla \mathcal{E}(\ul{c}^{*})=\lambda_{1}\nabla\varphi(\ul{c}^{*})$:
\begin{equation} \label{corr:1B}
\begin{array}{l}
\frac{1}{2}  (\ul{H}^{\ul{s}}(g)+ (\ul{H}^{\ul{s}}(g))^T) \ul{c}^{*}(g)= \lambda_{1}\, \mathbf{M}_{\mu^2} \ul{c}^{*}(g) \Leftrightarrow \\[7pt]
\frac{1}{2}  \mathbf{M}_{\mu}^{\!-\!1} (\ul{H}^{\ul{s}}(g) \! + \! (\ul{H}^{\ul{s}}(g))^T) \mathbf{M}_{\mu}^{\!-\!1} (\mathbf{M}_{\mu} \ul{c}^{*}(g))\\
=
\lambda_{1}(\mathbf{M}_{\mu} \ul{c}^{*}(g)),
\end{array}
\end{equation}
which boils down to finding the eigenvector with minimal absolute eigenvalue $|\lambda_{1}|$ which gives our result.
$\hfill \Box$
\begin{theorem}[Second Order Fit via Symmetric \\
\mbox{Product} Hessian] \label{th:0c}
Let $\rho_p,\rho_o, s_p,s_o>0$.
The normalized eigenvector $\ul{M}_{\mu}\ul{c}^{*}(g)$ with smallest eigenvalue of matrix
\begin{equation}
\begin{array}{ll}
\ul{M}_{\mu}^{-1} \int \limits_{SE(2)} \tilde{G}_{\boldsymbol{\rho}}(h^{-1}g) \cdot  \tilde{\ul{R}}_{h^{-1} g}^T (\ul{H}^{\ul{s}}(h))^T \\ \ul{M}_{\mu}^{-2} \ul{H}^{\ul{s}}(h) \tilde{\ul{R}}_{h^{-1} g} \; {\rm d}\overline{\mu}(h) \,\ul{M}_{\mu}^{-1}
\end{array}
\end{equation}
provides the solution $\ul{c}^{*}(g)$ of optimization problem (\ref{SE2secondOrderOptimizationGradient}).
\end{theorem}
\textbf{Proof }
First we use the definition of the gradient (\ref{gradientmu2}) and then we again rewrite the time-derivative as a directional derivative:
\begin{equation}\label{ddtgradientSE2}
\begin{array}{ll}
	\left. \frac{d}{dt} \nabla \tilde{V}(\tilde{\gamma}_{h,g}^{\ul{c}}(t))\right|_{t=0} &=
 \frac{d}{dt} \sum \limits_{i=1}^{3} \mathcal{A}_{i} \tilde{V}(\tilde{\gamma}_{g}^{\ul{c}}(t)) \mu_i^{-2} \mathcal{A}_i|_{\tilde{\gamma}_{g}^{\ul{c}}(t)} \bigg|_{t=0} \\
  &=\sum \limits_{i,j=1}^3 \tilde{\ul{c}}^j \mathcal{A}_j \mathcal{A}_i \tilde{V}(h) \mu_i^{-2} \mathcal{A}_i|_h
  \end{array}
\end{equation}
for $\tilde{\ul{c}} = \tilde{\ul{R}}_{h^{-1} g} \ul{c}$, recall (\ref{RtildeSE2}), and where $\mu_i=\mu$ for $i=1,2$ and $\mu_i=1$ for $i=3$. Here we use $\tilde{\gamma}_{h,g}^{\ul{c}}(0)=h$, and the formula for left-invariant vector fields (\ref{LINV}). Now insertion of (\ref{ddtgradientSE2}) into the metric tensor $\gothic{G}_{\mu}$ (\ref{metrictensor2}) yields
\begin{equation}
\begin{array}{l}
	\gothic{G}_{\mu} \left(\left.\frac{d}{dt} \nabla \tilde{V}(\tilde{\gamma}_{h,g}^{\ul{c}}(t))\right|_{t=0},\left.\frac{d}{dt} \nabla \tilde{V}(\tilde{\gamma}_{h,g}^{\ul{c}}(t))\right|_{t=0}\right) \\
	\,\,\,\,=\tilde{\ul{c}}^T (\ul{H}^{\ul{s}}(h))^T \ul{M}_{\mu}^{-2} \ul{H}^{\ul{s}}(h) \tilde{\ul{c}} \\
	\,\,\,\,= \ul{c}^T \tilde{\ul{R}}_{h^{-1} g}^T (\ul{H}^{\ul{s}}(h))^T \ul{M}_{\mu}^{-2} \ul{H}^{\ul{s}}(h) \tilde{\ul{R}}_{h^{-1} g} \ul{c}.
\end{array}
\end{equation}
Finally, the convex optimization functional in (\ref{SE2secondOrderOptimizationGradient}) can be written as
\begin{equation}
\begin{array}{l}
 \mathcal{E}(\ul{c}) :=
\ul{c}^T \bigg( \int \limits_{SE(2)} \tilde{G}_{\boldsymbol{\rho}}(h^{-1}g) \cdot  \tilde{\ul{R}}_{h^{-1} g}^T (\ul{H}^{\ul{s}}(h))^T \\
\qquad \ul{M}_{\mu}^{-2} \ul{H}^{\ul{s}}(h) \tilde{\ul{R}}_{h^{-1} g} \; {\rm d}\overline{\mu}(h) \bigg) \,\ul{c}.
\end{array}
\end{equation}
Again we have the boundary condition $\varphi(\ul{c})=\ul{c}^T \ul{M}_{\mu^2} \ul{c}-1=0$ and the result follows by application of the Euler-Lagrange formalism: $\nabla \mathcal{E}(\ul{c}^{*})=\lambda_{1}\nabla\varphi(\ul{c}^{*})$.
 $\hfill \Box$ 

\section{Exponential Curve Fits in $SE(3)$ \label{ch:SE3}}
In this section we generalize the exponential curve fit theory from the preceding chapter on $SE(2)$ to $SE(3)$. Because our data on the group $SE(3)$ was obtained from data on the quotient $\R^{3}\rtimes S^{2}$ we will also discuss projections of exponential curve fits on the quotient.

We start in Subsection~\ref{ch:introquotient} with some prelimenaries on the quotient structure (\ref{PositionsAndOrientations}). Here we will also introduce the concept of projected exponential curve fits.
Subsequently, in Subsection~\ref{ch:FamilyExponentialCurvesSE3}, we provide basic theory on how to obtain the appropriate family of neighboring exponential curves. More details can be found in Appendix~\ref{app:B}.
In Subsection~\ref{ch:1} we formulate exponential curve fits of the first order as a variational problem.
For that we define the structure tensor on $SE(3)$, which we use to solve the variational problem in Theorems~\ref{th:2} and \ref{th:3b}.
Then we present the two-fold algorithm for achieving torsion-free exponential curve fits.
In Subsection~\ref{ch:SE3-order2} we formulate exponential curve fits of the second order as a variational problem.
Then we define the Hessian tensor on $SE(3)$, which we use to solve the variational problem in Theorem~\ref{th:6}.
Again torsion-free exponential curve fits are accomplished via a two-fold algorithm.

Throughout this section we will rely on the differential geometrical tools of Section~\ref{ch:geomTools}, listed in subtables~\ref{app:TableOfNotations}\!.1 and \ref{app:TableOfNotations}\!.2 in Appendix~\ref{app:TableOfNotations}.
We also generalize concepts on exponential curve fits introduced in the previous section to the case $d=3$ (requiring additional notation). They are listed in subtable~\ref{app:TableOfNotations}\!.3 in Appendix~\ref{app:TableOfNotations}.

\subsection{Preliminaries on the quotient $\R^{3} \rtimes S^2$.}\label{ch:introquotient}

Now let us set $d=3$, and let us assume input $U$ is given and let us first concentrate on its domain. This domain equals the
joint space $\R^{3} \rtimes S^{2}$ of positions and orientations of dimension $5$, which we identified with a 5-dimensional group quotient of $SE(3)$, where $SE(3)$ is of dimension $6$ (recall (\ref{PositionsAndOrientations})). For including a notion of alignment it is crucial to include the non-commutative relation in (\ref{groupproduct}) between rotations and translation, and not to consider the space of positions and orientations as a flat Cartesian product. Therefore we model the joint space of positions and orientations as the Lie group quotient (\ref{PositionsAndOrientations}), where
\[
SO(2)\equiv \textrm{Stab}(\ul{a})= \{\ul{R} \in SO(3)\;|\; \ul{R} \ul{a}=\ul{a}\}
\]
for reference axis $\ul{a}=\ul{e}_{z}=(0,0,1)^T$.
Within this quotient structure two rigid body motions $g=(\ul{x},\ul{R}), g'=(\ul{x}',\ul{R}') \in SE(3)$ are  equivalent if
\[
\begin{array}{l}
g' \sim  g \desda (g')^{-1}g \in \{\ul{0}\} \times SO(2)  \desda \\[6pt] 
\ul{x}-\ul{x}'=\ul{0} \textrm{ and }
\exists_{\ul{R}_{\ul{e}_{z},\alpha} \in SO(2)} \; :\;
(\ul{R}')^{-1}\ul{R}=\ul{R}_{\ul{e}_{z},\alpha}.
\end{array}
\]
Furthermore, one has the action $\odot$ of $g=(\ul{x},\ul{R}) \in SE(3)$ onto $(\ul{y},\ul{n})  \in \R^{3} \times S^{2}$, which is defined by
\begin{equation} \label{odot}
\begin{array}{l}
g \odot (\ul{y},\ul{n})=(\ul{x},\ul{R}) \odot(\ul{y},\ul{n}):=(\ul{x}+ \ul{R}\ul{y}, \ul{R}\ul{n}).
\end{array}
\end{equation}
As a result we have
\[
g' \sim g \desda g' \odot (\ul{0},\ul{a})= g \odot (\ul{0},\ul{a}).
\]
Thereby, a single element in $\R^{3} \rtimes S^{2}$ can be considered as the equivalence class of all rigid body motions that map reference position and orientation $(\ul{0},\ul{a})$ onto $(\ul{x},\ul{n})$.
Similar to the common identification of $S^{2}\equiv SO(3)/SO(2)$, we denote elements of the Lie group quotient $\R^{3}\rtimes S^{2}$
by $(\ul{x},\ul{n})$.

\subsubsection{Legal Operators}
Let us recall from Section~\ref{ch:togauge} that exponential curve fits induce gauge frames.
Note that both the induced gauge frame
$\{\mathcal{B}_1,\ldots, \mathcal{B}_{6}\}$ and the non-adaptive frame
$\{\mathcal{A}_{1},\ldots,\mathcal{A}_{6}\}$ are defined on the Lie group $SE(3)$, and cannot be defined on the quotient. Nevertheless, combinations of them can be well-defined on $\R^{3}\rtimes S^2$ (e.g.
$\Delta_{\R^{3}}=\mathcal{A}_{1}^{2}+\mathcal{A}_{2}^{2}+\mathcal{A}_{3}^{2}$ is well-defined on the quotient).
This brings us to the definition of so-called \emph{legal} operators, as shown in \cite[Thm.1]{DuitsJMIV}.
In short, the operator $\tilde{U} \mapsto \tilde{\Phi}(\tilde{U})$ is legal (left-invariant and well-defined on the quotient) if and only if
\begin{equation} \label{legal}
\begin{array}{l}
\tilde{\Phi}= \tilde{\Phi} \circ \mathcal{R}_{h_\alpha} \textrm{ for all }\alpha \in [0,2\pi). \\
\tilde{\Phi} \circ \mathcal{L}_{g} = \mathcal{L}_{g} \circ \tilde{\Phi} \textrm{ for all }g \in SE(3),
\end{array}
\end{equation}
recall (\ref{eq:leftrightregularrepr}), where
\begin{equation}\label{HALPHA}
h_{\alpha}:=(\ul{0},\ul{R}_{\ul{e}_{z},\alpha}).
\end{equation}
 with the $\ul{R}_{\ul{e}_{z},\alpha}$ the counterclockwise rotation about $\ul{e}_{z}$. Such legal operators relate one-to-one to operators
$\Phi: \mathbb{L}_{2}(\R^{3} \rtimes S^{2}) \to \mathbb{L}_{2}(\R^{3} \rtimes S^{2})$ via
\[
U \mapsto \Phi (U) \; \; \leftrightarrow \; \; \tilde{U} \mapsto \tilde{\Phi}(\tilde{U})=\widetilde{\Phi(U)},
\]
relying consequently on (\ref{tildeU}).

\subsubsection{Projected Exponential Curve Fits}
Action (\ref{odot}) allows us to map a curve $\tilde{\gamma}(\cdot)=(\ul{x}(\cdot),\ul{R}(\cdot))$ in $SE(3)$ onto a curve $(\ul{x}(\cdot),\ul{n}(\cdot))$ on
$\R^{3} \rtimes S^2$ via
\begin{equation} \label{final2}
\begin{array}{ll}
(\ul{x}(t), \ul{n}(t)) &:= \tilde{\gamma}(t)\odot (\ul{0},\ul{e}_{z})=(\ul{x}(t), \ul{R}(t)\, \ul{e}_{z}).
 \end{array}
\end{equation}
This can be done with exponential curve fits $\tilde{\gamma}^{\ul{c}=\ul{c}^{*}(g)}_{g}(t)$ to define projected exponential curve fits.
\begin{definition}
We define for $g=(\ul{x},\ul{R}_{\ul{n}})$ the projected exponential curve fit
\begin{equation} \label{projectedCurveFit}
\gamma_{(\ul{x},\ul{n})}^*(t):=
\tilde{\gamma}_{g}^{\ul{c}^{*}(g)}(t) \odot (\ul{0},\ul{e}_{z}).
\end{equation}
\end{definition}

\begin{lemma}\label{lemmaProjectedCurveFit}
The projected exponential curve fit is well-defined on the quotient, i.e. the right-hand side of (\ref{projectedCurveFit}) is independent of the choice of $\ul{R}_{\ul{n}}$ s.t. $\ul{R}_{\ul{n}} \ul{e}_z = \ul{n}$, if the optimal tangent found in our fitting procedure satisfies:
\begin{equation} \label{constraint}
\ul{c}^{*}(g h_{\alpha})= \ul{Z}_{\alpha}^T \ul{c}^{*}(g), \qquad \textrm{for all } \alpha \in [0,2\pi],
\end{equation}
and for all $g\in SE(3)$, with
\begin{equation} \label{Zalphadef}
\ul{Z}_{\alpha}:= \begin{pmatrix}
\ul{R}_{\ul{e}_{z},\alpha} && \ul{0} \\
\ul{0} && \ul{R}_{\ul{e}_{z},\alpha}
\end{pmatrix} \in SO(6).
\end{equation}
\end{lemma}
\textbf{Proof }
For well-posed projected exponential curve fits we need the right-hand side of (\ref{projectedCurveFit}) to be independent of  $\ul{R}_{\ul{n}}$ s.t. $\ul{R}_{\ul{n}} \ul{e}_z = \ul{n}$ i.e. it should be invariant under $g \rightarrow g h_{\alpha}$. Therefore we have the following constraint on the fitted curves:
\begin{equation} \label{projectedCurveFitConstraint}
\tilde{\gamma}_{g}^{\ul{c}^{*}(g)}(t) \odot (\ul{0},\ul{e}_{z})=\tilde{\gamma}_{g h_{\alpha}}^{\ul{c}^{*}(g h_{\alpha})}(t) \odot (\ul{0},\ul{e}_{z}).
\end{equation}
Then the constraint on the optimal tangent (\ref{constraint}) follows from fundamental identity
\begin{equation} \label{fundamental}
(\tilde{\gamma}^{\ul{c}}_{g h_{\alpha}}(\cdot)) = \tilde{\gamma}^{\ul{Z}_{\alpha}\ul{c}}_g(\cdot) \, h_{\alpha},
\end{equation}
which holds\footnote{Eq.~\!(\ref{fundamental}) follows from (\ref{niceID}) in App.\!~\ref{app:B}, by setting \mbox{$\ul{Q}=\ul{Z}_{\alpha}$.}} for all $h_{\alpha}$.
We apply this identity (\ref{fundamental}) to the right-hand side of (\ref{projectedCurveFitConstraint}) and use the definition of $\odot$  defined in (\ref{odot}) yielding:
\begin{equation} \label{projectedCurveFitConstraint2}
\begin{array}{rcl}
\tilde{\gamma}_{g}^{\ul{c}^{*}(g)}(t) \odot (\ul{0},\ul{e}_{z})
&=&\tilde{\gamma}_{g}^{\ul{Z}_{\alpha} \ul{c}^{*}(g h_{\alpha})}(t) h_{\alpha} \odot (\ul{0},\ul{e}_{z})  \vspace{-1mm} \\
&\Updownarrow&  \vspace{-1mm} \\
\tilde{\gamma}_{g}^{\ul{c}^{*}(g)}(t) \odot (\ul{0},\ul{e}_{z})
&=&\tilde{\gamma}_{g}^{\ul{Z}_{\alpha} \ul{c}^{*}(g h_{\alpha})}(t) \odot (\ul{0},\ul{e}_{z}) \vspace{-1mm} \\
&\Updownarrow&  \vspace{-1mm} \\
 \ul{c}^{*}(g)&=& \ul{Z}_{\alpha} \ul{c}^{*}(g h_{\alpha}) .
\end{array}
\end{equation}
Finally our constraint  (\ref{constraint}) follows from $\ul{Z}_{\alpha}^T=\ul{Z}_{\alpha}^{-1}$.
$\hfill \Box$




\subsection{Neighboring Exponential Curves in $SE(3)$ \label{ch:FamilyExponentialCurvesSE3}}

\begin{figure}
\hspace{-0.3cm}\mbox{}
\includegraphics[width=1.0\hsize]{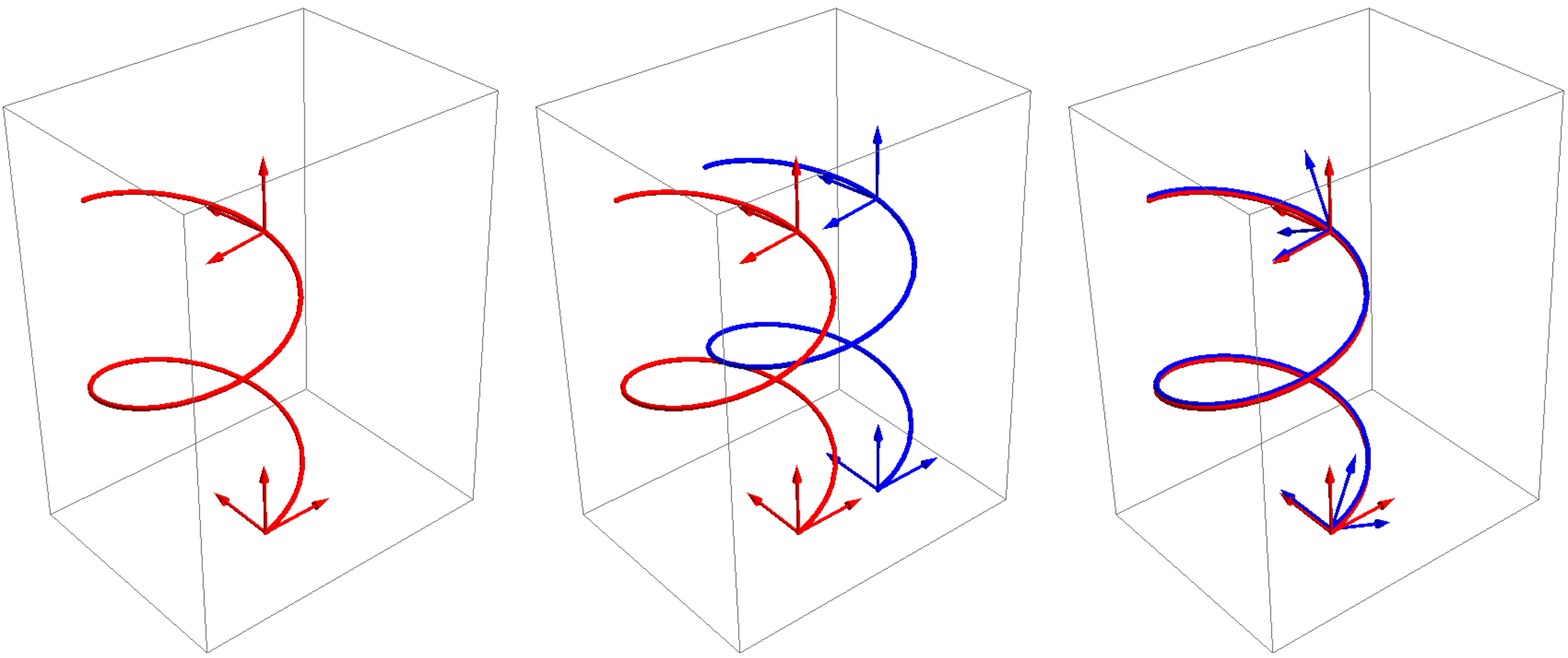}
\caption{Illustration of the family of curves $\tilde{\gamma}_{h,g}^{\ul{c}}$ in $SE(3)$. Left: The (spatially projected) exp-curve $t \mapsto P_{\R^{3}}\tilde{\gamma}^{\ul{c}}_{g}(t)$, with $g=(\ul{x},\ul{R}_{\ul{n}})$ in red. The frames indicate the rotation part $P_{SO(3)}\tilde{\gamma}^{\ul{c}}_{g}(t)$, which for clarity we depicted only at two time instances $t$.
Middle: neighboring exp-curve $t \mapsto \tilde{\gamma}^{\ul{c}}_{g,h}(t)$ with $h=(\ul{x}',\ul{R}_{\ul{n}})$, $\ul{x}\neq \ul{x}'$ in blue, i.e. neighboring exp-curve departing with same orientation and different position.
Right: exp-curve $t \mapsto \tilde{\gamma}^{\ul{c}}_{g,h}(t)$ with $h=(\ul{x},\ul{R}_{\ul{n}'})$, $\ul{n}'\neq \ul{n}$, i.e. the neighboring exp-curve departing with same position and different orientation.
\label{fig:pt}}
\end{figure}

Here we generalize the concept of family of neighboring exponential curves (\ref{paralleltransportRd}) in the $\R^d$-case, and Definition~\ref{def:SE2neighbors} in the $SE(2)$-case, to the $SE(3)$-case.
\begin{definition} \label{def:neigbors}
Given a fixed reference point $g \in SE(3)$ and velocity component $\ul{c}=\ul{c}(g)=(\ul{c}^{(1)}(g),\ul{c}^{(2)}(g)) \in \R^6$, we define the family $\{\tilde{\gamma}_{h,g}^{\ul{c}}(\cdot)\}$ of neighboring exponential curves by
\begin{equation} \label{paralleltransportSE3version2}
t \mapsto \tilde{\gamma}_{h,g}^{\ul{c}}(t):=\tilde{\gamma}_{h}^{\tilde{\ul{R}}_{h^{-1}g} \ul{c}}(t),
\end{equation}
with rotation matrix $\tilde{\ul{R}}_{h^{-1}g} \in SO(6)$ defined by
\begin{equation}\label{RtildeSE3}
              \tilde{\ul{R}}_{h^{-1}g} := \begin{pmatrix}
              (\ul{R}')^T \ul{R} && \ul{0} \\
              \ul{0} && (\ul{R}')^T \ul{R}
              \end{pmatrix} ,
\end{equation}
for all $g=(\ul{x},\ul{R}), h=(\ul{x}',\ul{R}') \in SE(3)$.
\end{definition}
The next lemma motivates our specific choice of neighboring exponential curves. The geometric idea is visualized in Fig.~\!\ref{fig:pt} and is in accordance with Fig.~\ref{fig:gonz} on the $SE(2)$ case.

\begin{lemma}\label{lemma:SE3}
Exponential curve $\tilde{\gamma}_{h,g}^{\ul{c}}$ departing from $h=(\ul{x}',\ul{R}')\in SE(3)$ given by (\ref{paralleltransportSE3version2}) has the same spatial and rotational velocity as exponential curve $\tilde{\gamma}^{\ul{c}}_g$ departing from \mbox{$g=(\ul{x},\ul{R}) \in SE(3)$}.

On the Lie algebra level; we have that the initial velocity component vectors of the curves $\tilde{\gamma}^{\ul{c}}_g$ and $\tilde{\gamma}_{h,g}^{\ul{c}}$ relate via $\ul{c} \mapsto \tilde{\ul{R}}_{h^{-1}g} \ul{c}$.

On the Lie group level; we have that the curves themselves $\tilde{\gamma}^{\ul{c}}_g(\cdot)=(\ul{x}_{g}(\cdot),\ul{R}_{g}(\cdot))$, $\tilde{\gamma}_{h,g}^{\ul{c}}(\cdot)=
(\ul{x}_{h}(\cdot),\ul{R}_{h}(\cdot))$ relate via
\begin{equation} \label{3star}
\begin{array}{ll}
\ul{x}_{h}(t) &=\ul{x}_{g}(t)-\ul{x}+\ul{x}', \\
\ul{R}_{h}(t) &=\ul{R}_{g}(t) \ul{R}^{-1}\ul{R}'.
\end{array}
\end{equation}
\end{lemma}
\textbf{Proof }
See Appendix~\ref{app:B}.
\begin{remark}
Lemma~\ref{lemma:SE3} extends Lemma~\ref{lemma:SE2} to the $SE(3)$ case. When projecting the curves $\tilde{\gamma}_{g}^{\ul{c}}$ and $\tilde{\gamma}_{h,g}^{\ul{c}}$ into the quotient, one has that curves
$\tilde{\gamma}_{g}^{\ul{c}} \odot (\ul{0},\ul{a})$, and $\tilde{\gamma}_{h,g}^{\ul{c}} \odot (\ul{0},\ul{a})$ in $\R^{3}\rtimes S^2$ carry the same spatial and angular velocity.
\end{remark}
\begin{remark}\label{rem:2}
In order to construct the family of neighboring exponential curves in $SE(3)$ one applies the transformation
$\ul{c} \mapsto \tilde{\ul{R}}_{h^{-1}g} \ul{c}$
in the Lie algebra. Such a transformation preserves the
left-invariant metric:
\begin{equation} \label{unitary}
1=\left.\gothic{G}\right|_{\tilde{\gamma}_{g}^{\ul{c}}(t)}(\dot{\tilde{\gamma}}_{g}^{\ul{c}}(t),
\dot{\tilde{\gamma}}_{g}^{\ul{c}}(t))
=\left.\gothic{G}\right|_{\tilde{\gamma}_{h,g}^{\ul{c}}(t)}(\dot{\tilde{\gamma}}_{h,g}^{\ul{c}}(t),
\dot{\tilde{\gamma}}_{h,g}^{\ul{c}}(t)),
\end{equation}
for all $h \in SE(3)$ and all $t \in \R$.
For further differential geometrical details see Appendix~\ref{app:B}.
\end{remark}

\subsection{Exponential Curve Fits in $SE(3)$ of the 1st Order \label{ch:1}}

Now let us generalize the first-order exponential curve fits of Theorem~\ref{th:0} to the setting of $\R^{3}\rtimes S^{2}$.
Here we first consider the following optimization problem on $SE(3)$ 
(generalizing (\ref{optR2})):
\begin{equation} \label{optSE3}
\boxed{
\ul{c}^{*}(g)  = \!\!
\argmin_{\scriptsize \begin{array}{c}
\ul{c} \in \R^{6}, \\ \|\ul{c}\|_{\mu}=1, \\ c^6=0 \end{array} }
\int \limits_{SE(3)}\!\!\!\! \tilde{G}_{\boldsymbol{\rho}}(h^{-1}\!g)
\left|\left. \frac{d}{dt}\! \tilde{V}(\tilde{\gamma}_{h,g}^{\ul{c}}(t))\right|^2_{t=0} \right|
{\rm d}\overline{\mu}(h),
}
\end{equation}
 Recall that $\|\cdot\|_{\mu}$ was defined in (\ref{normmu}), $\tilde{V}$ in (\ref{VEET}) and $\mu$ in (\ref{Haarmeasure}). The reason for including the condition $c^6=0$ will become clear after defining the structure matrix.

\subsubsection{The Structure Tensor on $SE(3)$ \label{ch:3-2b}}

We define structure matrices $\mathbf{S}^{\ul{s},\boldsymbol{\rho}}$ of $\tilde{U}$ by
\begin{equation} \label{structSE3}
\begin{array}{l}
(\mathbf{S}^{\ul{s},\boldsymbol{\rho}} (\tilde{U}))(g)=
\int \limits_{SE(3)}\!\! \tilde{G}_{\boldsymbol{\rho}}(h^{-1}\!g) \cdot \\
 \qquad \tilde{\ul{R}}_{h^{-1}g}^{T} \nabla^{\ul{s}} \tilde{U}(h)(\nabla^{\ul{s}} \tilde{U}(h))^{T} \tilde{\ul{R}}_{h^{-1}g} {\rm d}\overline{\mu}(h),
\end{array}
\end{equation}
where we use matrix $\tilde{\ul{R}}_{h^{-1}g}$ defined in Eq. (\ref{RtildeSE3}). Again we use short notation $\mathbf{S}^{\ul{s},\boldsymbol{\rho}} :=\mathbf{S}^{\ul{s},\boldsymbol{\rho}} (\tilde{U})$.
\begin{remark}
By construction (\ref{tildeU}) and (\ref{LINV}) we have
\[(\mathcal{A}_{6} \tilde{U})(g)=
\lim \limits_{h \downarrow 0} \frac{\tilde{U}(g \, e^{h A_{6}}) - \tilde{U}(g)}{h}=0,
 \]
so the null space of our structure-matrix includes
\begin{equation} \label{nullspace}
\mathcal{N}:=\textrm{span}\{(0,0,0,0,0,1)^{T}\}.
\end{equation}
\end{remark}
\begin{remark} \label{rem:ass}
We assume that $\ul{s}=(s_{p}, s_{o})$ and function $\tilde{U}$ are chosen in such a way that the null space of the structure matrix is precisely equal to $\mathcal{N}$
(and not larger).
\end{remark}
Due to the assumption in Remark~\ref{rem:ass} we need to impose the condition
\begin{equation} \label{c6iszero}
c^{6}=0 \ \desda \
\dot{\tilde{\gamma}}^{\ul{c}}_{g}(0) \, \cap \, \mathcal{N} = \emptyset
\end{equation}
in our exponential curve optimization to avoid
non-uniqueness of solutions. To clarify this, we note that the optimization functional in (\ref{optSE3}) can be rewritten as
\[
\mathcal{E}(\ul{c}):=\ul{c}^{T} \ul{M}_{\mu^2} \mathbf{S}^{\ul{s},\boldsymbol{\rho}}(g) \ul{M}_{\mu^2} \ul{c},
\]
as we will show in the
next theorem where we solve the optimization problem for first-order exponential curve fits. Indeed, for uniqueness we need (\ref{c6iszero}) as otherwise we would have
$\mathcal{E}(\ul{c}+\ul{M}_{\mu^{-2}}\ul{c}_0)=\mathcal{E}(\ul{c})$ for all $\ul{c}_0 \in \mathcal{N}$.

\begin{theorem}[First Order Fit via Structure Tensor] \label{th:2}
The normalized eigenvector $\ul{M}_{\mu}\ul{c}^{*}(g)$ with smallest non-zero eigenvalue of the rescaled structure matrix \\
$\ul{M}_{\mu}\ul{S}^{\ul{s},\boldsymbol{\rho}}(g)\ul{M}_{\mu}$
provides the solution $\ul{c}^{*}(g)$ to optimization problem (\ref{optSE3}).
\end{theorem}
\textbf{Proof } \
All steps (except for the final step of this proof, where the additional constraint $c^6=0$ enters the problem) are analogous to the proof of the first order method in the SE(2) case: the proof of Theorem \ref{th:0}. We will now shortly repeat these first steps. First we rewrite the time derivative as a directional derivative which is then rewritten to the gradient
\begin{equation}
  \begin{array}{ll}
  \left| \!\left.\frac{d}{dt} \! \left( \tilde{V}(\tilde{\gamma}^{\ul{c}}_{h,g}(t)) \right)\right|_{t=0} \right|^2 \!\!\!   &= \!
  \left| \langle {\rm d}\tilde{V}|_{
  \tilde{\gamma}^{\ul{c}}_{h,g}(0)}, \dot{\tilde{\gamma}}^{\ul{c}}_{h,g}(0)\rangle \right|^2 \\
    &=\! \left| \langle {\rm d}\tilde{V}|_{h}, \tilde{\mathbf{R}}_{h^{-1}g} \ul{c}  \rangle \right|^2 \\
     &=\! \left| \left.\gothic{G}_{\mu}\right|_{h}(\nabla \tilde{V}(h), \tilde{\mathbf{R}}_{h^{-1}g} \ul{c} ) \right|^{2} \!.
  \end{array}
\end{equation}
We then put this result in matrix-vector form:
\begin{equation}
  \begin{array}{l}
   \left| \left.\gothic{G}_{\mu}\right|_{h}(\nabla \tilde{V}(h), \tilde{\mathbf{R}}_{h^{-1}g} \ul{c} ) \right|^{2} \\
          = \ul{c}^{T} \ul{M}_{\mu^{2}} \tilde{\mathbf{R}}_{h^{-1}g}^{T} (\nabla \tilde{V}(h)) (\nabla \tilde{V}(h))^{T} \tilde{\mathbf{R}}_{h^{-1}g} \ul{M}_{\mu^{2}} \ul{c}.
  \end{array}
\end{equation}
This again yields the following optimization functional
\begin{equation}
\begin{array}{ll}
\mathcal{E}(\ul{c}) &= \!\int_{SE(3)} \tilde{G}_{\boldsymbol{\rho}}(h^{-1}g) \left| \left. \frac{d}{dt} \tilde{V}(\tilde{\gamma}^{\ul{c}}_{h,g}(t)) \right|_{t=0} \right|^2  \!{\rm d}\overline{\mu}(h)    \\
 &= \ul{c}^{T}\ul{M}_{\mu^{2}}  \mathbf{S}^{\ul{s},\boldsymbol{\rho}}(g)  \ul{M}_{\mu^{2}} \ul{c}.
\end{array}
\end{equation}
So, just as in the $SE(2)$-case we have the following Euler-Lagrange equations:
\begin{equation}\label{EulerLagrangreSE3}
	\begin{array}{rcl}
\mathbf{M}_{\mu^{2}} \ul{S}^{\ul{s},\boldsymbol{\rho}}(g) \mathbf{M}_{\mu^{2}} \ul{c}^{*}(g) &=& \lambda_{1}\, \mathbf{M}_{\mu^2} \ul{c}^{*}(g)  \vspace{-1mm} \\
 & \Updownarrow &  \vspace{-1mm} \\
\! \mathbf{M}_{\mu} \ul{S}^{\ul{s},\boldsymbol{\rho}}(g) \mathbf{M}_{\mu} (\mathbf{M}_{\mu} \ul{c}^{*}(g)) &=&
\lambda_{1}\, (\mathbf{M}_{\mu} \ul{c}^{*}(g)).
\end{array}
\end{equation}
Again the second equality in (\ref{EulerLagrangreSE3}) follows from the first by multiplication by $\ul{M}_{\mu}^{-1}$.

Finally, the constraint $c^{6}=0$ is included in our optimization problem (\ref{optSE3}) to excluded the null space (\ref{nullspace}) from the optimization, therefore we take the eigenvector with the smallest \emph{non-zero} eigenvalue providing us the final result.
$\hfill \Box$

\subsubsection{Projected Exponential Curve Fits in $\R^{3}\rtimes S^{2}$ \label{ch:3-3}}

In reducing the problem to $\R^{3} \rtimes S^{2}$ we first note that
\begin{equation} \label{STSE3Zalpha}
\begin{array}{l}
\mathbf{S}^{\ul{s},\boldsymbol{\rho}}(g \, h_{\alpha})=
\ul{Z}_{\alpha}^{T} \mathbf{S}^{\ul{s},\boldsymbol{\rho}}(g) \ul{Z}_{\alpha},
\end{array}
\end{equation}
with $\ul{Z}_{\alpha}$ defined in Eq. (\ref{Zalphadef}), and
where we recall $h_{\alpha}=(\ul{0},\ul{R}_{\ul{e}_{z},\alpha})$.

In the following theorem we summarize the well-posedness of our projected curve fits on data $U:\R^{3}\rtimes S^2 \to \R$ and use the quotient structure to simplify the structure tensor.

\begin{theorem}[First Order Fit and Quotient Structure] \label{th:3b}
Let $g=(\ul{x},\ul{R}_{\ul{n}})$ and $h=(\ul{x}',\ul{R}_{\ul{n}'})$ where $\ul{R}_{\ul{n}}$ and $\ul{R}_{\ul{n}'}$ denote any rotation which maps $\ul{e}_z$ onto $\ul{n}$ and $\ul{n}'$ respectively. Then, the structure tensor defined by (\ref{structSE3}) can be expressed as
\begin{equation}\label{structR3XS2}
\begin{array}{l}
 \mathbf{S}^{\ul{s},\boldsymbol{\rho}}(g)= 2\pi
\int \limits_{\R^{3}} \int \limits_{S^{2}}  G_{s_{p}}^{\R^{3}}(\ul{x}\!-\!\ul{x}') \; G_{s_0}^{S^{2}}(\ul{R}_{\ul{n}'}^{T}\ul{n}) \\[6pt]
\tilde{\mathbf{R}}_{h^{-1}g}^T
\nabla \tilde{V}(h) \, (\nabla \tilde{V}(h))^{T}
\tilde{\mathbf{R}}_{h^{-1}g}
{\rm d}\sigma(\ul{n}') {\rm d}\ul{x}'.
\end{array}
\end{equation}
The normalized eigenvector $\ul{M}_{\mu}\ul{c}^{*}(\ul{x},\ul{R}_{\ul{n}})$ with smallest non-zero eigenvalue of the rescaled structure matrix $\ul{M}_{\mu}\ul{S}^{\ul{s},\boldsymbol{\rho}}(g)\ul{M}_{\mu}$provides the solution of
 (\ref{optSE3}) and defines a projected curve fit in $\R^{3}\rtimes S^{2}$:
\begin{equation} \label{optcurve}
\gamma_{(\ul{x},\ul{n})}^*(t)=(\,\tilde{\gamma}^{\ul{c}^{*}(\ul{x}, \ul{R}_{\ul{n}})}_{(\ul{x}, \ul{R}_{\ul{n}})}(t) \,)\odot (\ul{0},\ul{e}_{z}),
\end{equation}
which is independent of the choice of $\ul{R}_{\ul{n}'}$ and $\ul{R}_{\ul{n}}$.
\end{theorem}
\textbf{Proof }
The proof consists of two parts. First we prove that (\ref{structR3XS2}) follows from the structure tensor defined in (\ref{structSE3}). Then we use Lemma \ref{lemmaProjectedCurveFit} to prove that our projected exponential curve fit (\ref{optcurve}) is well-defined. For both we use Theorem~\ref{th:2} as our venture point.

For the first part of the proof we note that the integrand in the structure tensor definition Eq.~(\ref{structSE3}) is invariant under $h \mapsto h h_{\alpha}=h (\ul{0},\ul{R}_{\ul{e}_{z},\alpha})$ on the integration variable. To show this we first note that $\ul{Z}_{\alpha}$ defined in (\ref{Zalphadef}), satisfies
\mbox{$\ul{Z}_{\alpha}(\ul{Z}_{\alpha})^T=I$}. Furthermore, we have
\[
\begin{array}{l}
\nabla \tilde{V}(h h_{\alpha})\equiv
\ul{Z}_{\alpha}^{T} \nabla \tilde{V}(h), \
\tilde{\ul{R}}^{T}_{(h h_{\alpha})^{-1} g}= \tilde{\ul{R}}^{T}_{h^{-1} g} \ul{Z}_{\alpha}.
\end{array}
\]
and $\tilde{G}_{\boldsymbol{\rho}}(h h_{\alpha})=\tilde{G}_{\boldsymbol{\rho}}(h)$. Therefore integration over third Euler-angle $\alpha$ is no longer needed in the definition of the structure tensor (\ref{structSE3}) as it just produces a constant $2\pi$ factor.

For the second part we apply Lemma \ref{lemmaProjectedCurveFit} and thereby it remains to be shown that condition $\ul{c}^{*}(g h_{\alpha})= \ul{Z}_{\alpha}^T \ul{c}^{*}(g)$ is satisfied.  This directly follows from (\ref{STSE3Zalpha}):
\begin{equation}\label{ZalphaOptimalTangent}
  \begin{array}{rcl}
  \mathbf{S}^{\ul{s},\boldsymbol{\rho}}(g \, h_{\alpha}) \ul{c}^* (g \, h_{\alpha}) &=& \lambda_1 \ul{c}^* (g \, h_{\alpha})    \vspace{-1mm} \\
   &\Updownarrow& \vspace{-1mm} \\
  \ul{Z}_{\alpha}^{T} \mathbf{S}^{\ul{s},\boldsymbol{\rho}}(g) \ul{Z}_{\alpha}  \ul{c}^* (g \, h_{\alpha}) &=& \lambda_1 \ul{c}^* (g \, h_{\alpha})  \vspace{-1mm} \\
   &\Updownarrow& \vspace{-1mm} \\
   \mathbf{S}^{\ul{s},\boldsymbol{\rho}}(g)  \left( \ul{Z}_{\alpha}  \ul{c}^* (g \, h_{\alpha}) \right) &=& \lambda_1 \left( \ul{Z}_{\alpha} \ul{c}^* (g \, h_{\alpha}) \right)  \vspace{-1mm} \\
   &\Updownarrow& \vspace{-1mm} \\
  \ul{Z}_{\alpha}  \ul{c}^* (g \, h_{\alpha}) &=& \ul{c}^* (g ),
  \end{array}
\end{equation}
which shows our condition. $\hfill \Box$


\subsubsection{Torsion-free Exponential Curve Fits of the 1st Order via a Two-fold Approach \label{ch:3-4}}

Theorem~\ref{th:2} provides us exponential curve fits that possibly carry torsion. From Eq.\!~(\ref{torsion1}) we deduce that the torsion norm of such an exponential curve fit is given by $|\tau|=\frac{1}{\|\ul{c}^{(1)}\|} (c^{1}c^{4}+c^{2}c^5 +c^{3}c^6)|\kappa|$. Together with the fact that we exclude the null space
$\mathcal{N}$ from our optimization domain by including constraint $c^6=0$, this results in insisting on zero torsion along horizontal exponential curves where $c^{1}=c^{2}=0$. Along other exponential curves torsion appears if $c^{1}c^{4}+c^{2}c^5\neq 0$.

Now the problem is that insisting, a priori, on zero torsion for horizontal curves while allowing non-zero torsion for other curves is undesirable. On top of this, torsion is a higher order less-stable feature than curvature. Therefore we would like to exclude it altogether from
our exponential curve fits presented in Theorem~\ref{th:2} and Theorem~\ref{th:3b}, by a different theory and algorithm. The results of the algorithm show that even if structures do have torsion, the local exponential curve fits do not need to carry torsion in order to achieve good results in the local frame adaptation, see e.g.~Fig.~\!\ref{fig:examples3D}.

The constraint of zero torsion forces us to split our exponential curve fit into a two-fold algorithm: \\[6 pt]
 \textbf{Step 1} Estimate at $g \in SE(3)$ the spatial velocity part $\ul{c}^{(1)}(g)$ from the spatial structure tensor.\\[6 pt]
 \textbf{Step 2} Move to a different location $g_{new} \in SE(3)$ where a horizontal exponential curve fit makes sense and then estimate the angular velocity $\ul{c}^{(2)}$ from the rotation part of the structure tensor over there.\\[6 pt]
This forced splitting is a consequence of the next lemma.
\begin{lemma}
Consider the class of exponential curves with nonzero spatial velocity $\ul{c}^{(1)}\neq \ul{0}$ such that their spatial projections do not have torsion.
Within this class the constraint $c^{6}=0$ does not impose constraints on curvature if and only if the exponential curve is horizontal.
\end{lemma}
\textbf{Proof }
For a horizontal curve $\tilde{\gamma}^{\ul{c}}_{g}(t)$ we have
$\chi=0 \desda c^{1}=c^{2}=0$ and indeed $|\tau| = \frac{\ul{c}^{(1)} \cdot \ul{c}^{(2)} |\kappa|}{\|\ul{c}^{(1)}\|}=c^{6}|\kappa|=0$
and we see that constraints $c^{6}=0$ and $|\tau|=0$ reduce to only one constraint.
The curvature magnitude stays constant along the exponential curve and the curvature vector at $t=0$, recall Eq.~\!(\ref{curvature}), is in this case given by
\[
\KKK(0)= \frac{1}{|c^{3}|}
\left(
\begin{array}{c}
c^{5}c^{3} \\
-c^{4} c^{3} \\
0
\end{array}
\right),
\]
which can be any vector orthogonal to spatial velocity $\ul{c}^{(1)}=(0,0,c^{3})^T$. Now let us check whether the condition is necessary. Suppose $t \mapsto \tilde{\gamma}_{g}^{\ul{c}}(t)$ is \emph{not} horizontal, and suppose it is torsion free with $c^6=0$. Then we have
$c^{1}c^{4}+c^{2}c^{5}=0$,
as a result the initial curvature
\[
\KKK(0)= \frac{1}{|c^{3}|}
\left(
\begin{array}{c}
c^{5}c^{3} \\
-c^{4} c^{3} \\
c^4c^2-c^1 c^{5}
\end{array}
\right),
\]
is both orthogonal to vector $\ul{c}^{(1)}=(c^1,c^2,c^3)^T$ and orthogonal to $(-c^2,c^1,0)^T$,
and thereby constrained to a one dimensional subspace. $\hfill \Box$ \\ \\
From these observations we draw the following conclusion for our exponential curve fit algorithms.
\ \\

\textbf{Conclusion:} In order to allow for all possible curvatures in our torsion-free exponential curve fits we must relocate
the exponential curve optimization at $g \in SE(3)$ in $\tilde{U}:SE(3) \to \R$ to a position $g_{new} \in SE(3)$ where a horizontal exponential curve can be expected. Subsequently, we can use Lemma~\ref{lemma:SE3} to transport the horizontal and torsion-free curve through $g_{new}$, back to
a torsion-free exponential curve through $g$.
\ \\

This conclusion is the central idea behind our following two-fold algorithm for exponential curve fits.

\subsubsection*{Algorithm Two-fold Approach:}
The algorithm follows the subsequent steps:

\textbf{Step 1a: }Initialization. Compute structure tensor \\
$\mathbf{S}^{\ul{s},\boldsymbol{\rho}}(g)$ from input image
$U:\R^{3} \times S^{2} \to \R^{+}$ via Eq.\!~(\ref{structR3XS2}).

\textbf{Step 1b: }Find the optimal spatial velocity:
\begin{equation}\label{step1}
\!\! \ul{c}^{(\;\!\!1\;\!\!)}(g) = 
\!\!\!\!\!
 \argmin_{\scriptsize \begin{array}{c}
\ul{c}^{(\;\!\!1\;\!\!)} \in \R^{3},\\
 \|\ul{c}^{(\;\!\!1\;\!\!)}\|=\mu^{-1}
 \end{array}}
 \!\!\!\! \left\{ \!\!
 \begin{pmatrix}
 \ul{c}^{(\;\!\!1\;\!\!)} \\
 \ul{0}
 \end{pmatrix}
 ^{\!\!T} \!\! \!
 \ul{M}_{\mu^2}
  \mathbf{S}^{\ul{s},\boldsymbol{\rho}}(g)
  \ul{M}_{\mu^2} \!
  \begin{pmatrix}
 \ul{c}^{(\;\!\!1\;\!\!)} \\
 \ul{0}
 \end{pmatrix}
 \!\! \right\} \!,
\end{equation}
for $g=(\ul{x},\ul{R}_{\ul{n}}$), which boils down to finding the eigenvector with minimal eigenvalue
of the $3\times 3$ spatial sub-matrix of the structure tensor (\ref{structSE3}).

\textbf{Step 2a: }Given $\ul{c}^{(1)}(g)$ we aim for an auxiliary set of coefficients, where we also take into account rotational velocity. To achieve this in a stable way we move to a different location in the group:
\begin{equation} \label{nnew}
g_{new}=(\ul{x},\ul{R}_{\ul{n}_{new}}), \ \ul{n}_{new}=\ul{R}_{\ul{n}} \ul{c}^{(1)},
\end{equation}
and apply the transport of Lemma~\ref{lemma:SE3} afterwards.
At $g_{new}$, we enforce horizontality, see Remark~\ref{rem:jump} below, and we consider the auxiliary optimization problem
\begin{equation} \label{aux}
\ul{c}_{new}(g_{new}) =
\!\!\!\!\!\! \argmin_{\scriptsize
\begin{array}{c}
\ul{c} \in \R^{6},\\
 \|\ul{c}\|_{\mu}=1, \\
 c^{1}\!=\!c^{2}\!=\!c^{6}\!=\!0
 \end{array}} \!\!\!\!\!\! \left\{
 \ul{c}^{T}
 \ul{M}_{\mu^2}
  \mathbf{S}^{\ul{s},\boldsymbol{\rho}}(g_{new})
  \ul{M}_{\mu^2}
 \ul{c}
 \right\}\! .
\end{equation}
Here zero deviation from horizontality (\ref{dH}) and zero torsion (\ref{torsion1}) is equivalent to the imposed constraint:
\[
\chi=0 \textrm{ and }|\tau|=0 \desda c^{1}=c^{2}=c^{6}=0.
\]

\textbf{Step 2b: }The auxiliary coefficients $\ul{c}_{new}(g_{new})=(0,0,c^{3}(g_{new}),c^{4}(g_{new}),c^{5}(g_{new}),0)^{T}$ of a torsion-free, horizontal exponential curve fit $\tilde{\gamma}_{g_{new}}^{\ul{c}_{new}}$ through $g_{new}$. Now we apply transport (via Lemma~\ref{lemma:SE3}) of this horizontal exponential curve fit towards the corresponding exponential curve through $g$:
\begin{equation} \label{cFinalTwoStep}
\ul{c}^{*}_{final}(g)= \begin{pmatrix}
\ul{R}_{\ul{n}}^{T}\ul{R}_{\ul{n}_{new}} && \ul{0} \\
\ul{0} && \ul{R}_{\ul{n}}^{T} \ul{R}_{\ul{n}_{new}} \end{pmatrix} \ul{c}_{new}(g_{new}).
\end{equation}
This gives the final, torsion-free, exponential curve fit
$t \mapsto \tilde{\gamma}_{g}^{\ul{c}^{*}(g)}(t)$ in $SE(3)$, yielding the final
output projected curve fit
\begin{equation} \label{finalresult}
t \mapsto (\tilde{\gamma}_{g}^{\ul{c}^{*}_{final}(g)}(t)) \odot (\ul{0},\ul{e}_{z}) \in \R^{3} \times S^{2},
\end{equation}
with $g=(\ul{x},\ul{R}_{\ul{n}})$, recall Eq.~\!(\ref{odot}).
\begin{remark}\label{rem:jump}
In \textbf{step 2a} of our algorithm we jump to a new location $g_{new}=(\ul{x},\ul{R}_{\ul{n}_{new}})$ with possibly different orientation $\ul{n}_{new}$
such that spatial tangent vector \[\sum_{i=1}^{3}c^{i}\left.\mathcal{A}_{i}\right|_{(\ul{x},\ul{R}_{\ul{n}_{new}})},
\]
 points in the same direction as $\ul{n}_{new} \in S^2$, recall Eq.~\!(\ref{lifttogroup}), from which it follows
that $\ul{n}_{new}$ is indeed given by (\ref{nnew}). If $\ul{c}^{(1)}=(c^{1},c^{2},c^{3})^T=\ul{a}=(0,0,1)^T$ then $\ul{n}_{new}=\ul{n}$.
\end{remark}
\begin{lemma}
The preceding algorithm is well-defined on the quotient $\R^{3}\rtimes S^{2}=SE(3)/(\{\ul{0}\}\times SO(2))$.
\end{lemma}
\textbf{Proof }
To show that the preceeding algorithm is well-defined on the quotient we need to show that the final result (\ref{finalresult}) is independent on both the choice of  of $\ul{R}_{\ul{n}} \in SO(3)$ s.t. $\ul{R}_{\ul{n}}\ul{e}_{z}=\ul{n}$ and the choice of  $\ul{R}_{\ul{n}_{new}} \in SO(3)$ s.t. $\ul{R}_{\ul{n}_{new}}\ul{e}_{z}=\ul{n}_{new}$.

First, we show independence on the choice of $\ul{R}_{\ul{n}}$. We apply Lemma \ref{lemmaProjectedCurveFit} and thereby it remains to be shown that condition $\ul{c}^{*}_{final}(g h_{\alpha})= \ul{Z}_{\alpha}^T \ul{c}^{*}_{final}(g)$ is satisfied. This follows directly from Eq. (\ref{cFinalTwoStep}) if as long as $\ul{n}_{new}$ found in Step 2a is independent of the choice of $\ul{R}_{\ul{n}}$. This property indeed follows from $\ul{c}^{(1)}(g h_{\alpha}) = \ul{R}_{\ul{e}_{z},\alpha}^T \ul{c}^{(1)}(g)$ which can be proven analogously to (\ref{ZalphaOptimalTangent}). Then we have
\begin{equation}
\begin{array}{rl}
  \ul{n}_{new} (g h_{\alpha} ) &=\ul{R}_{\ul{n}} \ul{R}_{\ul{e}_{z},\alpha} \ul{c}^{(1)} (g h_{\alpha} ) \\
  &= \ul{R}_{\ul{n}} \ul{R}_{\ul{e}_{z},\alpha} \ul{R}_{\ul{e}_{z},\alpha}^T \ul{c}^{(1)} (g )=   \ul{n}_{new} (g).
  \end{array}
\end{equation}
So we conclude that (\ref{finalresult}) is indeed independent on the choice of $\ul{R}_{\ul{n}}$.

\begin{figure}
\includegraphics[width=\hsize]{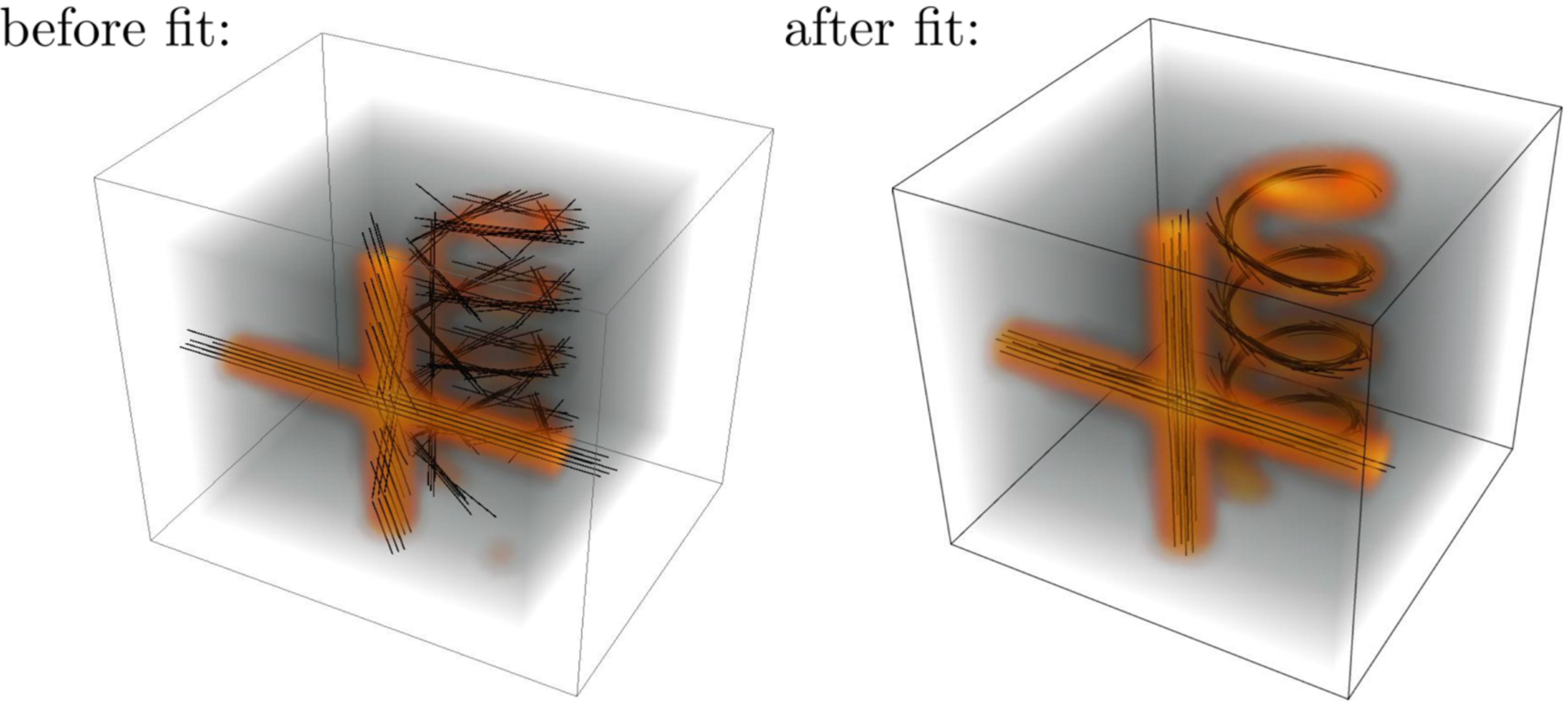}
\caption{Volume rendering of a 3D test-image.
The real part of the orientation score (cf.~Section~\ref{sec:ios})
provides us a density $U$ on $\R^{3}\rtimes S^2$. Left: spatial parts of exponential curves (in black) aligned with spatial generator $\left.\mathcal{A}_{3}\right|_{(\ul{x},\ul{R}_{\ul{n}_{max}(\ul{x})})}$, where $\ul{n}_{max}(\ul{x})= \underset{\ul{n}\in S^{2}}{\textrm{argmax}} |U(\ul{x},\ul{n})|$. Right: spatial parts of our exponential curve fits \mbox{Eq.~\!(\ref{finalresult})} computed via the algorithm in Section~\ref{ch:3-4},
which better follow the curvilinear structures.
 \label{fig:examples3D}}
\end{figure}

Finally, Eq. (\ref{finalresult}) is independent of the choice of $\ul{R}_{\ul{n}_{new}}$.
This follows from $\ul{c}_{new}(g h_{\alpha}) = \ul{Z}_{\alpha}^T \ul{c}_{new}(g)$ in Step 2a. Then $\ul{c}^{*}_{final}$ in  Eq.~\!(\ref{cFinalTwoStep}) is independent of the choice of $\ul{R}_{\ul{n}_{new}}$ because $\ul{Z}_{\alpha}^{T}$ in $\ul{c}_{new} \mapsto \ul{Z}_{\alpha}^{T} \ul{c}_{new}$ is canceled by $\ul{R}_{new} \mapsto \ul{R}_{new} \ul{R}_{\ul{e}_{z},\alpha}$ in Eq.~\!(\ref{cFinalTwoStep}).
$\hfill \Box$ \\
\\
In Fig.~\!\ref{fig:examples3D} we provide an example of spatially projected exponential curve fits in $SE(3)$ via the twofold approach. Here we see that the resulting gauge frames better
follow the curvilinear structures of the data (in comparison to the normal left-invariant frame).

\subsection{Exponential Curve Fits in $SE(3)$ of the 2nd Order \label{ch:SE3-order2}}

In this section we will generalize Theorem~\ref{th:0b} to the case $d=3$, where again
we include the restriction to torsion-free exponential curves.

\subsubsection{The Hessian on SE(3) \label{ch:Hessian}}

For second order curve fits we consider the following optimization problem:
\begin{equation} \label{optHessianSumSE3}
\boxed{
\begin{array}{l}
\ul{c}^{*}(g) =
\argmin_{\ul{c}\in \R^{6}, \|\ul{c}\|_{\mu}=1, c^6=0} 
\left|\;\left.\frac{d^2}{dt^2} \tilde{V}(\gamma_{g}^{\ul{c}}(t)) \right|_{t=0}\;\right|\ ,
\end{array}
}
\end{equation}
with $\tilde{V}=\tilde{G}_{\ul{s}}*\tilde{U}$. Before solving this optimization problem in Theorem \ref{th:6} we first
define the $6\times 6$ non-symmetric Hessian matrix by
\begin{equation} \label{fullhess}
(\ul{H}^{\ul{s}} (\tilde{U}))(g)= [\mathcal{A}_{j}\mathcal{A}_{i} (\tilde{V})](g) \, , \;\;\;
\textrm{with } \tilde{V}=\tilde{G}_{\ul{s}}*\tilde{U}
\end{equation}
and where $i=1,\ldots,6$ denotes the row index, and $j=1,\ldots, 6$ denotes the column index. Again we write $\ul{H}^{\ul{s}}:=\ul{H}^{\ul{s}} (\tilde{U})$.

\begin{theorem}[Second Order Fit via Symmetric Sum Hessian] \label{th:6}
Let $g \in SE(3)$ be such that the symmetrized Hessian matrix
$\frac{1}{2} \ul{M}_{\mu}^{-1}(
\ul{H}^{\ul{s}}(g)+ (\ul{H}^{\ul{s}}(g))^T) \ul{M}_{\mu}^{-1}$
has eigenvalues with the same sign.
Then the normalized eigenvector $\ul{M}_{\mu}\ul{c}^{*}(g)$ with smallest absolute non-zero eigenvalue of the symmetrized Hessian matrix
provides the solution $\ul{c}^{*}(g)$ of optimization problem (\ref{optHessianSumSE3}).

\end{theorem}
\textbf{Proof }
Similar to the proof of Theorem~\ref{th:0b} (only now with summations from 1 to 5). Again we include our additional constraint $c^6=0$ by taking the smallest \emph{non-zero} eigenvalue. $\hfill \Box$

\begin{remark}
The restriction to $g \in SE(3)$ such that the
eigenvalues of the symmetrized Hessian carry the same sign is necessary for
a unique solution of the optimization. Note that in case of our first order
approach via the positive definite structure tensor, no such cumbersome constraints arise.
In case $g \in SE(3)$ is such that the eigenvalues of the symmetrized Hessian have
different sign there are 2 options:
\begin{enumerate}
\item Move towards a neighboring point where the Hessian eigenvalues have the same sign
and apply transport (Lemma~\ref{lemma:SE3}, Fig.~\!\ref{fig:pt}) of the exponential curve fit at the neighboring point.
\item  Take $\ul{c}^{*}(g)$ still as the eigenvector with smallest absolute eigenvalue (representing minimal absolute principal curvature).
though this no longer solves (\ref{optHessianSumSE3}).
\end{enumerate}
\end{remark}

\subsubsection{Torsion-free Exponential Curve Fits of the 2nd Order via a Two-Fold Algorithm \label{ch:2fold2ndorder}}

In order to obtain torsion-free exponential curve fits of the second order via our two-fold algorithm,
we follow the same algorithm as in Subsection~\ref{ch:3-4},
but now with the Hessian field $\ul{H}^{\ul{s}}$ (\ref{fullhess}) instead of the structure tensor field.

\textbf{Step 1a: }Initialization. Compute Hessian \\
$\mathbf{H}^{\ul{s}}(g)$ from input image
$U:\R^{3} \times S^{2} \to \R^{+}$ via Eq.\!~(\ref{fullhess}).

\textbf{Step 1b: }Find the optimal spatial velocity by (\ref{step1}) where we replace $\ul{M}_{\mu^2} \mathbf{S}^{\ul{s},\boldsymbol{\rho}}(g)   \ul{M}_{\mu^2} $ by $\mathbf{H}^{\ul{s}}(g)$.

\textbf{Step 2a: } We again fit a horizontal curve at $g_{new}$ given by (\ref{nnew}). The procedure is done via (\ref{aux}) where we again replace $\ul{M}_{\mu^2} \mathbf{S}^{\ul{s},\boldsymbol{\rho}}(g)   \ul{M}_{\mu^2} $ by $\mathbf{H}^{\ul{s}}(g)$.

\textbf{Step 2b: } Remains unchanged. We again apply Eq.~\!(\ref{cFinalTwoStep}) and Eq.~\!(\ref{finalresult}).

There are some serious computational technicalities in the efficient computation of the entries of the Hessian for discrete input data, but this is outside the scope of this article and will be pursued in future work.

\begin{figure}
\centerline{
\includegraphics[width=0.5\hsize]{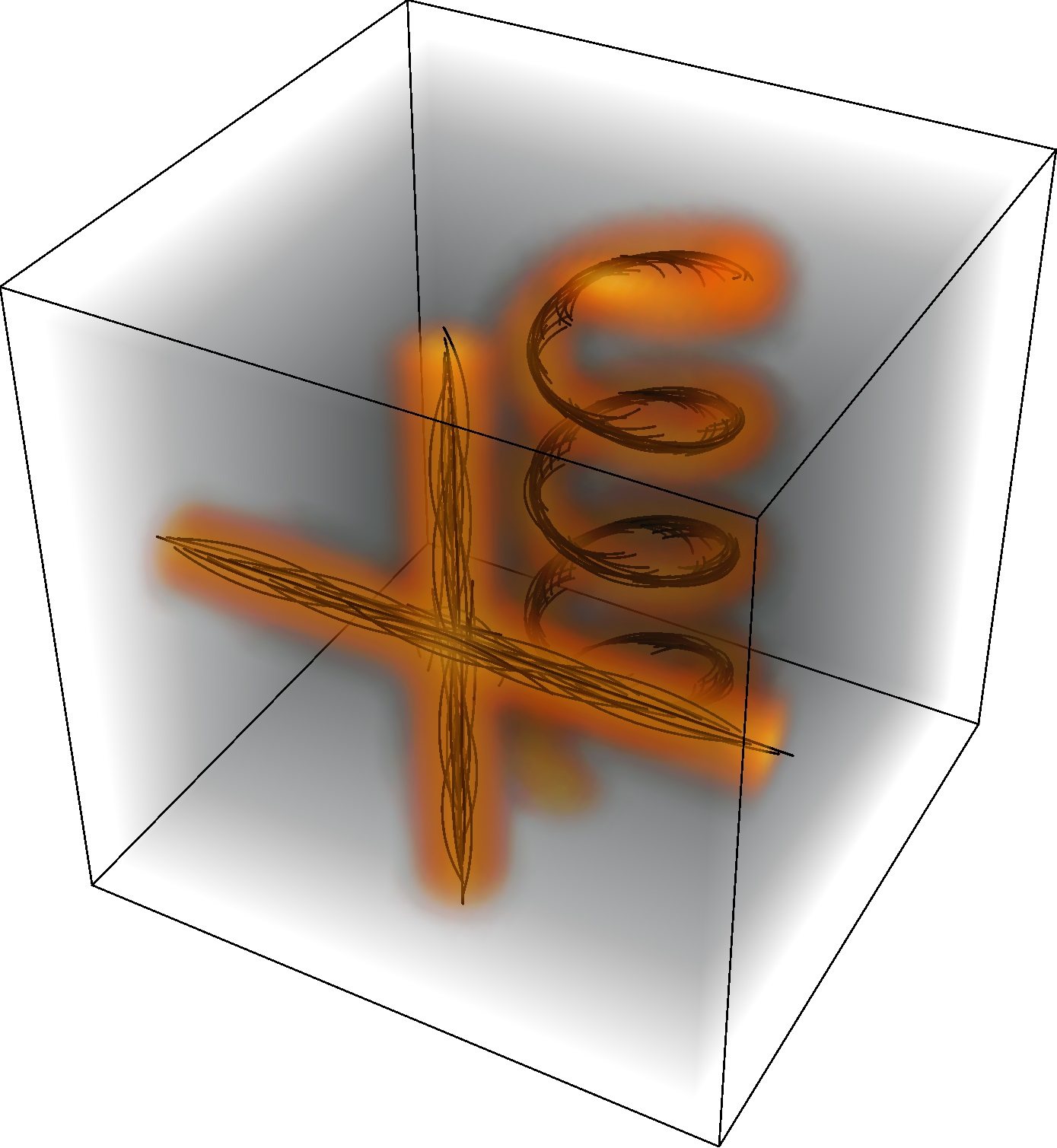}
\includegraphics[width=0.5\hsize]{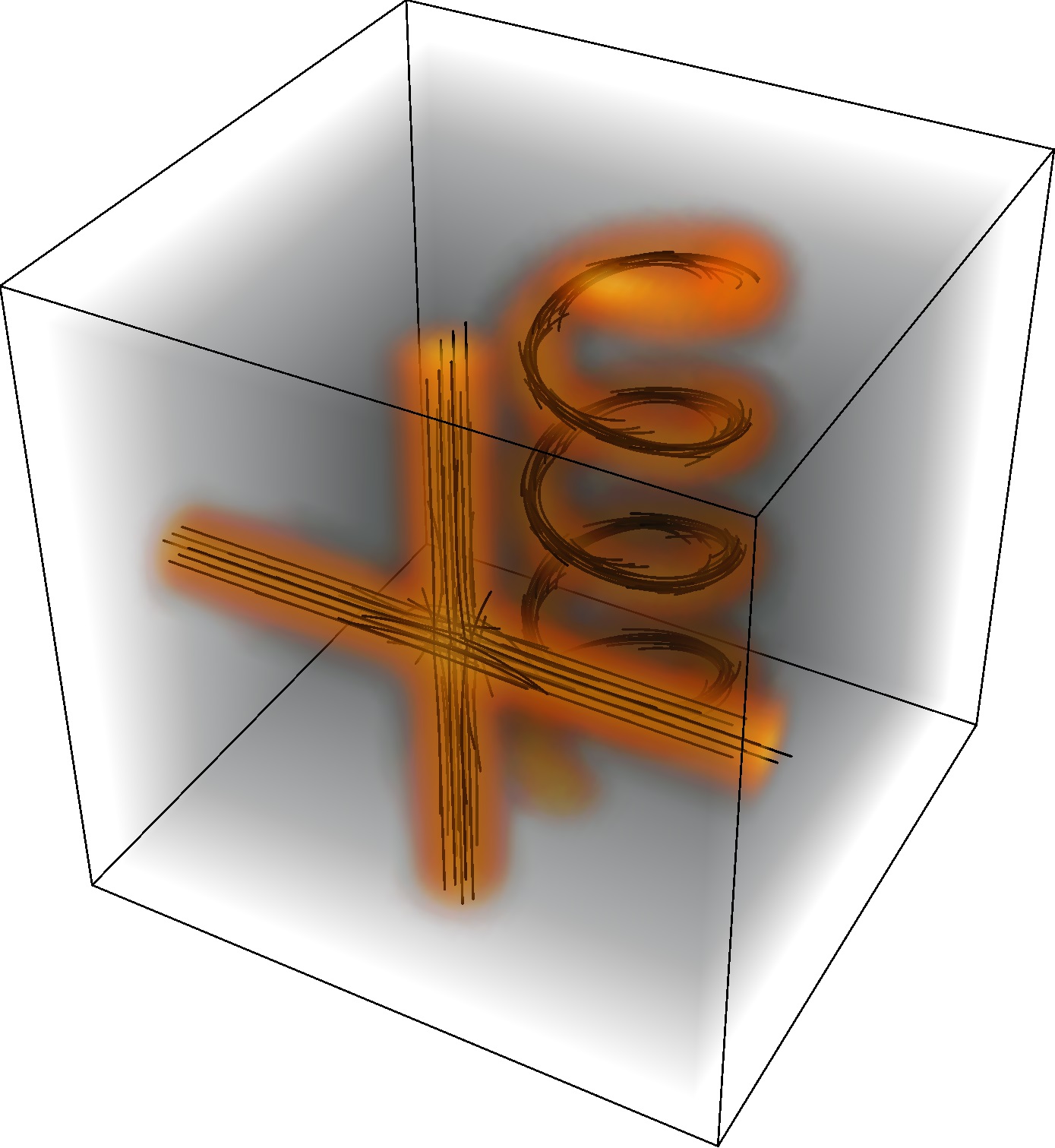}
}
\caption{In black the spatially projected part of exponential curve fits $t \mapsto \gamma^{\ul{c}}_{g}(t)$ \emph{ of the second kind }
(fitted to the real part of the 3D invertible orientation score, for details see~\!Fig.\!~\ref{fig:OS3D})
of the 3D image visualized via volume rendering.
Left: Output of the 2-fold approach outlined in Subsection~\ref{ch:2fold2ndorder}.
Right: Output of the 2-fold approach outlined in
Appendix~\ref{app:C}
with $s_p=\frac{1}{2}$, $s_o=\frac{1}{2} (0.4)^2$, $\mu=10$.
\label{fig:SecondOrderFitSE3}}
\end{figure}

\begin{remark}
In Appendix~\ref{app:C} we propose another two-fold second-order exponential curve fit method. Here one solves a variational problem for exponential curve fits
where exponentials are factorized over respectively spatial and angular part. Empirically, this approach performs good (see e.g. Fig.~\!\ref{fig:SecondOrderFitSE3}).
\end{remark}

\section{Image Analysis Applications \label{ch:app}}
In this section we present examples of applications where the use of gauge frame in $SE(d)$ obtained via exponential curve fits is used for defining data-adaptive left invariant operators. Before presenting the applications, we start by briefly summarizing the invertible orientation score theory in Sec. \ref{sec:ios}.

In case $d=2$ the application presented is the enhancing of the vascular tree structure in 2D retinal images via differential invariants based on gauge frames. This is achieved by extending the classical Frangi vesselness filter \cite{Frangi} to distributions $\tilde{U}$ on $SE(2)$. Gauge frames in $SE(2)$ can also be used in non-linear multiple-scale crossing preserving diffusions as demonstrated in \cite{Sharma}, but we will not discuss this application in this paper.

In case $d=3$ the envisioned applications include blood vessel detection in 3D MR-angiography, e.g. the detection of the Adamkiewicz vessel, relevant for surgery planning. Also in extensions towards fiber-enhancement of diffusion-weighted MRI \cite{DuitsIJCV2010,DuitsJMIV} the non-linear diffusions are of interest. Some preliminary practical results have been conducted on such 3D-datasets \cite{Janssen,DuitsIJCV,Creusen}, but here we shall restrict ourselves to very basic artificial 3D-datasets to show a proof of concept, and leave these three applications for future work.

\subsection{Invertible Orientation Scores} \label{sec:ios}

\begin{figure}
\includegraphics[width=0.99\hsize]{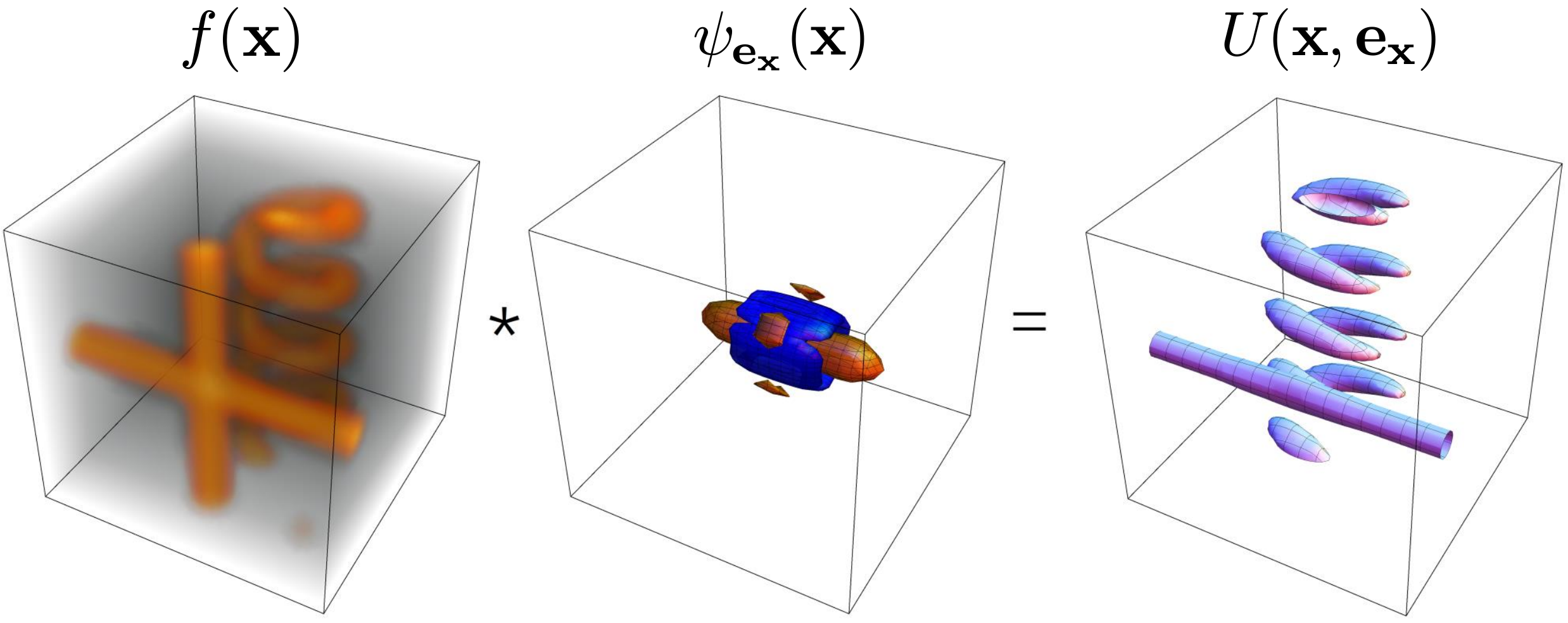}
\includegraphics[width=0.99\hsize]{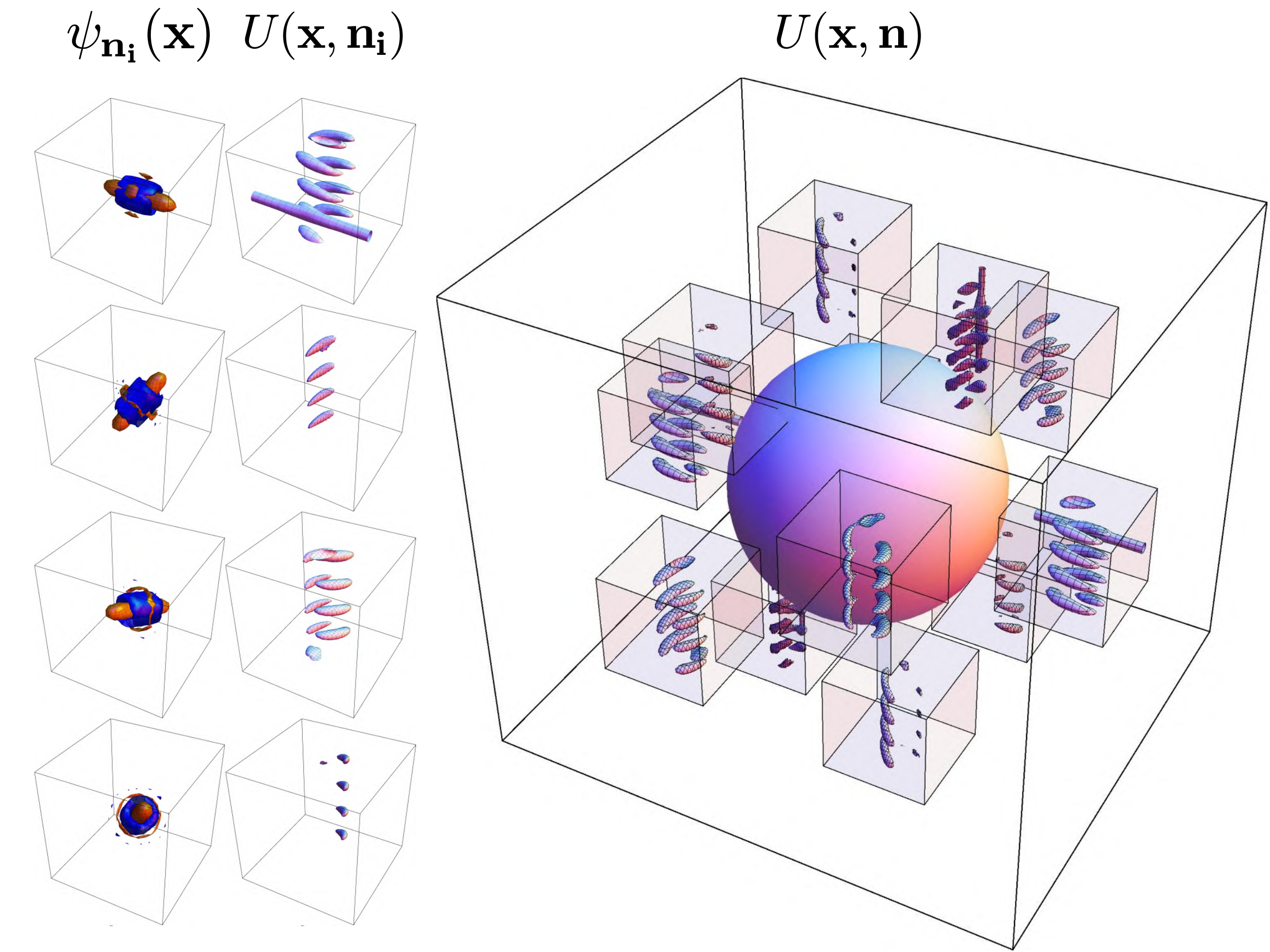}
\caption{Visualization of one iso-level of the real part of an invertible orientation score of a 3D-image created via the cake wavelet in Fig. \ref{fig:cakewavelets}.
\label{fig:OS3D}}
\end{figure}

In the image analysis applications discussed in this section our function $U:\R^{d} \rtimes S^{d-1} \to \R$ is given by the real part of an invertible orientation score: 
\[
U(\ul{x},\ul{n})= \textrm{Re}\{\mathcal{W}_{\psi}f(\ul{x},\ul{R}_{\ul{n}})\},
\]
where $\ul{R}_{\ul{n}}$ is any rotation mapping reference axis $\ul{a}$ onto $\ul{n} \in S^{d-1}$, where $f\in \mathbb{L}_2(\R^d)$ denotes a input image, and where $\psi$ is a so-called 'cake-wavelet' and with
\begin{equation} \label{OS}
\mathcal{W}_{\psi}f(\ul{x},\ul{R}_{\ul{n}})= \int \limits_{\R^d} \overline{\psi(\ul{R}^{-1}_{\ul{n}}(\ul{y}-\ul{x}))}
f(\ul{y})\; {\rm d}\ul{y}.
\end{equation}
For $d >2$ we restrict ourselves to wavelets $\psi$ satisfying
\begin{equation} \label{cigar}
\psi(\ul{R}_{\ul{a},\alpha}^{-1} \ul{x})=\psi(\ul{x}), \textrm{ for all }\ul{x}\in \R^d
\end{equation}
and for all rotations $\ul{R}_{\ul{a},\alpha} \in \textrm{Stab}(\ul{a})$ (for $d=3$ this means for all rotations about axis $\ul{a}$, Eq.~\!(\ref{conventie})). As a result $U$ is well-defined on the left cosets
$\R^{d}\rtimes S^{d-1}=SE(d)/(\{\ul{0}\} \times SO(d-1))$ as the choice of $\ul{R}_{\ul{n}} \in SO(d)$ mapping $\ul{a}$ onto $\ul{n}$ is irrelevant. See Fig.~\!\ref{fig:OS3D} for an example of a 3D orientation score.

\begin{figure}
\centering
\includegraphics[width=\hsize]{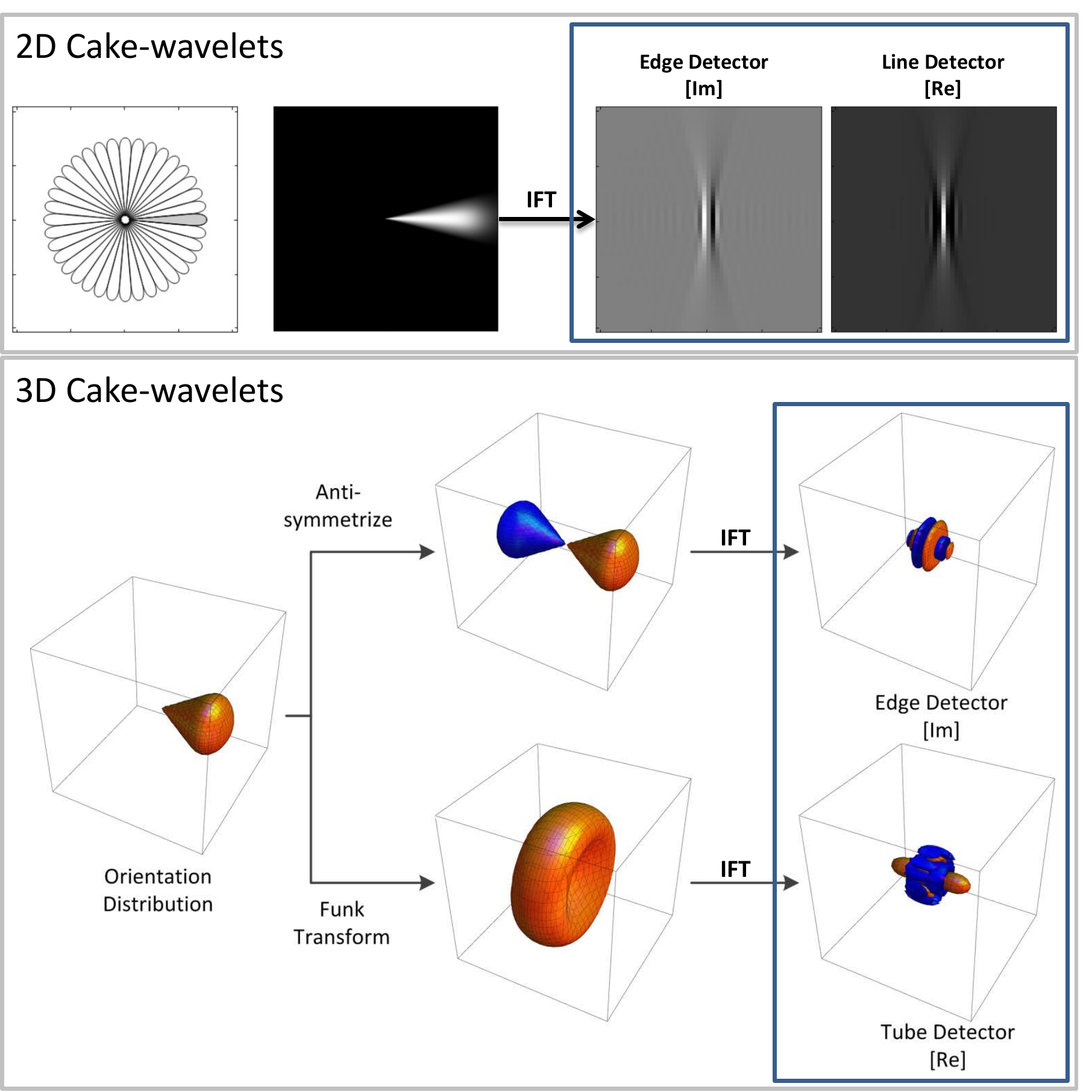}
\caption{Visualization of cake-wavelets in 2D (top) and 3D (bottom). In 2D we fill up the `pie' of frequencies with overlapping ``cake pieces'', and application of an inverse DFT (see \cite{Bekkers}) provides wavelets whose real and imaginary parts are respectively line and edge detectors. In 3D we include anti-symmetrization and the Funk transform  \cite{Descoteaux} on $\mathbb{L}_{2}(S^2)$ to obtain the same, see \cite{Janssen}.  The idea is to redistribute spherical data from orientations towards  circles laying in planes orthogonal to those orientations. Here, we want the real part of our wavelets to be line detectors (and not plate detectors) in the spatial domain. In the figure one positive iso-level is depicted in orange and one negative iso-level is depicted in blue.
\label{fig:cakewavelets}}
\end{figure}

If we restrict to disk-limited images, exact reconstruction is performed via the adjoint:
\begin{equation}\label{OSreconstruction}
\begin{array}{l}
f=\mathcal{W}_{\psi}^{*}\mathcal{W}_{\psi}f=
\mathcal{F}^{-1}\left[\www \mapsto \frac{1}{(2\pi)^{\frac{d}{2}}M_{\psi}(\www)}\right. \\
 \left.\!\int \limits_{SO(d)}\! \mathcal{F}[\mathcal{W}_{\psi}f(\cdot,\ul{R})](\www) \mathcal{F}\psi(\ul{R}^{-1}\www) {\rm d}{\mu}_{SO(d)}(\ul{R})\right].
\end{array}
\end{equation}
if $\psi$ is an admissible wavelet. The condition for admissibility of wavelets $\psi$ are given in \cite{DuitsRPHDThesis}. In this article, the wavelets $\psi$ are given either by the 2D `cake-wavelets' used in \cite{Bekkers,DuitsIJCV} or by their recent 3D-equivalents given in \cite{Janssen}. Detailed formulas and recipes to construct such wavelets efficiently can be found in \cite{Janssen} and in order to provide the global intuitive picture they are depicted in Fig.~\!\ref{fig:cakewavelets}. 

In the subsequent sections we consider two types of operators acting on the invertible orientation scores (recall $\Phi$ in the commuting diagram of Fig.~\!\ref{Fig:Intro}):
\begin{enumerate}
\item for $d=2$, differential invariants on orientation scores based on gauge frames $\{\mathcal{B}_{1},\mathcal{B}_{2},\mathcal{B}_{3}\}$.

\item for $d=2,3$, non-linear adaptive diffusions steered along the gauge frames, i.e.
    \begin{equation} \label{pregaugeflow}
    W(\ul{x},\ul{n},t)=\widetilde{W}(\ul{x},\ul{R}_{\ul{n}},t)=\Phi_{t}(\tilde{U})(\ul{x},\ul{R}_{\ul{n}}),
    \end{equation}
    where $\widetilde{W}(g,t)$, with $t\geq 0$, is the solution of:
    \begin{equation} \label{gaugeflow}
    \left\{
    \begin{array}{rl}
    \frac{\partial \tilde{W}}{\partial t}(g,t)&= \sum \limits_{i=1}^{n_d} D_{ii}\left.(\mathcal{B}_{i})^2\right|_{g} \tilde{W}(g,t),  \\
    \tilde{W}(g,0)&=\tilde{U}(g),
    \end{array}
    \right.
    \end{equation}
    where the gauge frame is induced by an exponential curve fit to data $\tilde{U}$ at location $g \in SE(d)$.
    \end{enumerate}

\subsection{Experiments in $SE(2)$}

We consider the application of enhancing and detecting the vascular tree structure in retinal images. Such image processing task is highly relevant as the retinal vasculature provides non--invasive observation of the vascular system. A variety of  diseases such as glaucoma, age--related macular degeneration, diabetes, hypertension, arteriosclerosis or Alzheimer's affect the vasculature and may cause functional or geometric changes \cite{Ikram2013}. Automated quantification of these defects promises massive screenings for vascular-related diseases on the basis of fast and inexpensive retinal photography. To automatically assess the state of the retinal vascular tree, vessel segmentation are needed. Because retinal images usually suffer from low contrast on small scales, the vasculature in the images needs to be enhanced prior to the segmentation. One well--established approach is the Frangi vesselness filter \cite{Frangi} which is used in robust retinal vessel segmentation methods \cite{Budai2013,Lupascu2010}. However, a drawback of the Frangi filter is that it can not handle crossings or bifurcations that make up an important part of the vascular network. This is precisely where the orientation score framework and the presented locally adaptive frame theory comes into play.

The $SE(2)$-vesselness filter, extending  Frangi vesselness \cite{Frangi} to $SE(2)$ (cf.\!~\cite{Hannink}) and based on the locally adapted frame  $\{\mathcal{B}_{1},\mathcal{B}_{2},\mathcal{B}_{3}\}$ is given by the following left invariant operator:
 \begin{equation} \label{SE2vesselness}
 \begin{array}{l}
 \Phi(\tilde{U})= \left\{
 \begin{array}{lll}
 e^{-\frac{\gothic{R}^2}{2\sigma_{1}^2}}\left(1- e^{-\frac{\gothic{S}}{2\sigma_{2}}} \right) & \textrm{ if }\gothic{Q}\geq 0, \\
 0 &\textrm{ if }\gothic{Q}<0.
 \end{array}
 \right. , \\
  \textrm{with anisotropy measure: }
 \gothic{R}= \frac{\mathcal{B}_{1}^{\, 2} \tilde{U}}{\mathcal{B}_{2}^{\, 2} \tilde{U} + \mathcal{B}_{3}^{\, 2} \tilde{U}}, \\
 \textrm{structureness: }
 \gothic{S}=(\mathcal{B}_{1}^{\, 2} \tilde{U})^2 + (\mathcal{B}_{2}^{\, 2} \tilde{U}+\mathcal{B}_{3}^{\, 2} \tilde{U})^2\!\!,   \\
 \textrm{convexity: }
 \gothic{Q}= \mathcal{B}_{2}^{\, 2}\mathcal{U} + \mathcal{B}_{3}^{\, 2}\mathcal{U},
 \end{array}
 \end{equation}
 with $\sigma_{1}=\frac{1}{2}$ and $\sigma_{2}=0.2 \|\mathcal{B}_{2}^{\,2}\tilde{U}+\mathcal{B}_{3}^{\, 2}\tilde{U}\|_{\infty}$. Here the decomposition of the vesselness in structureness, anisotropy and convexity follows the same general principles of the vesselness. As in vessels are line like structures we use the exponential curve fits of 2nd order obtained via the symmetric product of the Hessian (i.e. solving the optimization problem in Thm.~\ref{th:0c}).

Similarly to the vesselness filter \cite{Frangi}, we need a mechanism to robustly deal with vessels of different width. This is why for this application we extend the (all-scale) orientation scores to multiple-scale invertible orientation scores. Such multiple-scale orientation scores \cite{Sharma} coincide with wavelet transforms on the similitude group $SIM(2)=\R^{2} \rtimes SO(2) \times \R^{+}$, where one uses a B-spline \cite{Unser,Felsberg} basis decomposition along the log-radial axis in the Fourier domain. In our experiments we used $N=4,12$ or $20$ orientation layers and a decomposition centered around $M=4$ discrete scales $a_l$ given by
\begin{equation} \label{scales}
a_{l}=a_{min}e^{l\, (M-1)^{-1} \log (a_{max}/a_{min})},
\end{equation}
$l=0,\ldots,M-1$ where $a_{max}$ is inverse proportional to the Nyquist-frequency $\rho_n$ and $a_{min}$ close to the inner scale \cite{Florackbook} induced by sampling (see \cite{Sharma} for details). Then, the multiple-scale orientation score is given by the following wavelet transform $\mathcal{W}_{\psi}f:SIM(2) \to \mathbb{C}$:
\begin{equation} \label{SIM2WT}
\mathcal{W}_{\psi}f(\ul{x},\theta,a)= \int \limits_{\R^{d}} \overline{\psi(a^{-1}\ul{R}_{\theta}^{-1}(\ul{y}-\ul{x}))} f(\ul{y})\, {\rm d}\ul{y},
\end{equation}
and we again set $U:=\textrm{Re}\{\mathcal{W}_{\psi}f\}$. Finally we define the total integrated multiple scale $SIM(2)$-vesselness by:
\begin{equation} \label{MSvess}
\begin{array}{l}
(\Phi^{SIM(2)}(U))(\ul{x}):= \\
\mu_{\infty}^{-1} \sum \limits_{i=0}^{M-1} \mu_{i,\infty}^{-1}
\sum \limits_{j=1}^N (\Phi(U(\cdot,\cdot,\cdot,a_i)))(\ul{x},\theta_j),
\end{array}
\end{equation}
where $SE(2)$-vesselness operator $\Phi$ is given by Eq.~\!(\ref{SE2vesselness}), and where $\mu_{\infty}$ and $\mu_{i,\infty}$ denote maxima w.r.t. sup-norm $\|\cdot\|_{\infty}$ taken over the subsequent terms.

Note that another option for constructing a $SIM(2)$-vesselness is to use the non-adaptive left-invariant frame $\{\mathcal{A}_{1},\mathcal{A}_{2},\mathcal{A}_{3}\}$ instead of the gauge frame. This non-adaptive $SE(2)$-vesselness operator is obtained by simply replacing the $\mathcal{B}_{i}$ operators by the $\mathcal{A}_{i}$ operators in Eq. (\ref{SE2vesselness}) accordingly.

The aim of the experiments presented in this section is to show the following advantages:
\begin{description}
\item[\bfseries Advantage 1:] The improvement of considering the multiple-scale vesselness filter via gauge frames in $SE(2)$, compared to multiple-scale vesselness \cite{Frangi} acting directly on images.
\item[\bfseries Advantage 2:] Further improvement when using the gauge frames instead of using the left-invariant vector fields in $SE(2)$-vesselness (\ref{SE2vesselness}).
\end{description}



\begin{figure}
\includegraphics[width=\hsize]{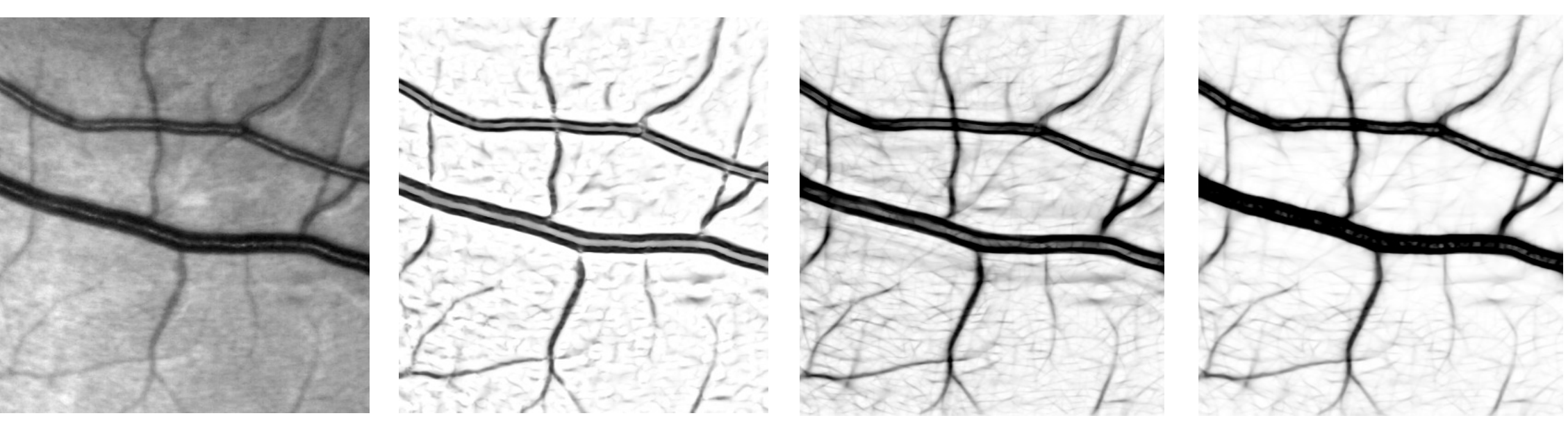}
\caption{From left to right:
Retinal image $f$ (from HRF-database), multi--scale vesselness filtering results for the multi-scale Frangi vesselness filter on $\R^{2}$, our
$SIM(2)$-vesselness via invertible multi--scale orientation score based on left-invariant frame
$\{\mathcal{A}_{1},\mathcal{A}_{2},\mathcal{A}_{3}\}$, and based on adaptive frame
$\{\mathcal{B}_{1},\mathcal{B}_{2},\mathcal{B}_{3}\}$. \label{fig:vess}}
\end{figure}

In the following experiment, we test these 3 techniques (Frangi vesselness \cite{Frangi}, $SIM(2)$-vesselness via the non-adaptive left invariant frame, and the newly proposed $SIM(2)$-vesselness via gauge frames) on the publically available\footnote{cf.\!~\mbox{{\small \emph{http://www5.cs.fau.de/research/data/fundus-images/}}}} High Resolution Fundus (HRF)-dataset \cite{Kohler2013},
containing manually segmented vascular trees by medical experts. The HRF-dataset consists of wide--field fundus photographs for a healthy, diabetic retinopathy and a glaucoma group (15 images each). A comparison of the 3 vesselness filters on a small patch is depicted Fig.\!~\ref{fig:vess}. Here, we see that our method performs better both at crossing and non-crossing structures. 


To perform a quantitative comparisson, we devised a simple segmentation algorithm to turn a vesselness filtered image $\mathcal{V}(f)$ into a segmentation. First an adaptive thresholding is applied, yielding a binary image
\begin{equation} \label{eq:IB}
f_B = \Theta\big([\mathcal{V}(f)-G_\gamma*\mathcal{V}(f)]-h\big),
\end{equation}
where $\Theta$ is the unit step function, $G_\gamma$ is a Gaussian of scale $\gamma=\frac{1}{2}\sigma^2 \gg1$ and $h$ is a threshold parameter. In a second step, the connected morphological components in $f_B$ are subject to size and elongation constraints. Components counting less than $\tau$ pixels or showing elongations below a threshold $\nu$ are removed. Parameters $\gamma, \tau$ and $\nu$ are fixed at 100 px, 500 px and 0.85 respectively. The vesselness map $\mathcal{V}(f):\R^{2} \to \R$ is one of the 3 methods considered.

\begin{figure}
\centerline{\includegraphics[width=1.\hsize]{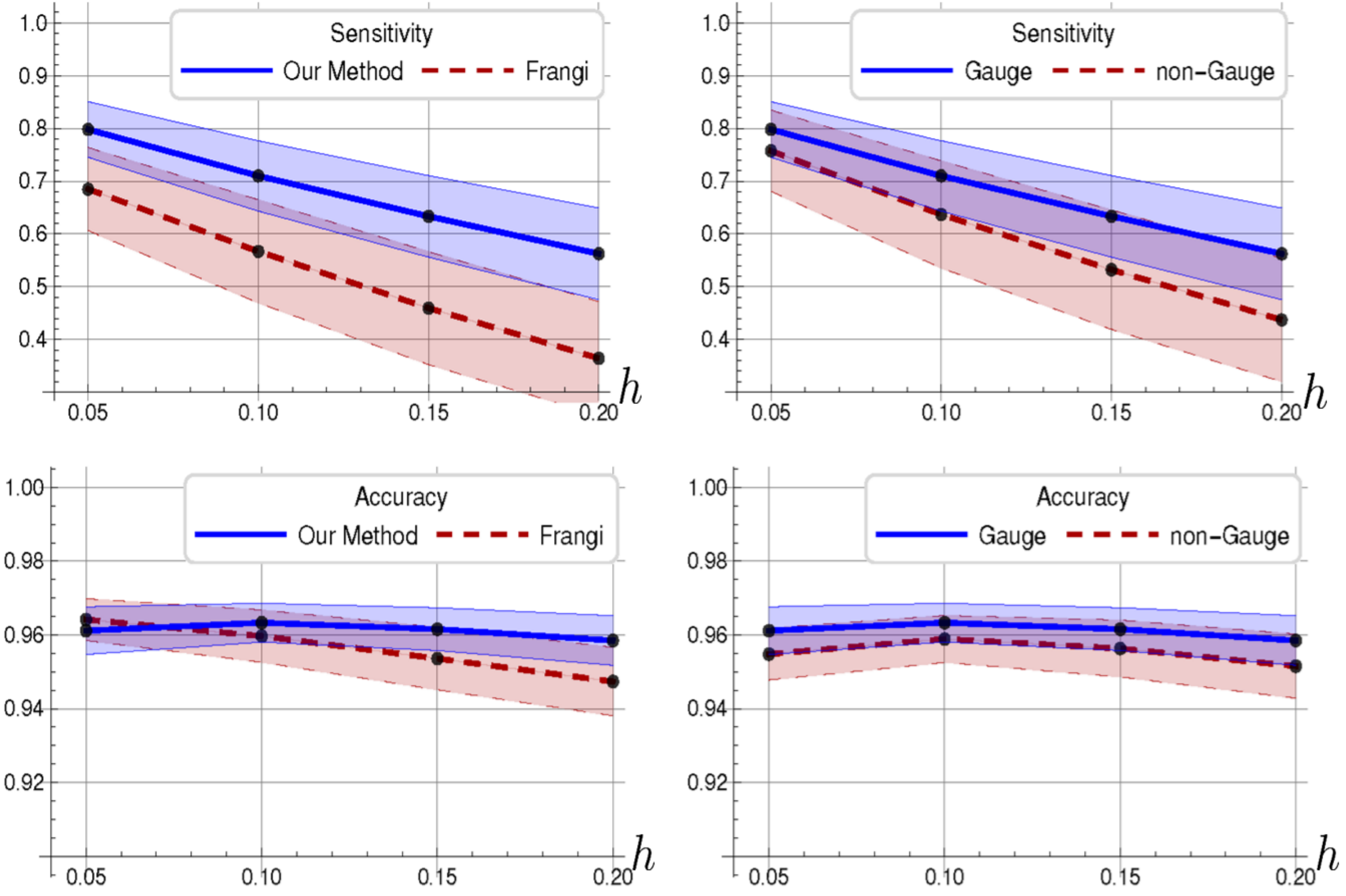}}
\caption{Left: Comparison of multiple scale Frangi vesselness and $SIM(2)$-vesselness via gauge frames. Average accuracy and sensitivity on the HRF dataset over threshold values $h$. Shaded regions correspond to $\pm1 \,\sigma$. Right: comparing of $SIM(2)$-vesselness with and without including the gauge frame (i.e. using $\{\mathcal{A}_{1},\mathcal{A}_{2},\mathcal{A}_{3}\}$ in Eq.\!~(\ref{SE2vesselness})). \label{fig:julius}}
\end{figure}%
%

The segmentation algorithm described above is evaluated on the HRF dataset. Average sensitivity and accuracy over the whole dataset are shown in Fig.~\ref{fig:julius} as a function of the threshold value $h$. It can be observed that our method performs considerably better than the one based on the multi--scale Frangi filter. The segmentation results obtained with $SIM(2)$-vesselness (\ref{MSvess}) based on gauge frames are more stable w.r.t variations in the threshold $h$ and the performance on the small vasculature has improved as measured via the sensitivity. Average sensitivity and accuracy at a threshold of $h=0.05$ compare well with other segmentation methods evaluated on the HRF dataset for the healthy cases (see \cite[Tab.~5]{Budai2013} and \cite{Hannink}). On the diabetic retinopathy and glaucoma group, our method even outperforms existing segmentation methods. 

\begin{figure*}
\centerline{\includegraphics[width=0.7\hsize]{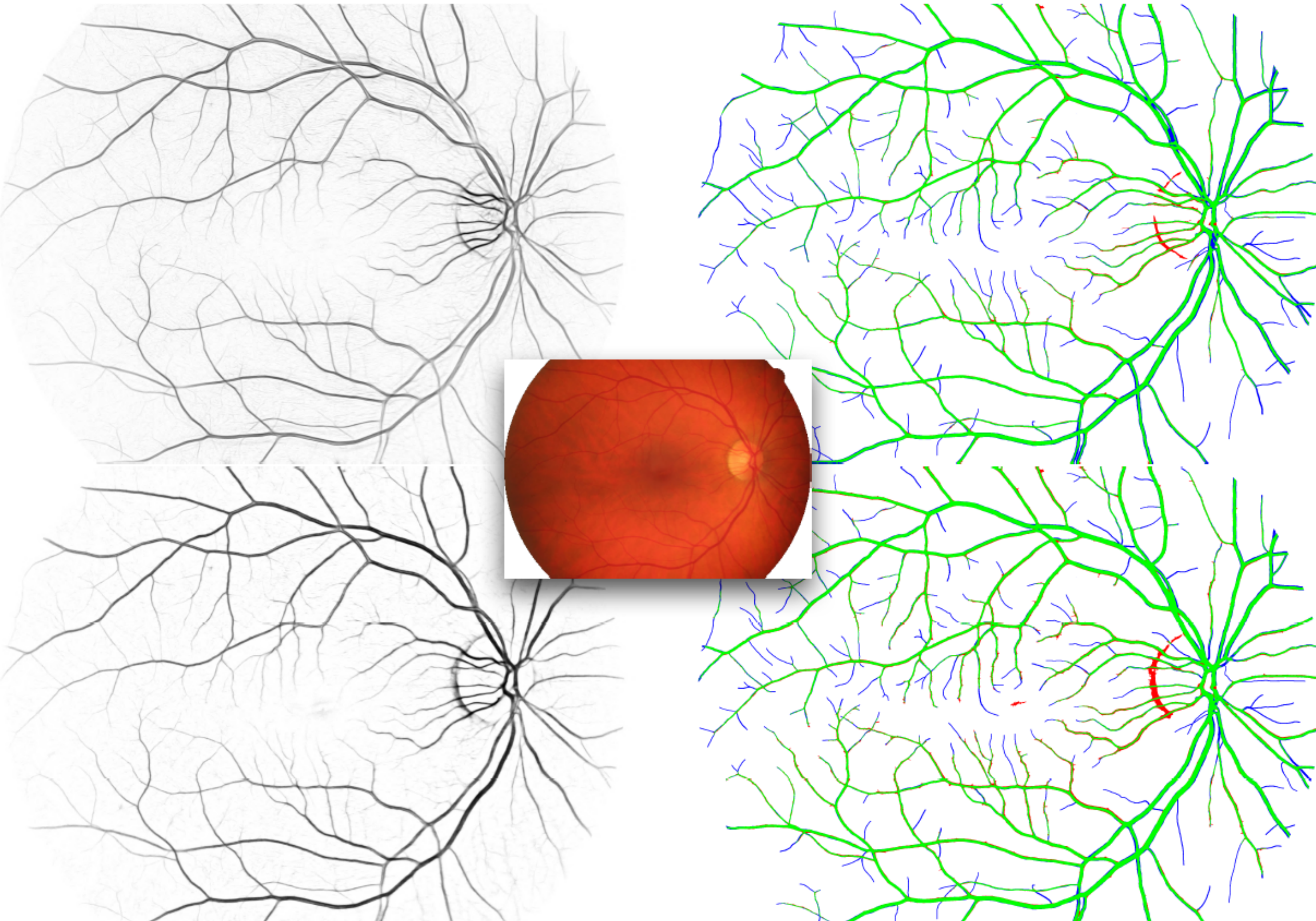}}
\caption{ Center: original image from HRF-dataset (healthy subject nr.~5).
Rows: the soft-segmentation (left) and the corresponding performance maps (right) based on the hard segmentation (\ref{eq:IB}). In green: true positives, in blue true negatives, in red false positives, compared to manual segmentation by expert. 1st row: $SIM(2)$-vesselness (\ref{MSvess}) based on non-adaptive frame $\{\mathcal{A}_{1},\mathcal{A}_{2},\mathcal{A}_{3}\}$. 2nd row: $SIM(2)$-vesselness (\ref{MSvess}) based on the gauge frame. 
\label{fig:claim2-0b}}
\end{figure*}

Finally, regarding the second advantage we refer to  Fig.~\!\ref{fig:claim2-0b}, where the $SIM(2)$-vesselness-filtering  via the locally adaptive frame produces a visually much more appealing soft-segmentation of the blood vessels than $SIM(2)$-vesselness filtering via the non-adaptive frame. It therefore also produces a more accurate segmentation as can be deducted from the comparison in Fig.~\!\ref{fig:julius}. For comparison, the multi{scale Frangi vesselness filter is also computed via summation over single scale results and max-normalized. Generally, we conclude from the experiments that the locally adaptive frame approach better reduces background noise, showing much less false positives in the final segmentation results. This can be seen from the typical segmentation results on relatively challenging patches in Fig.~\!\ref{fig:examples}.

\begin{figure}
\includegraphics[width=\hsize]{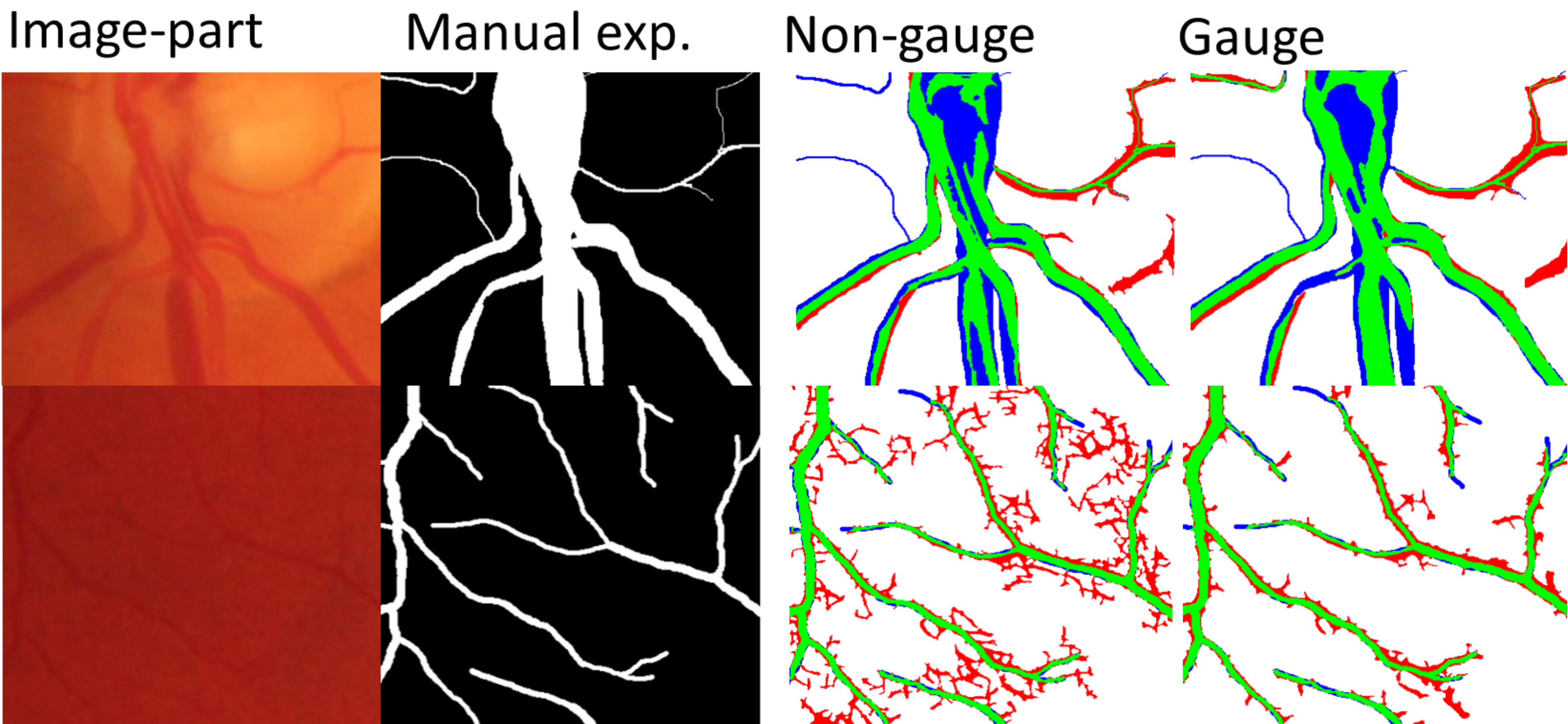}
\caption{Two challenging patches (one close to the optic disk and one far away from the optic disk processed with the same parameters). The gauge frame approach typically reduces false positives (red) on the small vessels, and increases false negatives (blue) at the larger vessels. The top patch shows a missing hole at the top in the otherwise reasonable segmentation by the expert. \label{fig:examples}}
\end{figure}

\subsection{Experiments in SE(3)}

We now show first results of the extension of coherence enhancing diffusion via invertible orientation scores (CEDOS \cite{FrankenIJCV}) of 2D images to the 3D setting. Again, data is processed according to Fig.\!~\ref{Fig:Intro}. First, we construct an orientation score according to \eqref{OS}, using the 3D cake wavelets (Fig.~\!\ref{fig:cakewavelets}). For determining the gauge frame 
we use the first order structure tensor method in combination with Eq.~\!~(\ref{gaugedD}) in Appendix~\ref{app:a}. In CEDOS we have $\Phi =\Phi_t$, as defined in \eqref{pregaugeflow} and \eqref{gaugeflow}, which is a diffusion along the gauge frame.

The diffusion in CEDOS can enhance elongated structures in 3D data while preserving the crossings as can be seen in the two examples in Fig.~\!\ref{fig:CED3}. In these experiments as well as in the example used in Fig.~\!\ref{fig:examples3D}, \ref{fig:SecondOrderFitSE3} and \ref{fig:OS3D}, we used the following 3D cake-wavelet parameters for constructing the 3D-invertible orientation scores: $N_0=42,s_\phi=0.7,k=2,N=20,\gamma=0.85,L=16$ evaluated on a grid of 21x21x21 pixels, for details see \cite{Janssen}. The settings for tangent vector estimation using the structure tensor are $s_p=\frac{1}{2}(1.5)^2,s_0=0,\rho_o=\frac{1}{2}(0.8)^2$ and $\mu=0.5$. We used $\rho_p=\frac{1}{2}(2)^2$ for the first dataset (Fig.\!~\ref{fig:CED3} top), and $\rho_p=\frac{1}{2}(3.5)^2$ for the second dataset (Fig.\!~\ref{fig:CED3} bottom). For the diffusion we used $t=2.5,D_{11}=D_{22}=0.01,D_{33}=1,D_{44}=D_{55}=D_{66}=0.04$, where the diffusion matrix is given w.r.t. gauge frame $\{\mathcal{B}_{1},\mathcal{B}_{2},\mathcal{B}_{3},\mathcal{B}_{4},\mathcal{B}_{5},\mathcal{B}_{6}\}$, and normalized frame $\{\mu^{-1} \mathcal{A}_{1},\mu^{-1} \mathcal{A}_{2},\mu^{-1} \mathcal{A}_{3},\mathcal{A}_{4},\mathcal{A}_{5},\mathcal{A}_{6}\}$.

\begin{figure}%
\centerline{\subfigure[{\small 3D Data}]{\includegraphics[width=0.33 \hsize]{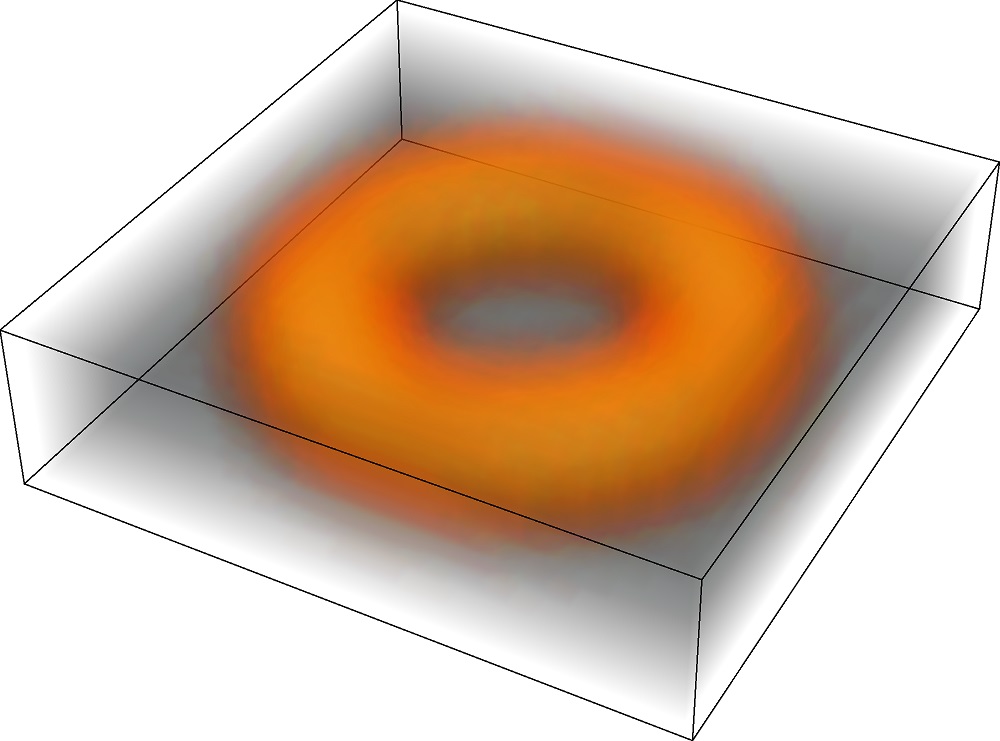}}
\subfigure[{\small Slice of Data}]{\includegraphics[width=0.33\hsize]{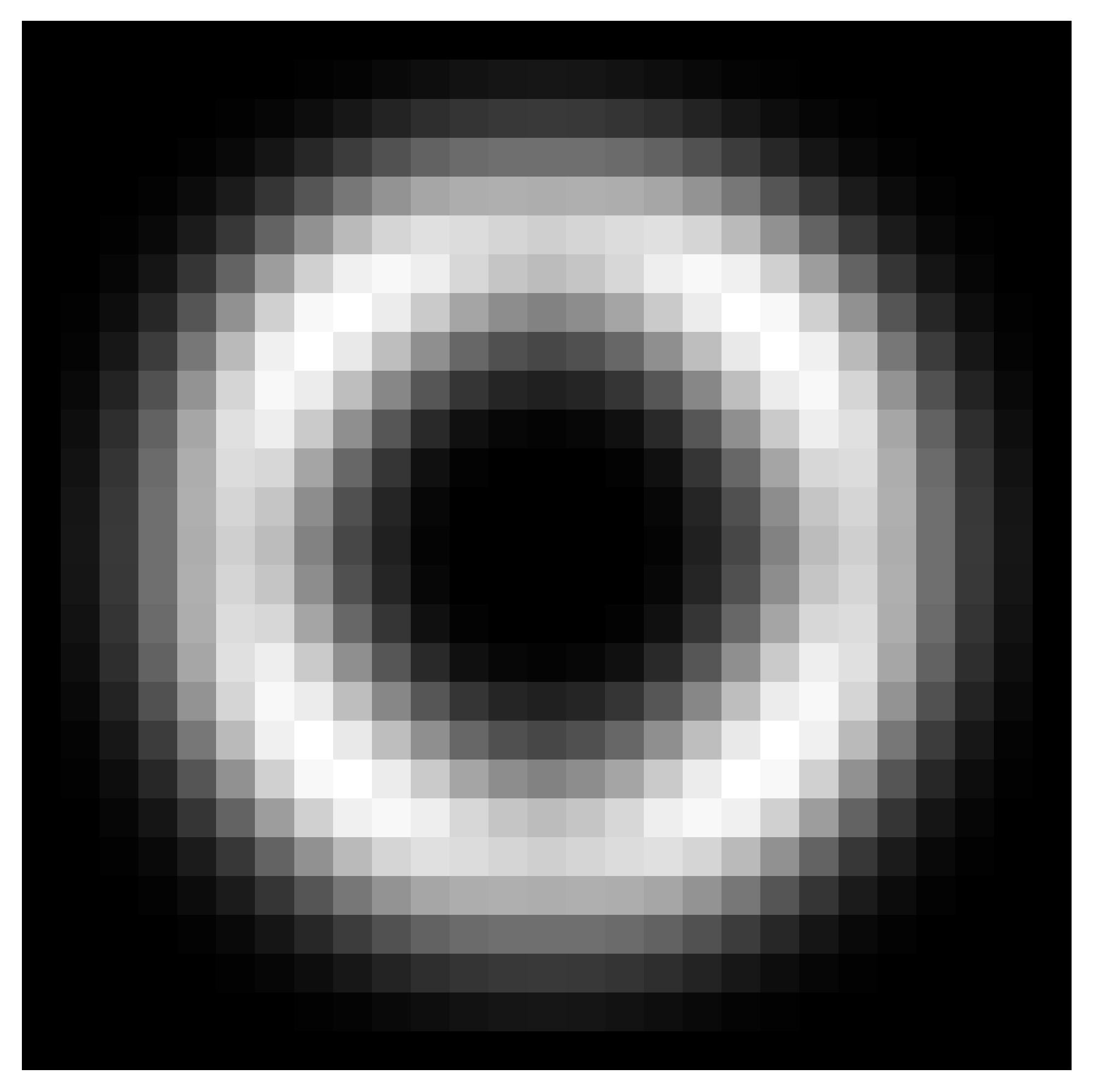}}
\subfigure[{\small Curve fits}]{\includegraphics[width=0.33\hsize]{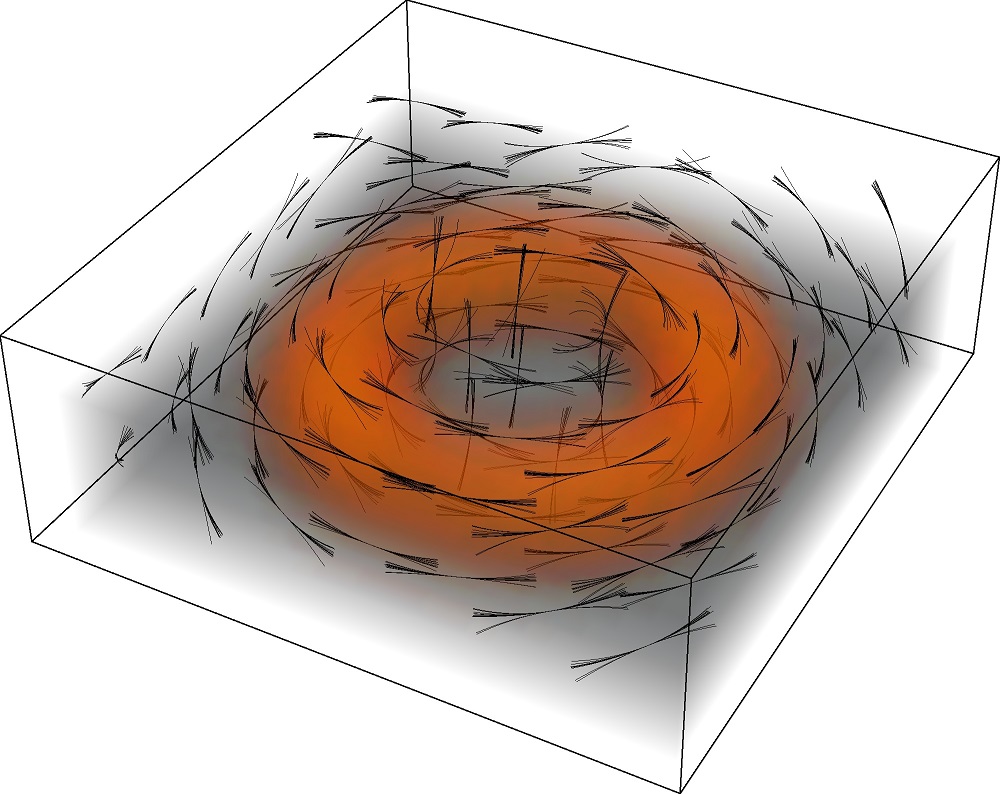}}}
\centerline{\subfigure[{\small Slice+Noise}]{\includegraphics[width=0.33\hsize]{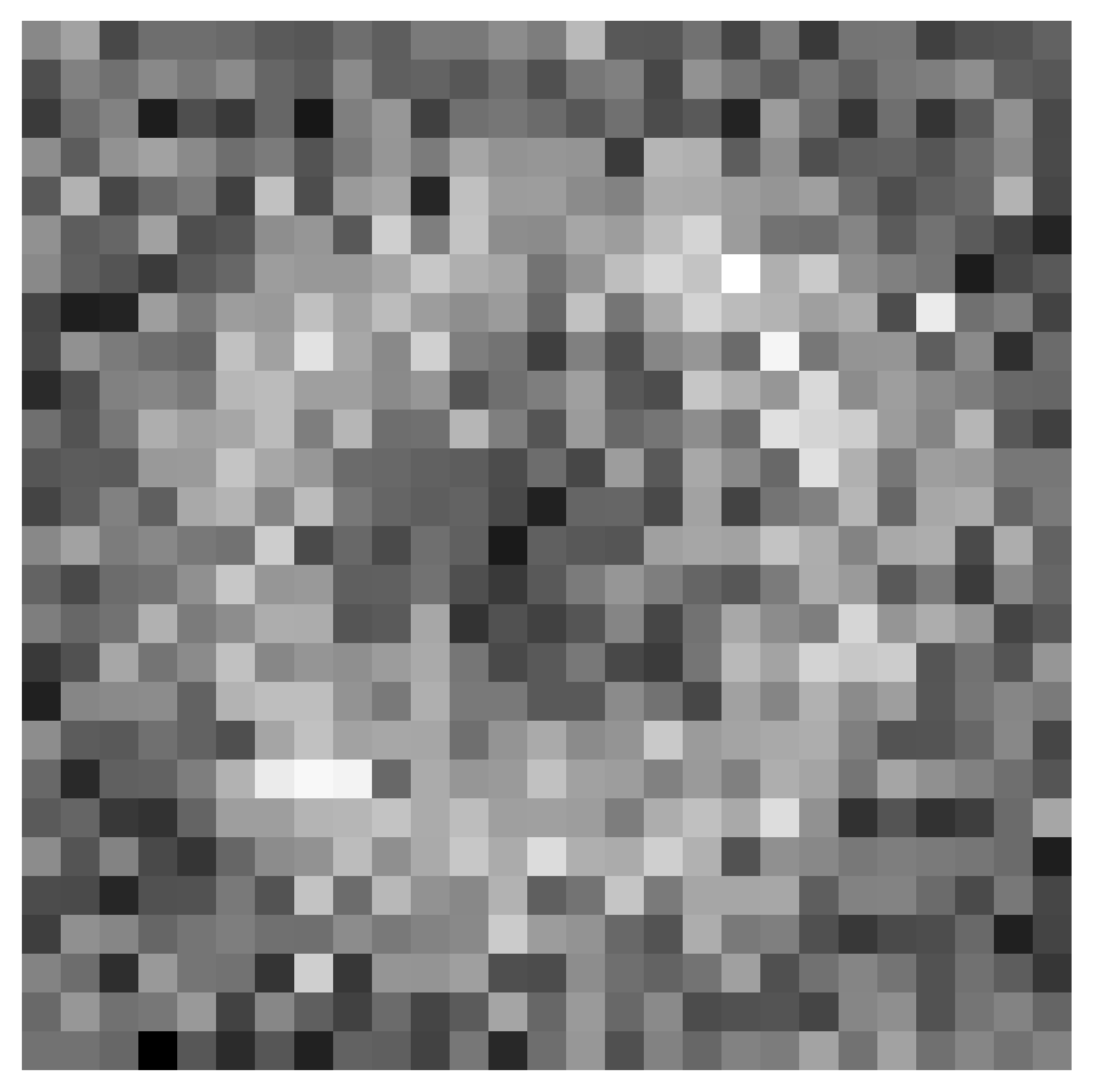}}
\subfigure[{\small Gauge}]{\includegraphics[width=0.33\hsize]{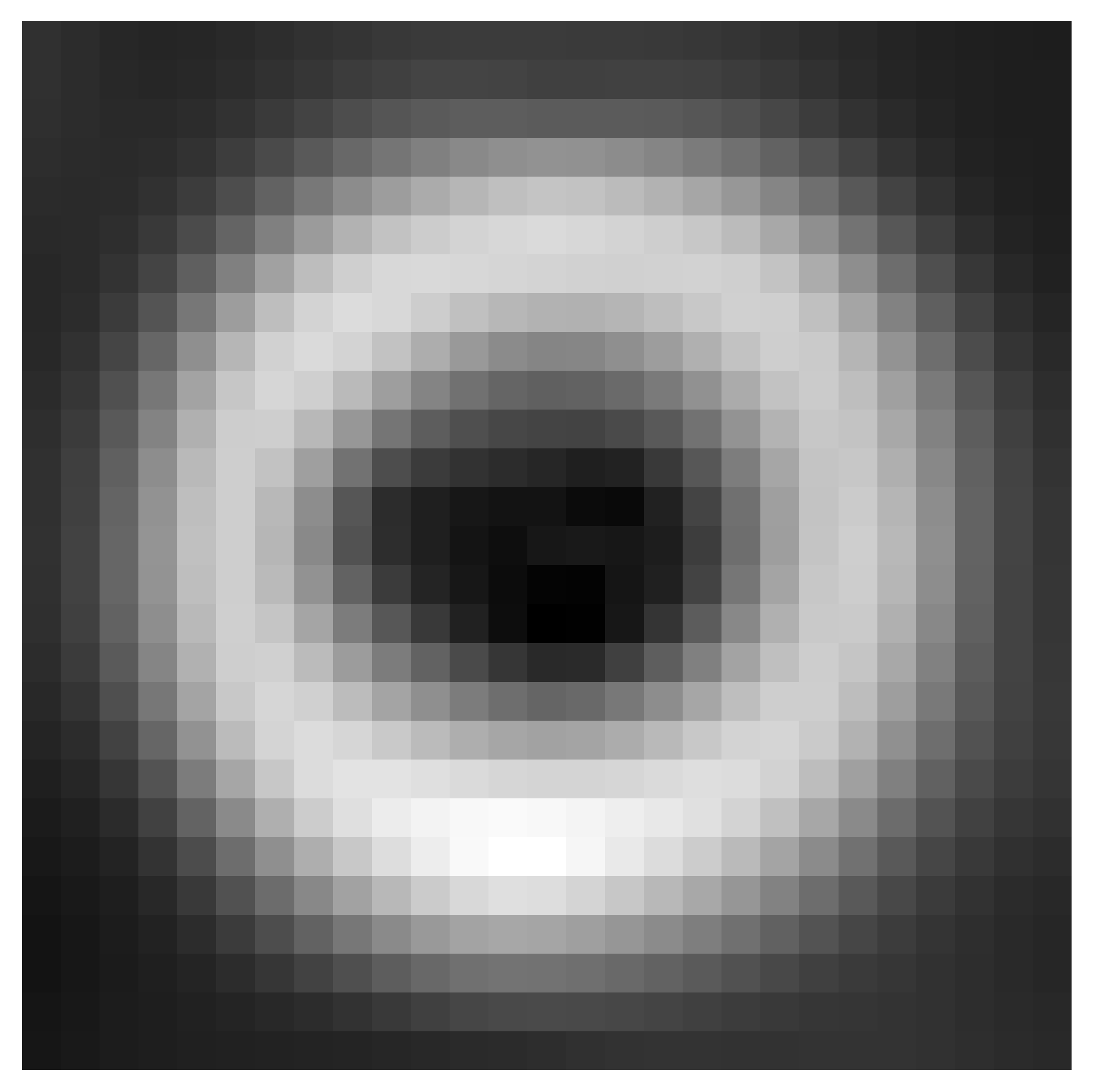}}
\subfigure[{\small No Gauge}]{\includegraphics[width=0.33\hsize]{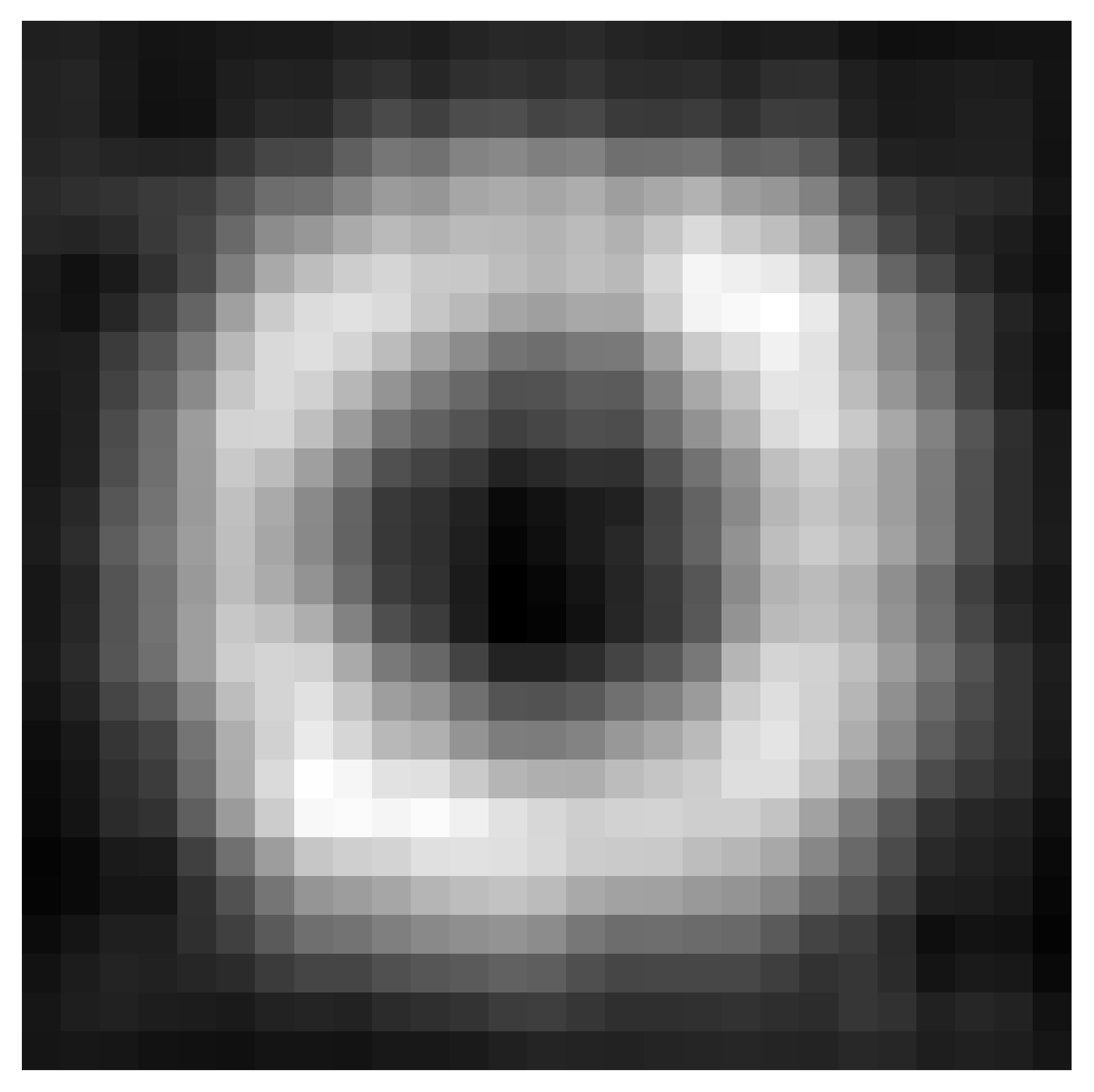}}}
\centerline{\subfigure[{\small 3D Data}]{\includegraphics[width=0.34 \hsize]{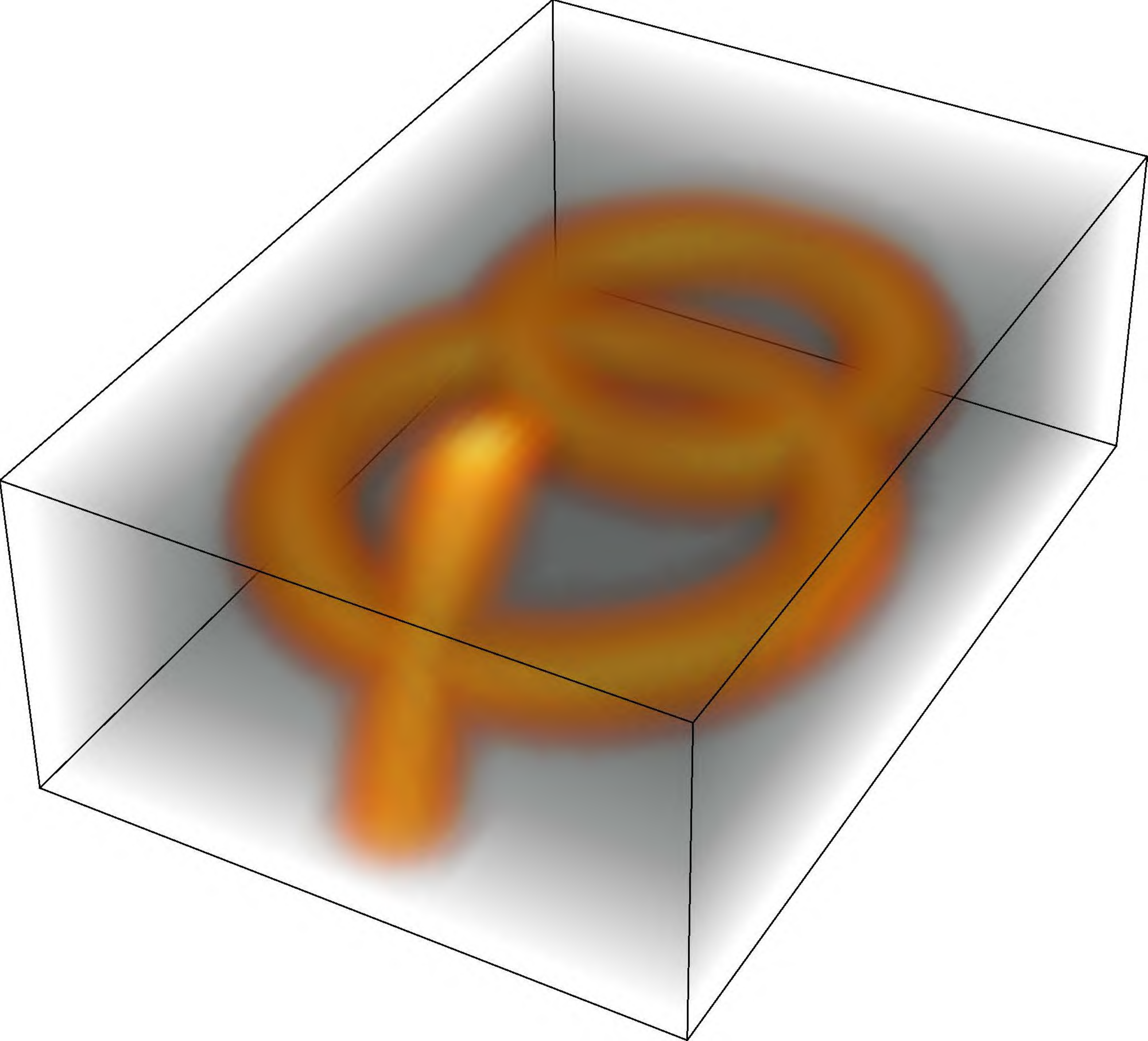}}
\subfigure[{\small Data Slice}]{\includegraphics[width=0.24\hsize]{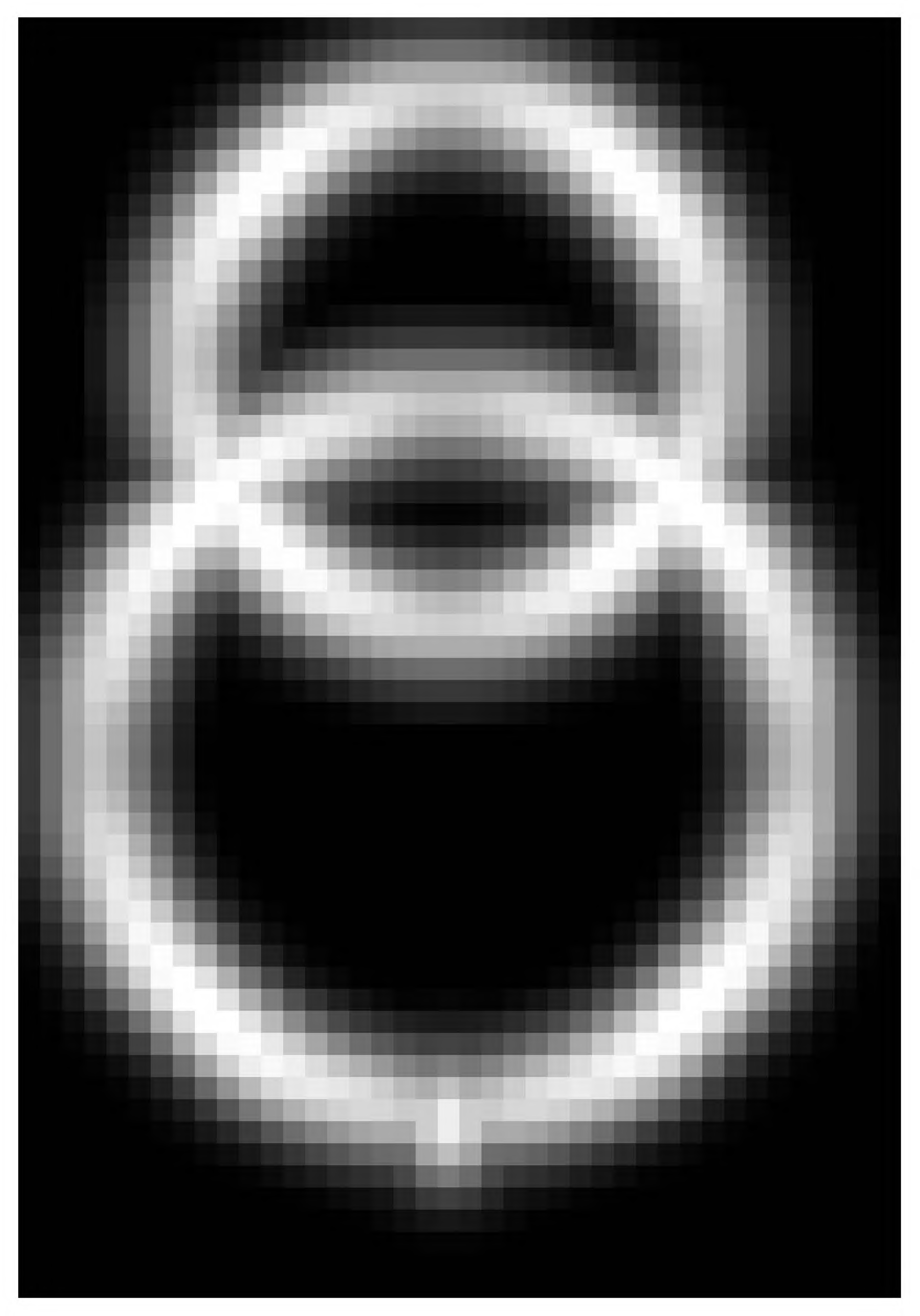}}
\subfigure[{\small Curve fits}]{\includegraphics[width=0.37\hsize]{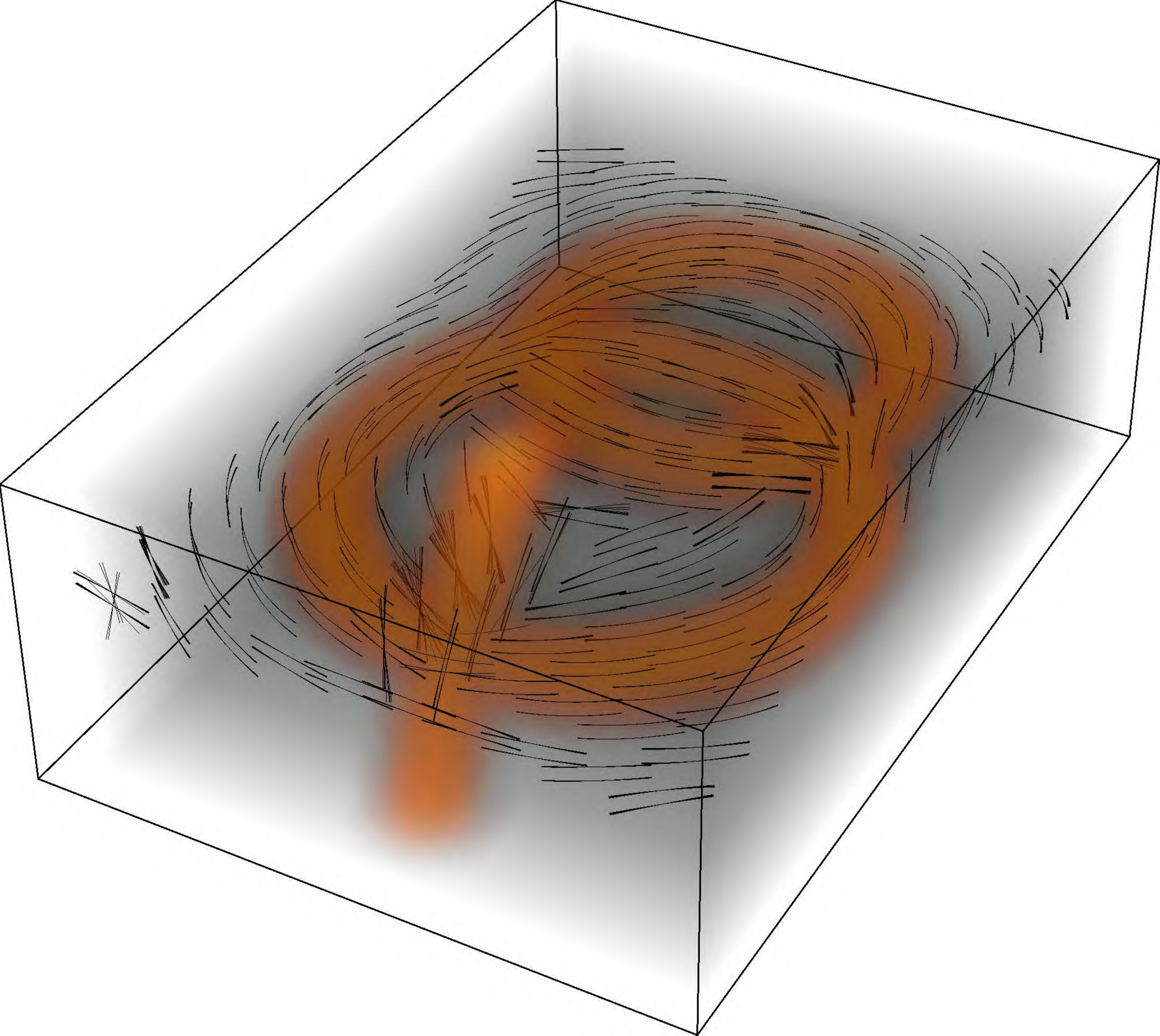}}}
\centerline{ \subfigure[{\small Slice+Noise}]{\includegraphics[width=0.33\hsize]{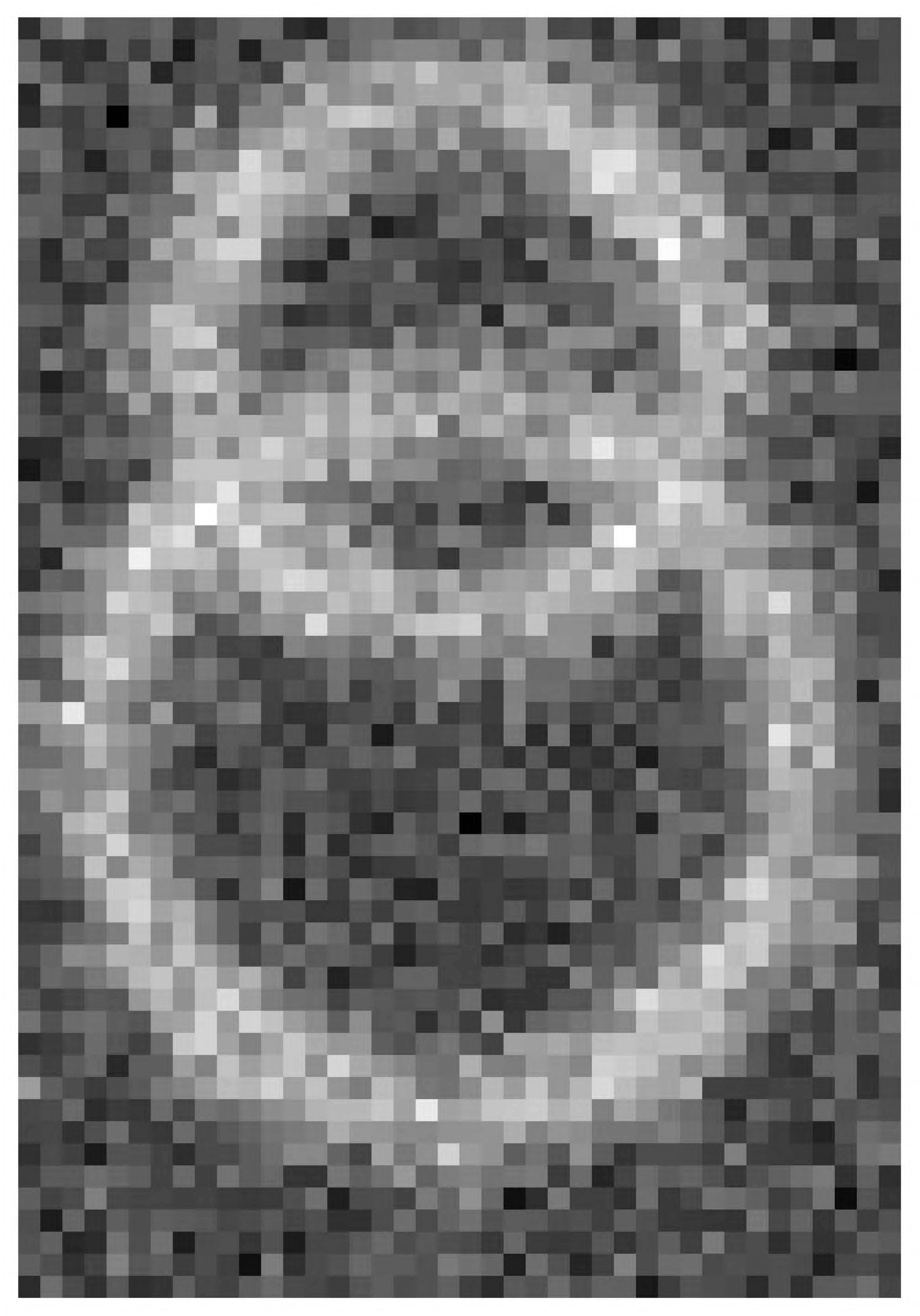}}
\subfigure[{\small Gauge}]{\includegraphics[width=0.33\hsize]{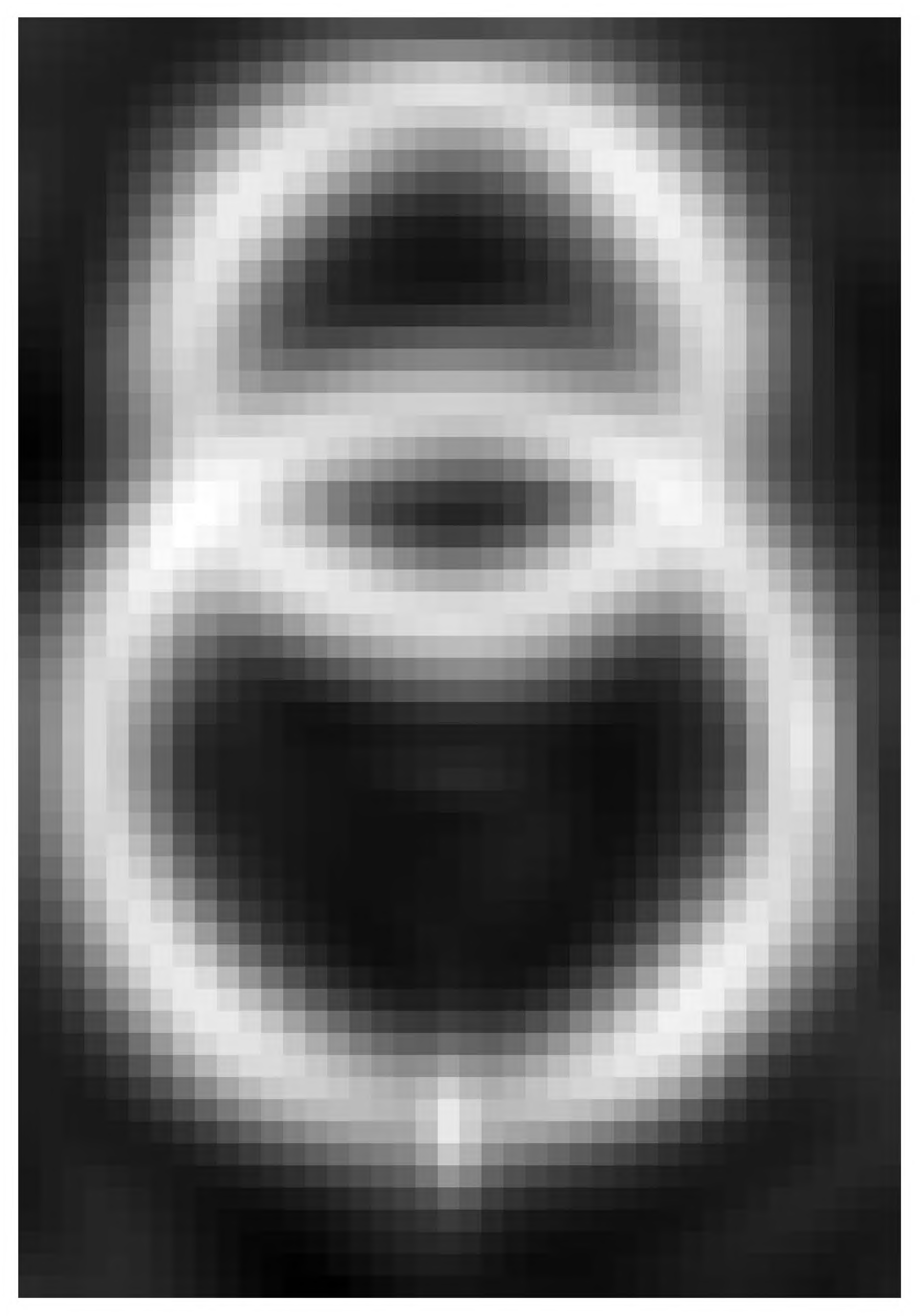}}
\subfigure[{\small No Gauge}]{\includegraphics[width=0.33\hsize]{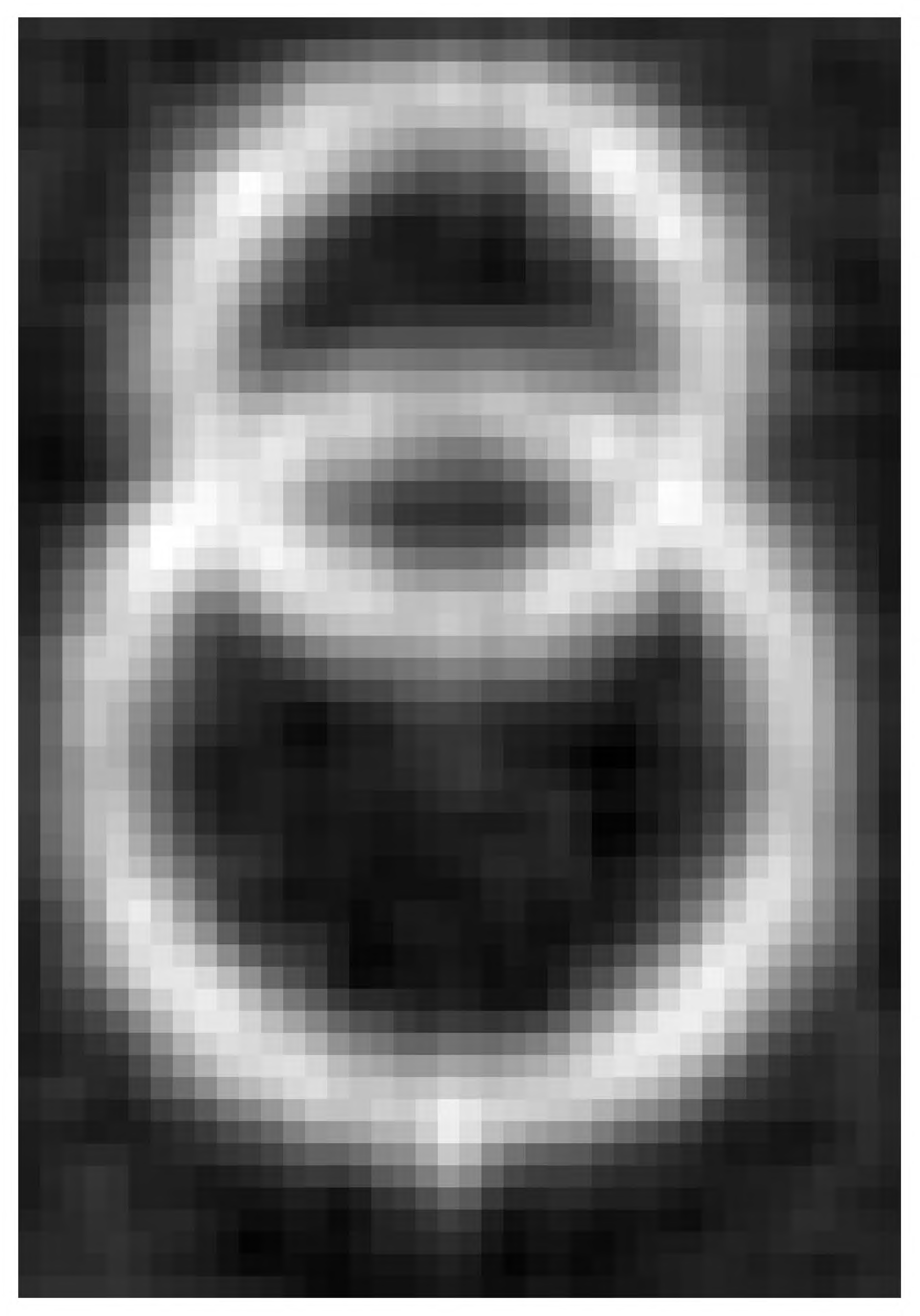}}}
\caption{Results of CEDOS with and without the use of gauge frames, on 3D artificial datasets containing highly curved structures. Gauge frames are obtained, see Appendix~\ref{app:a}, via 1st order exponential curve fits using the two-fold algorithm of Subsection~\ref{ch:3-4}.} \label{fig:CED3}
\end{figure}

The advantages of including the gauge frames w.r.t. the non adaptive frame can be better appreciated
in Fig.~\!\ref{fig:CED2OS}. Here, we borrow from the neuroimaging community the {\em{glyph visualization}}, a standard technique for displaying distributions $U:\mathbb{R}^3\times S^2\to \mathbb{R}^+$. In such visualizations every voxel contains a spherical surface plot (a glyph) in which the radial component is proportional to the output-value of the distribution at that orientation, and the colors indicate the orientations.  One can observe that diffusion along the gauge frames include better adaptation for curvature. This is mainly due to the angular part in the $\mathcal{B}_{3}$-direction, cf.~\!Fig.~\ref{fig:gauge3D}, which includes curvature, in contrast to $\mathcal{A}_{3}$-direction. The angular part in $\mathcal{B}_{3}$ causes some additional angular blurring leading to more isotropic glyphs.


\begin{figure}
\centerline{
\subfigure[{\small Data}]{\includegraphics[width=0.21\textwidth]{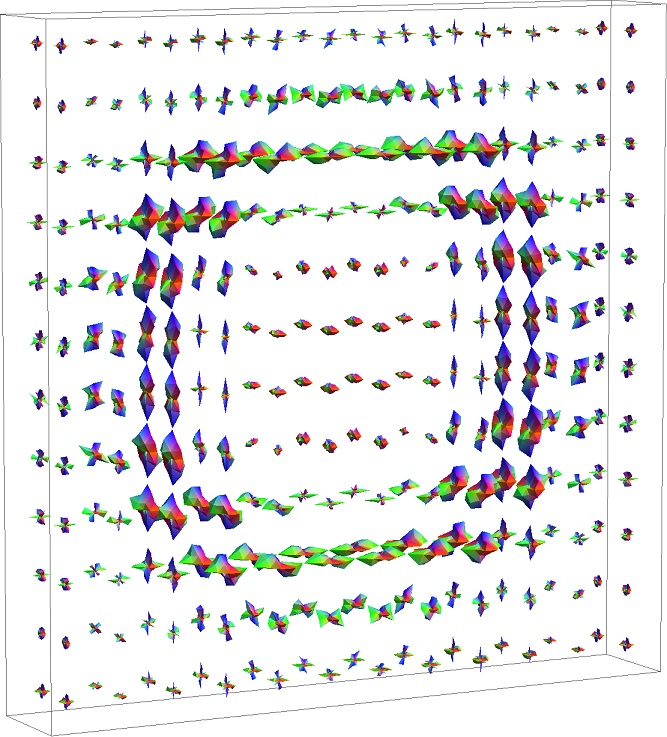}}
\subfigure[{\small Data \& Rician Noise}]{\includegraphics[width=0.21\textwidth]{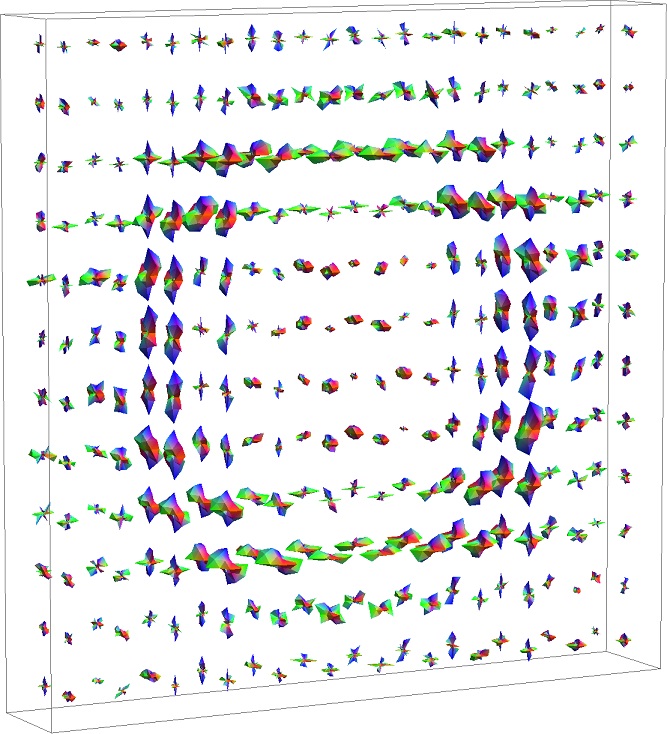}}}
\centerline{\subfigure[{\small Enhanced using frame $\{\mathcal{B}_{1},\ldots,\mathcal{B}_{6}\}$}]{\includegraphics[width=0.21\textwidth]{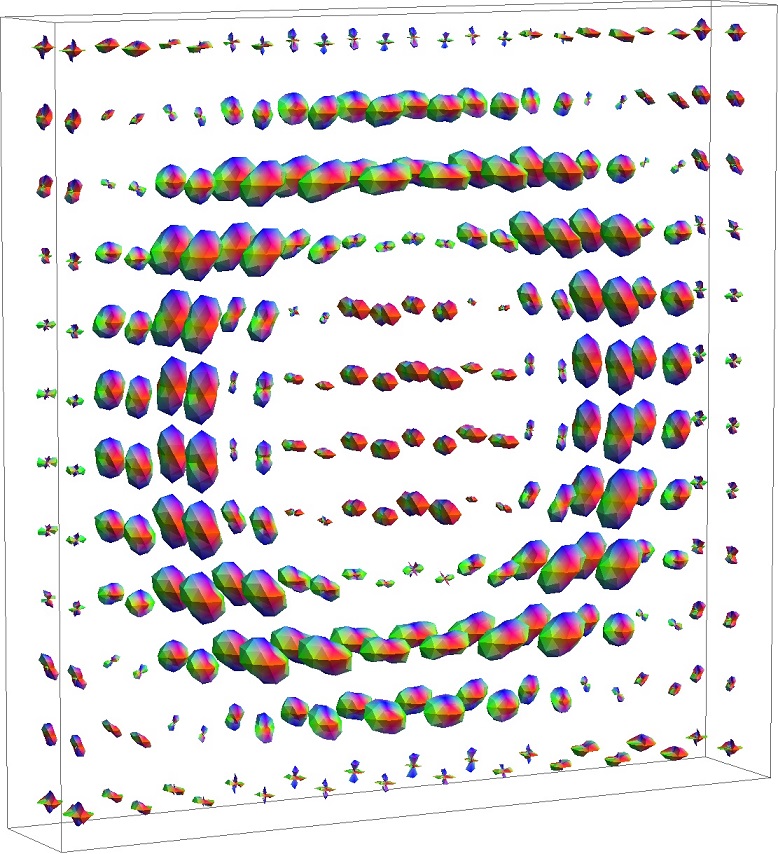}}
\subfigure[{\small Enhanced using frame $\{\mathcal{A}_{1},\ldots,\mathcal{A}_{6}\}$}]{\includegraphics[width=0.21\textwidth]{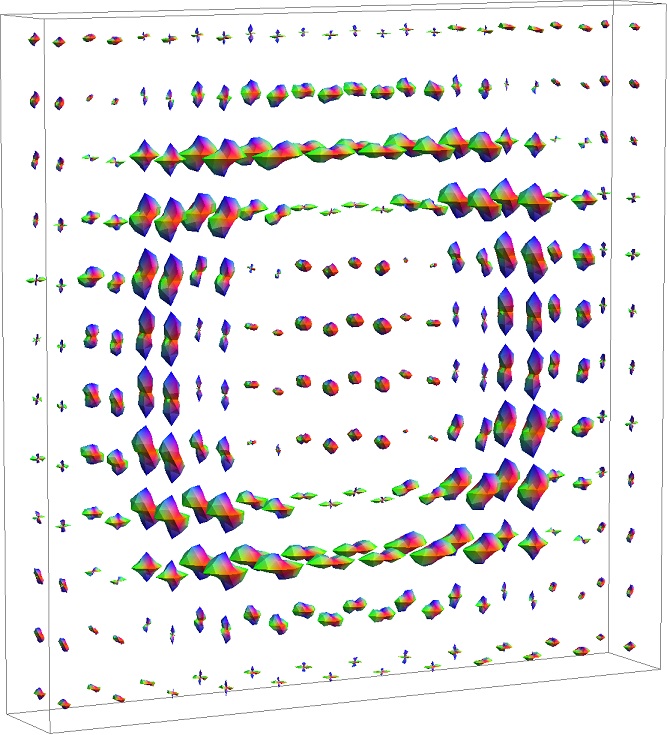}}}
\caption{Glyph visualization (see text) of the absolute value of the diffused orientation scores with and without the use of gauge frames in the the artificial dataset depicted in Fig.~\!\ref{fig:CED3} top.\label{fig:CED2OS}}%
\end{figure}

\section{Conclusion \label{ch:conlusion}}
Locally adaptive frames (`gauge frames') on images based on the structure tensor or Hessian of the images are ill-posed at the vicinity of complex structures. Therefore we create locally adaptive frames on distributions on $SE(d)$, $d=2,3$ that extend the image domain (with positions and orientations). This gives rise to a whole family of local frames per position, enabling us to deal with crossings and bifurcations. In order to generalize gauge frames in the image domain to gauge frames in $SE(d)$, 
we have shown that exponential curve fits gives rise to suitable gauge frames. 
We distinguished between exponential curve fits of the 1st order and of the 2nd order:
\begin{enumerate}
\item Along the 1st order exponential curve fits, the 1st order variation of the data (on $SE(d)$) along the exponential curve is locally minimal. The Euler-Lagrange equations are solved by finding the eigenvector of the structure tensor of the data, with smallest eigenvalue.
\item Along the 2nd order exponential curve fits, a 2nd order variation of the data (on $SE(d)$) along the exponential curve is locally minimal. The Euler-Lagrange equations are solved by finding the eigenvector of the Hessian of the data, with smallest eigenvalue.
\end{enumerate}
In $SE(2)$, the 1st order approach is new while the 2nd order approach formalizes previous results. In $SE(3)$, these two approaches are presented for the first time. Here, it is necessary to include a restriction to torsion-free exponential curve fits in order to be both compatible with the null-space of the structure/Hessian tensors and the quotient structure of $\R^{3} \rtimes S^{2}$. We have presented an effective two-fold algorithm to compute such torsion-free exponential curve fits. Experiments on artificial datasets show that even if the elongated structures have torsion, the gauge frame 
is well-adapted to the local structure of the data.


Finally, we considered the application of a differential invariant for enhancing retinal images. Experiments show clear advantages over the classical vesselness filter \cite{Frangi}. Furthermore, we also show clear advantages of including the gauge frame over the standard left-invariant frame in $SE(2)$.
Regarding 3D image applications, we managed to construct and implement crossing-preserving coherence enhancing diffusion via invertible orientation scores (CEDOS), for the first time. However, it has only been tested on artificial datasets. Therefore, in future work we will study the use of locally adaptive frames in real 3D medical imaging applications, e.g. in 3D MR angiography \cite{MichielMaster}. Furthermore, in future work we will apply the theory of this work and focus on the explicit algorithms, where we plan to release \emph{Mathematica}-implementations of locally adaptive frames in $SE(3)$.

\section*{Acknowledgements}
The authors wish to thank J.M.~Portegies for fruitful discussions on the construction of gauge frames in $SE(3)$ and T.C.J.~ Dela Haije for help in optimizing code for the $SE(3)$-case. Finally, we would like to thank Dr.~A.J.E.M.~Janssen for careful reading and valuable suggestions on the structure of the paper. The research leading to these results has received funding from the European Research Council under the European Community's Seventh Framework Programme (FP7/2007-2013) / ERC grant \emph{Lie Analysis}, agr.~nr.~335555. \\
\includegraphics[width=0.3\hsize]{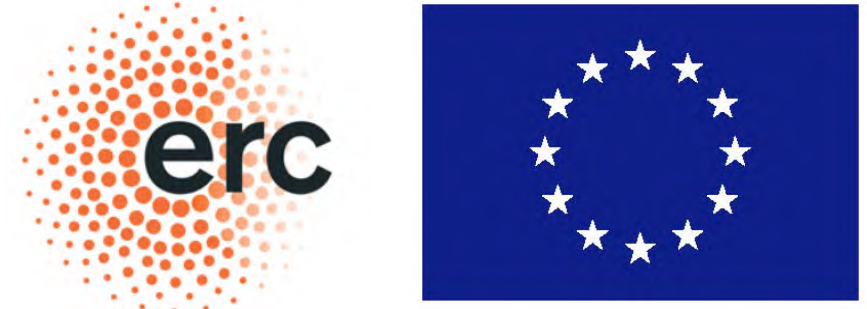}

\appendix
\section{Construction of the Locally Adaptive Frame from an Exponential Curve Fit\label{app:a}}

Let $\tilde{\gamma}_{g}^{\ul{c}}(t)= g \, e^{t\sum \limits_{i=1}^{n_d}c^{i}A_{i}}$ be an exponential curve through $g$ that fits data $\tilde{U}:SE(d) \to \mathbb{R}$ at $g \in SE(d)$ in
Lie group $SE(d)$ of dimension $n_d=d(d+1)/2$. In Section~\ref{ch:SE2} ($d=2$), and in Section~\ref{ch:SE3} ($d=3$), we provide theory and algorithms to derive such curves.
In this section we assume $\gamma_{g}^{\ul{c}}(\cdot)$ is given.

Recall from (\ref{ec}) that the (physical) velocity at time $t$ of the exponential curve $\tilde{\gamma}_{g}^{\ul{c}}$ equals
$(\tilde{\gamma}_{g}^{\ul{c}})'(t)=\sum \limits_{i=1}^{n_d} c^{i} \left.\mathcal{A}_{i}\right|_{g=\tilde{\gamma}_{g}^{\ul{c}}(t)}$.
Recall that the spatial and respectively rotational components of the
velocity are stored in the vectors
\[
\begin{array}{l}
\ul{c}^{(1)}=(c^1,\ldots, c^{d})^T \in \R^{d}, \\
\ul{c}^{(2)}=(c^{d+1}, \ldots, c^{n_d})^T \in \mathbb{R}^{r_d}
\end{array}
\]
Let us write $\ul{c}=\left(\begin{array}{c}\ul{c}^{(1)} \\\ul{c}^{(2)}\end{array}\right) \in \R^{n_d}$, with $n_d=d+r_d$.

Akin to the case $d=2$ discussed in the introduction we define the Gauge frame via $\underline{\mathcal{B}}:= (\ul{R}^{\ul{c}})^T \ul{M}_{\mu}^{-1}\underline{\mathcal{A}}$, but now with
\begin{equation} \label{gaugedD}
\begin{array}{l}
\underline{\mathcal{B}}=(\mathcal{B}_{1},\ldots,\mathcal{B}_{n_d})^T, \
\underline{\mathcal{A}}=(\mathcal{A}_{1},\ldots,\mathcal{A}_{n_d})^T, \\[6pt]
\ul{M}_{\mu}= \begin{pmatrix}
   \mu I_d  && \ul{0} \\
   \ul{0} && I_{r_d} \end{pmatrix},
   \textrm{ and }
\ul{R}^{\ul{c}}=\ul{R}_{2} \ul{R}_{1} \in SO(n_d).
\end{array}
\end{equation}
For explicit formulae of the left-invariant vector fields in the $d$-dimensional case we refer to \cite{bookchapter}.

Now $\ul{R}_{1}$ is the counter-clockwise rotation that
rotates the spatial reference axis $\left(\begin{array}{c}\ul{a} \\\ul{0} \end{array}\right)$, recall our convention (\ref{conventie}),
onto $\left(
\begin{array}{c}
\mu\|\ul{c}^{(1)}\|\ul{a} \\
\ul{c}^{(2)}
\end{array}
\right)$
strictly within the 2D-plane spanned by these two vectors. Rotation $\ul{R}_{2}$ is the counter-clockwise rotation that rotates $\left(
\begin{array}{c}
\mu\|\ul{c}^{(1)}\|\ul{a} \\ \ul{c}^{(2)}
\end{array}\right)$ onto
$\left(\begin{array}{c}\mu \ul{c}^{(1)}\\ \ul{c}^{(2)} \end{array}\right)$ strictly within the 2D-plane spanned by these two vectors. As a result one has
\begin{equation} \label{R1R2construction}
\begin{array}{l}
\left(
\begin{array}{c}
\ul{a} \\
\ul{0}
\end{array}
\right)
 \;\overset{\ul{R}_{1}}{\mapsto}\;
\left(
\begin{array}{c}
\mu \|\ul{c}^{(1)}\|\ul{a} \\
\ul{c}^{(2)}
\end{array}
\right)
\;\overset{\ul{R}_{2}}{\mapsto} \;
\left(
\begin{array}{c}
\mu\ul{c}^{(1)} \\
\ul{c}^{(2)}
\end{array}
\right)=\ul{M}_{\mu}\ul{c} \desda \\
\ul{c}=\ul{M}_{\mu}^{-1}\, \ul{R}^{\ul{c}} \left(
\begin{array}{c}
\ul{a} \\
\ul{0}
\end{array}
\right).
\end{array}
\end{equation}
In particular we have that the preferred spatial direction
$\left(\begin{array}{c}
\ul{a} \\
\ul{0}
\end{array}\right)
 \cdot \underline{\mathcal{A}}$ is mapped onto $\left(\begin{array}{c}
\ul{a} \\
\ul{0}
\end{array}\right) \cdot \underline{\mathcal{B}}= \ul{c} \cdot \underline{\mathcal{A}}$.

The next theorem shows us that our choice of assigning an entire gauge frame to a single exponential curve fit is the right one for our applications. 
\begin{theorem} \label{app:th} (construction of the gauge frame)
Let $\ul{c}(g)$ denote the local tangent components of exponential curve fit $t \mapsto \tilde{\gamma}^{\ul{c}(g)}_g(t)$ at $g=(\ul{x},\ul{R}) \in SE(d)$ in the data given by
$
\tilde{U}(\ul{x},\ul{R})=U(\ul{x},\ul{R}\ul{a})$.
Consider the mapping of the frame of left-invariant vector fields $\left.\underline{\mathcal{A}}\right|_{g}$ to the locally adaptive frame:
\begin{equation} \label{Baction}
\left.\underline{\mathcal{B}}\right|_{g}:= (\ul{R}^{\ul{c}(g)})^T \ul{M}_{\mu}^{-1}\left.\underline{\mathcal{A}}\right|_{g},
\end{equation}
with $\ul{R}^{\ul{c}}=\ul{R}_{2} \ul{R}_{1} \in SO(n_d)$,
with subsequent counter-clockwise \emph{planar} rotations $\ul{R}_{1}, \ul{R}_{2}$ given by (\ref{R1R2construction}).
Then the mapping $\underline{\mathcal{A}}_{g} \mapsto \underline{\mathcal{B}}_{g}$ has the following properties:
\begin{itemize}
\item The main spatial tangent direction $(L_{g})_* {\small \left(\begin{array}{c}\ul{a}\\ \ul{0}\end{array}\right)} \cdot \left.\underline{\mathcal{A}}\right|_{e}$
is mapped to exponential curve fit direction $\ul{c}^{T}(g) \cdot \left.\underline{\mathcal{A}}\right|_{g}$.
\item Spatial left-invariant vector fields that are $\gothic{G}_{\mu}$-orthogonal to this main spatial direction stay in the spatial part of the tangent space $T_{g}(SE(d))$ under rotation $\ul{R}^{\ul{c}}$ and they are invariant up to
normalization under the action
(\ref{Baction}) if and only if
the exponential curve fit is horizontal.
\end{itemize}
\end{theorem}
\textbf{Proof }
Regarding the first property we note that
\[(L_{g})_* \left(\begin{array}{c}\ul{a}\\ \ul{0}\end{array}\right) \cdot \left.\underline{\mathcal{A}}\right|_{e}=
\left(\begin{array}{c}\ul{a}\\ \ul{0}\end{array}\right) \cdot \left.\underline{\mathcal{A}}\right|_{g}
\]
as left-invariant vector fields are obtained by push-forward of the left multiplication.
Furthermore, by Eq.~\!(\ref{Baction}) and Eq.~\!(\ref{R1R2construction}) we have 
\[
\begin{array}{l}
\left(
\begin{array}{c}
\ul{a} \\
\ul{0}
\end{array}
\right)
\cdot \underline{\mathcal{B}}
=
\ul{M}_{\mu}^{-1} \ul{R}^{\ul{c}}
\left(
\begin{array}{c}
\ul{a} \\
\ul{0}
\end{array}
\right)
\cdot \underline{\mathcal{A}}
=
\ul{c} \cdot \underline{\mathcal{A}}.
\end{array}
\]
Regarding the second property, we note that if $\ul{b} \cdot \ul{a}=0 \Rightarrow$ 
\[
\ul{R}^{\ul{c}}
\left(
\begin{array}{c}
\ul{b} \\
\ul{0}
\end{array}
\right)=\ul{R}_{2}\ul{R}_{1}\left(
\begin{array}{c}
\ul{b} \\
\ul{0}
\end{array}
\right)
=\ul{R}_{2}
\left(
\begin{array}{c}
\ul{b} \\
\ul{0}
\end{array}
\right)
\]
and $\tilde{\gamma}_{g}^{\ul{c}}$ is horizontal iff $\frac{\ul{c}^{(1)}}{\|\ul{c}^{(1)}\|}=\ul{a}$ in which case
the planar rotation $\ul{R}_{2}$ reduces to the identity and $\ul{R}_{2}\ul{R}_{1}\left(
\begin{array}{c}\ul{b} \\
\ul{0}
\end{array}
\right)=\left(\begin{array}{c}
\ul{b} \\\ul{0}
\end{array}
\right)^T$ and only spatial normalization by $\mu^{-1}$ is applied.
$\hfill \Box$
\begin{remark}
For $d=2$ and $\ul{a}=(1,0)^T$ the above theorem can be observed in Fig.~\!\ref{fig:gauge}, where main spatial direction $\mathcal{A}_{1}=\cos \theta \,\partial_{x} +\sin \theta \,\partial_{y}$ is mapped onto $\mathcal{B}_{1}=
\ul{c} \cdot \underline{\mathcal{A}}$ and where $\mathcal{A}_{2}$ is mapped onto $\mathcal{B}_{2}=\mu^{-1}(-\sin \chi \mathcal{A}_{1} + \cos \chi \mathcal{A}_{2})$.
\end{remark}
\begin{remark}
For $d=3$ and $\ul{a}=(0,0,1)^T$ the above theorem can be observed in Fig.~\!\ref{fig:gauge3D}, where main spatial direction $\mathcal{A}_{3}=\ul{n} \cdot \nabla_{\R^{3}}$ is mapped onto $\mathcal{B}_{3}=\ul{c} \cdot \underline{\mathcal{A}}$, and where $\mathcal{A}_{1}$ and $\mathcal{A}_{2}$ are mapped to the strictly spatial generators $\mathcal{B}_{1}$ and $\mathcal{B}_{2}$.
For further details see \cite{MichielMaster}.
\end{remark}
\begin{figure}
\centerline{
\includegraphics[width=\hsize]{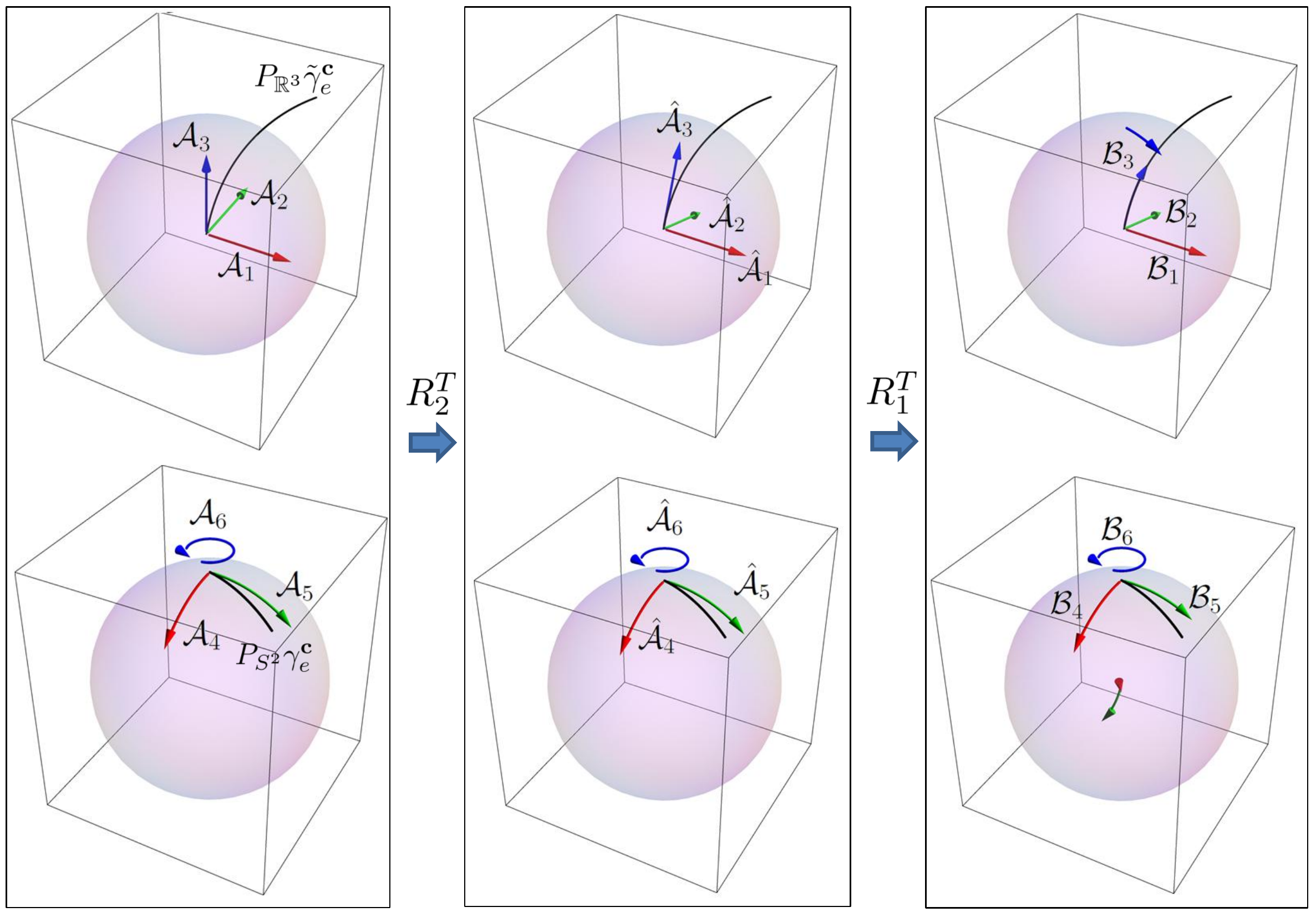}
}
\caption{{\small Visualization of the mapping of
left-invariant frame $\left.\{\mathcal{A}_{1},\ldots,\mathcal{A}_{6}\}\right|_{g}$ onto
locally adaptive spatial frame $\left.\{\mathcal{B}_{1},\ldots, \mathcal{B}_6\}\right|_{g}$
and $\tilde{\gamma}_{g}^{\ul{c}}(\cdot)$ a \emph{non-horizontal and torsion-free} exponential curve passing trough $g=\gothic{e}=(\ul{0},I)$.
The top row indicates the spatial part in $\R^3$, whereas the bottom row indicates the angular part in $S^2$. The top black curve is the spatial projection
of $\tilde{\gamma}_{g}^{\ul{c}}(\cdot)$, and the bottom black curve is the angular projection of the exponential curve.
After application of $\ul{R}_{2}^T$ the
exponential curve is horizontal w.r.t. the frame $\{\hat{A}_{i}\}$,
subsequently $\ul{R}_{1}^{T}$ leaves spatial generators orthogonal to a horizontal curve invariant, so that $\hat{A}_{1}=\mathcal{B}_{1}$
and $\hat{A}_{2}=\mathcal{B}_{2}$ are strictly spatial,
as is in accordance with Theorem~\ref{app:th}.
The angular part of $\mathcal{B}_{3}$ is shown via the curvature at the blue arrow. The spatial part of $\mathcal{B}_{4}$ and $\mathcal{B}_{5}$ is depicted in the center of the ball.
}
\label{fig:gauge3D}}
\end{figure}

\section{The Geometry of Neighboring Exponential Curves \label{app:B}}

In this appendix we provide some differential geometry underlying the family of neighboring exponential curves.

First we prove Lemma~\ref{lemma:SE3} on the construction of the family $\{\tilde{\gamma}^{\ul{c}}_{h,g}\}$ of neighboring exponential curves in $SE(3)$, recall Fig.~\ref{fig:pt}, and then we provide an alternative coordinate free definition of $\tilde{\gamma}^{\ul{c}}_{h,g}$ in addition to our Definition~\ref{def:neigbors} .

For the proof of Lemma~\ref{lemma:SE3}
we will just show equalities (\ref{3star}) as from this equality it directly follows by differentiation w.r.t. $t$  that the exponential curves $\tilde{\gamma}^{\ul{c}}_{h,g}(\cdot)=(\ul{x}_{h}(\cdot),\ul{R}_{h}(\cdot))$ and $\tilde{\gamma}^{\ul{c}}_{g}(\cdot)=(\ul{x}_{g}(\cdot),\ul{R}_{g}(\cdot))$ have the same spatial and angular velocity. For the spatial velocities it is obvious, for the angular velocities, we note that rotational velocity matrices $\OOmega_{h}$ and $\OOmega_{g}$ are indeed equal:
\[
\begin{array}{rl}
\ul{R}_{h}(t)=&\ul{R}_{g}(t) \ul{R}^{-1} \ul{R}' \Rightarrow \\
\OOmega_h:=&\left.\frac{d}{dt} \ul{R}_{h}(t)\right|_{t=0} (\ul{R}_h(0))^{-1}\\
 =&\left. \frac{d}{dt} \ul{R}_{g}(t)\right|_{t=0} (\ul{R}_g(0))^{-1}=\OOmega_{g},
\end{array}
\]
where we note that $\ul{R}_g(t)=e^{t \OOmega_g} \ul{R}_g(0) = e^{t \OOmega_g} \ul{R}$, and $\ul{R}_h(t)=e^{t \OOmega_h} \ul{R}_h(0) = e^{t \OOmega_h} \ul{R}'$.

Regarding the remaining derivation of (\ref{3star}), we note that it is equivalent to
\begin{equation}\label{3starb}
  \tilde{\gamma}_{h,g}^{\ul{c}}(t) =h \, (\ul{0}, (\ul{R}')^{T}\!\ul{R}) \, (g^{-1} \tilde{\gamma}_{g}^{\ul{c}}(t)) \, (\ul{0},(\ul{R}')^{T}\!\ul{R})^{-1}\!\!,
\end{equation}
by group product (\ref{groupproduct}). So we focus on the derivation of this identity. Let us set
$\ul{Q}:=(\ul{R}')^T\ul{R}$. Now relying
on the matrix representation (\ref{thedw}) and matrix exponential we deduce the following identity for $h=\gothic{e}=(\ul{0},I)$:
\begin{equation} \label{niceID}
\tilde{\gamma}_{\gothic{e},g}^{\ul{c}}(t)= \tilde{\gamma}_\gothic{e}^{ {\tiny \begin{pmatrix}
\ul{Q} && \ul{0} \\
\ul{0} && \ul{Q} \end{pmatrix}} \ul{c}}(t)= (\ul{0},\ul{Q})\,
\tilde{\gamma}_{\gothic{e}}^{\ul{c}}(t)\, (\ul{0},\ul{Q}^{-1}),
\end{equation}
which holds for all $\ul{Q} \in SO(3)$, in particular for $\ul{Q}=(\ul{R}')^T\ul{R}$.
Consequently, we have
\[
\tilde{\gamma}_{h,g}^{\ul{c}}(t)=h \tilde{\gamma}_{\gothic{e},g}^{\ul{c}}(t)=h\, (\ul{0},\ul{Q})\, g^{-1}
\tilde{\gamma}_g^{\ul{c}}(t) \, (\ul{0},\ul{Q}^{-1}),
\]
from which the result follows. $\hfill \Box$
\\ \\
We conclude From Lemma~\ref{lemma:SE3} that our Definition~\ref{def:neigbors} is indeed the right definition for our purposes, but as it is a definition expressed in left-invariant coordinates it also leaves the question
what the underlying coordinate-free unitary map from $T_{g}(SE(3))$ to $T_{h}(SE(3))$ actually is. Next we answer this question where we keep Eq.~\!(\ref{3starb}) in mind.
\begin{definition}\label{def:Uo}
Let us define the unitary operator \\$\gothic{U}^{h,g}:T_{g}(SE(3)) \to T_{h}(SE(3))$ by
\[
\begin{array}{rl}
\gothic{U}^{h,g} := (L_{h})_{*} \tilde{\mathbf{R}}_{h^{-1}g} (L_{g^{-1}})_{*},
\end{array}
\]
for each pair $g=(\ul{x},\ul{R}), h=(\ul{x}',\ul{R}') \in SE(3)$.
\end{definition}

\begin{remark}
From (\ref{unitary}) it follows that
the unitary correspondence between $T_{\tilde{\gamma}^{\ul{c}}_{g}(t)}$ and $T_{\tilde{\gamma}^{\ul{c}}_{h,g}(t)}$ is preserved \emph{for all} $t \in \R$.
\end{remark}
\begin{definition}\label{def:neighbouringCF}
The coordinate free definition of $\tilde{\gamma}^{\ul{c}}_{h,g}$ is that it is the
unique exponential curve passing through $h$ at $t=0$ 
with
\[
(\tilde{\gamma}_{h,g}^{\ul{c}})'(0)=\gothic{U}^{h,g}\; \left((\tilde{\gamma}_{g}^{\ul{c}})'(0)\right).
\]
\end{definition}
\begin{remark}
The previous definition (Definition \ref{def:neigbors}) follows from the coordinate free definition (Definition \ref{def:neighbouringCF}). This can be shown via identity (\ref{3starb}) which can be rewritten as
\begin{equation}
\tilde{\gamma}_{h,g}^{\ul{c}}(t) =L_{h} \circ \textrm{conj}(0,\ul{Q}) \circ (L_{g^{-1}}) \tilde{\gamma}_{g}^{\ul{c}}(t)),
\end{equation}
which indeed yields
\begin{equation*}
\begin{array}{ll}
(L_{h} \circ \textrm{conj}(0,\ul{Q}) \circ (L_{g^{-1}}))_{*}
&= (L_{h})_{*} \circ \textrm{Ad}(0,\ul{Q}) \circ (L_{g^{-1}})_{*}\\
 &= (L_{h})_{*} {\scriptsize \begin{pmatrix} \ul{Q} && \ul{0} \\
  \ul{0} && \ul{Q}  \end{pmatrix}} (L_{g^{-1}})_{*} \\
  &=\gothic{U}^{h,g},
\end{array}
\end{equation*}
again with $\mathbf{Q}=(\ul{R}')^{T}\ul{R}$, $\textrm{conj}(g)h=g h g^{-1}$, and $\textrm{Ad}(g)=(\textrm{conj}(g))_*$ is the adjoint representation \cite{Jost}.
\end{remark}


\section{Exponential Curve Fits on $SE(3)$ of the 2nd Order via Factorization \label{app:C}}

Instead of applying a 2nd order exponential curve fit (\ref{optHessianSumSE3}) containing a single exponential
one can factorize exponentials, and consider the following optimization:
\begin{equation}\label{alternativeSecondOrderSE3}
\boxed{ \begin{array}{l}
\ul{c}^{*}(g) =
\argmin_{\ul{c}\in \R^{6}, \|\ul{c}\|_{\mu}=1, c^6=0} \\
\left. \left|\;\frac{d^2}{dt^2} \tilde{V}(g\, e^{t(c^{1}A_{1}+ c^{2}A_{2}+c^{3} A_{3})}
e^{t(c^{4}A_{4}+ c^{5}A_{5})})\;\right|_{t=0} \right| .
\end{array}
}
\end{equation}
As shown in Theorem~\ref{th:8} the Euler-Lagrange equations are solved by spectral decomposition of the symmetric Hessian
given by
\begin{equation}
  \overline{\ul{H}}^{\ul{s}}:=\overline{\ul{H}}^{\ul{s}}(\tilde{U})=
\left(
\begin{array}{ccc}
\mathcal{A}_{1} \mathcal{A}_{1} \tilde{V} & \ldots & \mathcal{A}_{1} \mathcal{A}_{6} \tilde{V} \\
\vdots & \ddots & \vdots \\
\mathcal{A}_{1} \mathcal{A}_{6} \tilde{V} &\ldots & \mathcal{A}_{6} \mathcal{A}_{6} \tilde{V}
\end{array}
\right),
\end{equation}
with $\tilde{V}=\tilde{G}_{\ul{s}}*\tilde{U}$. This Hessian differs from the consistent Hessian in Appendix~\ref{app:B}.

\begin{theorem}[Second Order Fit via Factorization] \label{th:8}
Let $g \in SE(3)$ be such that Hessian matrix
$\ul{M}_{\mu}^{-1}(\overline{\ul{H}}^{\ul{s}}(g)) \ul{M}_{\mu}^{-1}$
has two eigenvalues with the same sign.
Then the normalized eigenvector $\ul{M}_{\mu}\ul{c}^{*}(g)$ with smallest eigenvalue
provides the solution $\ul{c}^{*}(g)$ of the following optimization problem (\ref{alternativeSecondOrderSE3}).

\end{theorem}
\textbf{Proof }
Define $F_{1}:=\ul{c}^{(1)} \cdot \ul{A}^{(1)} \in T_{e}(SE(3))$ with $\ul{A}^{(1)}\!:=(A_{1},A_{2},A_{3})^T$\!\!.
Define $F_{2}:=\ul{c}^{(2)} \cdot \ul{A}^{(2)} \in T_{e}(SE(3))$ with $\ul{A}^{(2)}\!:=(A_{4},A_{5},A_{6})^T$\!\!.
Define vector fields $\mathcal{F}_{1}|_g:= (L_{g})_* F_{1}$, $\mathcal{F}_{2}|_g:= (L_{g})_* F_{2}$. Then
Then
{\small
\[
\begin{array}{rl}
\left.\frac{d^2}{dt^2} \tilde{V}(g\, e^{t F_{1}} e^{t F_{2}})\right|_{t=0} \!\!&=\!
\lim \limits_{h \rightarrow 0} \! \frac{\tilde{V}\!(g  e^{h F_{1}} e^{h F_{2}})- 2 \tilde{V}\!(g) + \tilde{V}\!(g e^{-h F_{1}} e^{-h F_{2}})}{h^2} \\[6pt]
&=\mathcal{F}_{1}\mathcal{F}_{1} \tilde{V}(g) + \mathcal{F}_{2}\mathcal{F}_{2} \tilde{V}(g) +
2 \mathcal{F}_{1} \mathcal{F}_{2} \tilde{V}(g) \\
&= (\ul{c}(g))^T  \overline{\ul{H}}^{\ul{s}}(g)  \ul{c}(g).
\end{array}
\]
}
This follows by direct computation and the formula
\[
\tilde{V}(q e^{h F_{k}})= \tilde{V}(q) + h \mathcal{F}_{k}\tilde{V}(q) + \frac{h^2}{2} \mathcal{F}_{k}^{2}\tilde{V}(q) + O(h^3),
\]
applied for $(q=g e^{h F_{1}}, k=2)$ and $(q=g, k=1)$.

Therefore we can express the optimization functional as 
{\small
\begin{equation}\label{notsoeasy}
\begin{array}{rl}
\mathcal{E}(\ul{c})\!:=&\left|\left.\frac{d^2}{dt^2}
\tilde{V}(g e^{t(c^{1}A_{1}\!+\! c^{2}A_{2}\!+\!c^{3} A_{3})}
e^{t(c^{4}A_{4}\!+\!c^{5}A_{5})}
)\right|_{t=0}\right| \\
=&\left|\ul{c}^{T} \overline{\ul{H}}^{\ul{s}}(g)  \ul{c} \right|, \end{array}
\end{equation}
}
with again boundary condition $\varphi(\ul{c})=\ul{c}^T \ul{M}_{\mu}^2 \ul{c}=1$,
from which the result follows via Euler-Lagrange $\nabla \mathcal{E}=\lambda \nabla \varphi$ and
left multiplication with $\ul{M}_{\mu}^{-1}$.
$\hfill \Box$ \\

This approach can again be decomposed in the two-fold approach. Effectively, this means
that in Section \ref{ch:2fold2ndorder} the upper triangle of the Hessian $\ul{H}^{\ul{s}}$ is replaced by the lower triangle, whereas the lower triangle is maintained.
This approach performs well in practice; see e.g. Fig.~\!\ref{fig:SecondOrderFitSE3}
where the results of the exponential curve fits of second order are similar to exponential curve fits of first order.


%
\section{The Hessian induced by the left Cartan connection \label{app:remcohappy}}

In this section we will provide a formal differential geometrical underpinning for our choice of Hessian-matrix
\begin{equation}\label{HessianAppendix}
  \ul{H} (\tilde{U})= [\mathcal{A}_{j}(\mathcal{A}_{i}\tilde{U})],
\end{equation}
where $i$ denotes the row-index and $j$ the column index on $SE(d)$, recall the case $d=2$ in (\ref{ourhessian}) and recall the case $d=3$ in (\ref{fullhess}).
Recall from Theorem~\ref{th:0b}, Theorem~\ref{th:0c} and Theorem~\ref{th:6} that this Hessian naturally appears via direct sums or products in our exponential curve fits of second order on $SE(d)$.

Furthermore we relate our exponential curve fit theory to the theory in \cite{Jost}, where the same idea of 2nd order fits of auto-parallel curves to a given smooth function $\tilde{U}:M \to \R$ in a Riemannian manifold is visible in \cite[Eq.3.3.50]{Jost}.
Here we stress that in the book of Jost \cite[Eq.3.3.50]{Jost} this is done in the very different context of the torsion-free Levi-Civita connections, instead of the left Cartan connection which does have non-vanishing torsion.

Let us start with the coordinate free definition of the Hessian induced by a given a connection $\nabla^*$ on the cotangent bundle.
\begin{definition} (coordinate free definition Hessian)
On a Riemannian manifold $(M,G)$ with connection $\nabla^*$ on $T^{*}(M)$, the Hessian of
smooth function $\tilde{U}:M \to \R$ is defined coordinate independently (\cite[Def.3.3.5]{Jost}) by
$\nabla^{*} d \tilde{U}$.
\end{definition}
In coordinate-free form one has
(cf.~\cite[Eq.3.3.50]{Jost})
\begin{equation}\label{CoordinateFreeDefinitionHessian}
  \nabla^* d \tilde{U}(X_p,X_p)=\left.\frac{d^2}{dt^2} \tilde{U}(\gamma(t))\right|_{t=0}
\end{equation}
for the auto-parallel (i.e. $\nabla_{\dot{\gamma}}\dot{\gamma}=0$) curve $\gamma(t)$ with tangent $\gamma'(0)=X_p$ passing through $c(0)=p$ at time zero.
\begin{remark}
In many books on differential geometry $\nabla^{*}$ is again denoted by $\nabla$ (we also did this in our previous works \cite{QAM2,DuitsJMIV2}).
In this appendix, however, we distinguish between the connection $\nabla$ on the tangent bundle and its adjoint connection $\nabla^*$ on the co-tangent bundle $T^*(SE(d))$.
\end{remark}

Let us recall that the structure constants of the Lie algebra are given by
\begin{equation} \label{struct}
[\mathcal{A}_{i},\mathcal{A}_{j}]=\mathcal{A}_{i}\mathcal{A}_{j}-\mathcal{A}_{j}\mathcal{A}_{i}=\sum \limits_{k=1}^{n_d}c^{k}_{ij} \mathcal{A}_{k}.
\end{equation}
As shown in previous work \cite{QAM2} the left Cartan connection\footnote{also known as minus Cartan connection}
$\nabla$ on $M=(SE(d),\gothic{G}_{\mu})$, is the (metric compatible) connection whose Christoffel symbols, expressed in the left-invariant moving
(co)frame of reference, are equal to the structure constants of the Lie algebra:
\[
\Gamma^{k}_{ij}=c_{ji}^k=-c_{ij}^k \in \{-1,0,1\}.
\]
More precisely, this means that if we compute the covariant derivative of a vector field $Y=\sum \limits_{k=1}^{n_d}y^{k}\mathcal{A}_{k}$ (i.e. a section in $T(SE(d))$
along the tangent $\dot{\tilde{\gamma}}(t)=\sum \limits_{i=1}^{n_d} \dot{\tilde{\gamma}}^i(t)\, \left.\mathcal{A}_{i}\right|_{\tilde{\gamma}(t)}$ of some smooth curve $t \mapsto \tilde{\gamma}(t)$ in $SE(d)$.
This is done as follows
\begin{equation} \label{CartanVF}
\nabla_{\dot{\tilde{\gamma}}} Y = \sum \limits_{k=1}^{n_d} \left(\dot{y}^k - \sum \limits_{i,j=1}^{n_d} c^{k}_{ij} \dot{\tilde{\gamma}}^i y^j\right) \mathcal{A}_{k},
\end{equation}
where we follow the notation in Jost's book \cite[p.108]{Jost} and define $\dot{y}^{k}(t):=\frac{d}{dt} y^{k}(\tilde{\gamma}(t))$.
By duality this induces the following (adjoint) covariant derivative of a covector field $\lambda$ (i.e. a section in $T^{*}(SE(d))$):
\begin{equation} \label{CartanCVF}
\nabla_{\dot{\tilde{\gamma}}}^{*} \lambda = \sum \limits_{i=1}^{n_d} \left(\dot{\lambda}_{i} + \sum \limits_{k,j=1}^{n_d} c^{k}_{ij} \lambda_{k} \dot{\tilde{\gamma}}^{j} \right) \omega^i,
\end{equation}
with $\dot{\lambda}_{i}(t) =\frac{d}{dt} \lambda_{i}(\tilde{\gamma}(t))$.
Then by antisymmetry of the structure constants it directly follows (see e.g.~\cite{DuitsJMIV2}) that the auto-parallel curves are the exponential curves:
\begin{equation}\label{AutoParallel}
  \nabla_{\dot{\tilde{\gamma}}}\dot{\tilde{\gamma}}=0 \textrm{ and }
\dot{\tilde{\gamma}}(0)=\ul{c} \textrm{ and } \tilde{\gamma}(0)=g
\desda \tilde{\gamma}=\tilde{\gamma}^{\ul{c}}_g.
\end{equation}
\begin{remark}
Due to torsion of the left Cartan connection, the auto-parallel curves do not coincide with the geodesics w.r.t. metric tensor $\gothic{G}_{\xi}$. This is in contrast to the Levi-Cevita connection (see for example Jost's book \cite[ch:3.3, ch:4]{Jost}) where auto-parallels are precisely the geodesics (see \cite[ch:4.1]{Jost}).

Intuitively speaking this means that in the curved geometry of the left Cartan connection on $SE(d)$ (that is present in the domain of an orientation score, see Figure~\ref{Fig:Intro}) the `straight curves' (i.e. the auto-parallel curves) do not coincide with the `shortest curves' (i.e. the Riemannian distance minimizers).
\end{remark}
The left Cartan connection is the consistent connection on $SE(d)$ in the sense that auto-parallel curves are the exponential curves studied in this article.
Therefore the consistent Hessian form on $SE(d)$ is induced by the left Cartan connection.
Expressing it in the left-invariant frame yields
\begin{equation}\label{CoordinateFreeHessianDerivation}
\boxed{
\begin{array}{rcl}
\nabla^* {\rm d}\tilde{U}(\mathcal{A}_{i},\mathcal{A}_{j}) & \overset{\textrm{def}}{=} &
(\nabla_{\mathcal{A}_{i}}^* {\rm d}\tilde{U})(\mathcal{A}_{j}) \\
  &\overset{(\ref{CartanCVF})}{=} &
\sum \limits_{j'=1}^{n_d} (\mathcal{A}_i \mathcal{A}_{j'} \tilde{U} + \sum \limits_{k=1}^{n_d} c^{k}_{j'i} \mathcal{A}_{k}\tilde{U}) \,\omega^{j'}(\mathcal{A}_{j}) \\
  &\overset{(\ref{coframe})}{=} & (\mathcal{A}_i \mathcal{A}_{j} \tilde{U} + \sum \limits_{k=1}^{n_d} c^{k}_{ji} \mathcal{A}_{k}\tilde{U}) \\
  &\overset{(\ref{struct})}{=} & (\mathcal{A}_i (\mathcal{A}_j \tilde{U}) + (\mathcal{A}_j\mathcal{A}_i-\mathcal{A}_i\mathcal{A}_j )\tilde{U})  \\
  & = & (\mathcal{A}_j (\mathcal{A}_i \tilde{U}),
\end{array}
}
\end{equation}
where $i$ denotes the row-index and $j$ the column index. So we conclude from this computation that (\ref{HessianAppendix}) is the correct consistent Hessian on $SE(d)$ for our purposes.

\begin{remark}
The  left Cartan connection has torsion and is not the same as the standard torsion-free Cartan-Schouten connection on Lie groups, which have also many applications in image analysis an statistics on Lie groups, cf. \cite{Pennec1,Pennec2}. Recall that within the orientation score framework, right invariance is undesirable.
\end{remark}

\begin{table*}
\renewcommand{\arraystretch}{1.3}
\section{Table of Notations \label{app:TableOfNotations}}
\begin{tabular}{|l|p{9.25cm}|p{5.4cm}|}
\hline
\textbf{Symbol}
& \textbf{Explanation}
& \textbf{Reference}
\\ \hline \hline
\multicolumn{3}{c}{} \vspace{-3mm} \\
\multicolumn{3}{c}{E.1 Spaces and Input Data}
\\ \hline
$SE(d)$
& The group of rotations and translations on $\mathbb{R}^d$
& Section \ref{ch:intro}, Section \ref{ch:SEd}, and (\ref{groupproduct})
\\ \hline
$\mathbb{R}^{d} \rtimes S^{d-1}$
& Space of positions $\&$ orientations as a group quotient in $SE(d)$
& (\ref{PositionsAndOrientations}), and Section \ref{ch:introquotient}.
\\ \hline
$\tilde{U}$
& Input data $\tilde{U}:SE(d) \to \mathbb{R}$
& (\ref{tildeU}), and Section \ref{ch:SEd}
\\ \hline
$U$
& Input data $U:\mathbb{R}^{d} \rtimes S^{d-1} \to \mathbb{R}$
& (\ref{tildeU}), (\ref{PositionsAndOrientations}), and Section \ref{ch:introquotient}.
\\ \hline
$\tilde{V}$
& Gaussian smoothed input data $\tilde{V}=\tilde{G}_{\mathbf{s}}* \tilde{U}$
& (\ref{VEET}), and (\ref{Gtilde})
\\ \hline
\multicolumn{3}{c}{} \vspace{-3mm} \\
\multicolumn{3}{c}{E.2 Tools from Differential Geometry}
\\ \hline
$\left.\mathcal{A}_{i}\right|_{g}$
& Left-invariant vector field $\mathcal{A}_{i}$ restricted to $g \in SE(d)$
& Section \ref{ch:LeftInvariantDerivatives}, and (\ref{LINV}),  (\ref{LINVDef})
\\ \hline
$\left.\mathcal{B}_{i}\right|_{g}$
& Gauge vector field $\mathcal{B}_{i}$ restricted to $g \in SE(d)$
& Section \ref{ch:togauge}, and (\ref{gauge})
\\ \hline
$\left.\gothic{G}_{\mu}\right|_{g}$
& Metric tensor $\gothic{G}_{\mu}$ restricted to $g \in SE(d)$
& Section \ref{ch:Metric}, and (\ref{metrictensor2})
\\ \hline
$\|\cdot\|_{\mu}$
& $\mu$-norm on $\mathbb{R}^{n_d}$, with \ $n_d=\textrm{dim}(SE(d))=\frac{d(d+1)}{2}$
& Section \ref{ch:Metric}, and (\ref{normmu})
\\ \hline
$\ul{M}_{\mu}$
& Matrix $\ul{M}_{\mu}:= \left( \begin{smallmatrix}
   \mu I_d  && \ul{0} \\
   \ul{0} && I_{r_d} \end{smallmatrix} \right)$ is used in definition of the $\mu$-norm $\|\cdot\|_{\mu}$
& Section \ref{ch:Metric}, and (\ref{normmu})
\\ \hline
${\rm d}\tilde{U}(g)$
& Derivative of $\tilde{U}$ at $g$ which is a covector in $T^*_{g}(SE(d))$
& (\ref{ExteriorDerivative})
\\ \hline
$\nabla\tilde{U}(g)$
& Gradient of $\tilde{U}$ at $g$ which is a vector in $T_{g}(SE(d))$
& Section~\ref{ch:GaussianSmoothing} and (\ref{gradientmu2})
\\ \hline
$\chi$
& Deviation from horizontality angle
& (\ref{dH})
\\ \hline
$\mathcal{L}$
& Left regular representation given by $\mathcal{L}_{g}\tilde{U}(h)=\tilde{U}(g^{-1}h)$
& (\ref{eq:leftrightregularrepr})
\\ \hline
$\mathcal{R}$
& Right regular representation given by $\mathcal{R}_{g}\tilde{U}(h)=\tilde{U}(h\, g)$
& (\ref{eq:leftrightregularrepr})
\\ \hline
$L$
& Left multiplication $L_{g} h= g h$
& (\ref{leftmultiplication})
\\ \hline
\multicolumn{3}{c}{} \vspace{-3mm} \\
\multicolumn{3}{c}{E.3, part I: Exponential Curves and Exponential Curve Fits on SE(d)}
\\ \hline
$\tilde{\gamma}^{\mathbf{c}}_{g}(\cdot)$
& Exponential curve starting from $g$ with velocity $\ul{c}=(\ul{c}^{(1)},\ul{c}^{(2)})$
& (\ref{expcurvend}), and Section \ref{ch:ExponentialCurves}
\\ \hline
$\tilde{\gamma}^{\mathbf{c}}_{h,g}(\cdot)$
& Neighboring exponential curve starting at $h \in SE(d)$ with the same velocity as curve $\tilde{\gamma}^{\mathbf{c}}_{g}$
&
(\ref{paralleltransportSE2}) and Section \ref{ch:FamilyExponentialCurvesSE2}, (\ref{paralleltransportSE3version2}) and Section \ref{ch:FamilyExponentialCurvesSE3}
\\ \hline
$\tilde{\gamma}^{{\tiny \mathbf{c}^*\!(g)}}_{g}(\cdot)$
& Exponential curve fit to data $\tilde{U}$ at $g \in SE(d)$
& (\ref{expcurvend}), and Theorem \ref{th:0},\ref{th:0b},\ref{th:0c},\ref{th:2},\ref{th:3b},\ref{th:6}, Fig. \ref{Fig:Intro}
\\ \hline
$\mathbf{c}^{*}(g)$
& Local tangent vector to exponential curve fit $\tilde{\gamma}^{\mathbf{c}^*}_{g}$ to data $\tilde{U}$
& (\ref{minproblemSE2firstorder}), (\ref{SE2secondOrderOptimization}), (\ref{SE2secondOrderOptimizationGradient}), (\ref{optSE3}), and (\ref{optHessianSumSE3})
\\ \hline
$\tilde{\mathbf{R}}_{h^{-1}\!g}$
& Rotation in $T_h(SE(d))$ arising in the construction of $\tilde{\gamma}^{\ul{c}}_{h,g}$
& (\ref{RtildeSE2}), and (\ref{RtildeSE3})
\\ \hline
$\mathbf{S}^{\mathbf{s},\boldsymbol{\rho}}$
& Structure tensor $\mathbf{S}^{\mathbf{s},\boldsymbol{\rho}}:=\mathbf{S}^{\mathbf{s},\boldsymbol{\rho}}(\tilde{U})$
of input data $\tilde{U}$
& (\ref{fig:SecondOrderFitSE3}), (\ref{structSE3}), and (\ref{structR3XS2})
\\ \hline
$\mathbf{H}^{\mathbf{s}}$
& Gaussian Hessian $\mathbf{H}^{\mathbf{s}}:= \mathbf{H}^{\mathbf{s}}(\tilde{U})=\mathbf{H}(\tilde{V})$ of input data $\tilde{U}$
& (\ref{ourhessian}), and (\ref{fullhess})
\\
\hline
\multicolumn{3}{c}{} \vspace{-3mm} \\
\multicolumn{3}{c}{E.3, part II: Exponential Curve Fits on SE(3) with Projections in $\R^{3}\rtimes S^2$}
\\ \hline
$\mathbf{R}_{\mathbf{a},\phi}$
& Counter-clockwise 3D rotation about axis $\mathbf{a}$ by angle $\phi$
& text below (\ref{lifttogroup})
\\ \hline
$\mathbf{R}_{\mathbf{n}}$
& Any 3D rotation that maps $\mathbf{a}=(0,0,1)^T$ onto $\mathbf{n} \in S^2$
& (\ref{lifttogroup}), and Theorem \ref{th:3b}
\\ \hline
$h_{\alpha}$
& Element $h_{\alpha}=(\ul{0},\ul{R}_{\ul{a},\alpha})$ of the subgroup $\equiv \{\ul{0}\} \times SO(2)$
& (\ref{HALPHA})
\\ \hline
$\odot$
& Symbol denoting action of SE(3) onto $\mathbb{R}^{3} \rtimes S^{2}$
& (\ref{odot}), and Section \ref{ch:introquotient}
\\ \hline
$\mathbf{Z}_{\alpha}$
& Rotation matrix in $SO(6)$ that arises in $\nabla \tilde{U}$ if $g \mapsto g h_{\alpha}$
& (\ref{Zalphadef})
\\ \hline
$\mathcal{N}$
& Null-space of the structure tensor $\mathbf{S}^{\mathbf{s},\boldsymbol{\rho}}$
& (\ref{nullspace})
\\ \hline
$\gamma^*_{(\mathbf{y},\mathbf{n})}(\cdot)$
& Projected exponential curve fit to data $U$ at $(\mathbf{y},\mathbf{n}) \in \mathbb{R}^{3} \rtimes S^{2}$
& (\ref{fig:SecondOrderFitSE3}), and (\ref{optcurve})
\\ \hline
$g_{new}$
& Location in $SE(3)$ for horizontal exponential curve fit
& (\ref{nnew})
\\ \hline
\multicolumn{3}{c}{} \vspace{-3mm} \\
\multicolumn{3}{c}{E.4 Applications}
\\ \hline
$f$
& Input greyscale image $f:\mathbb{R}^{d} \to \mathbb{R}$
& (\ref{OS}), and (\ref{OSreconstruction})
\\ \hline
$\mathcal{W}_{\psi}f$
& Orientation score of greyscale image $f$ via cakewavelet $\psi$
& (\ref{OS}), and Fig. \ref{Fig:Intro}, \ref{fig:OS3D}, \ref{fig:cakewavelets}
\\ \hline
$\Phi_t$
& Nonlinear diffusion operator (diagonal diffusion in gauge frame)
& (\ref{gaugeflow})
\\ \hline
$\tilde{W}(g,t)$
& Scale space representation of $\tilde{U}$ at $g \in SE(d)$ and scale $t>0$
& (\ref{gaugeflow})
\\ \hline
$\Phi$
& Vesselness operator
& (\ref{SE2vesselness})
\\ \hline
\multicolumn{3}{c}{} \vspace{-3mm} \\
\multicolumn{3}{c}{E.5 Appendix}
\\ \hline
$\nabla_{\dot{\tilde{\gamma}}} Y$
& Covariant derivative of vector field $Y$ along $\dot{\tilde{\gamma}}$ w.r.t. Left-Cartan connection $\nabla$ on $T(SE(d))$
& (\ref{CartanVF}), and Appendix \ref{app:remcohappy}
\\ \hline
$\nabla^*_{\dot{\tilde{\gamma}}} \omega$
& Covariant derivative of covector field $\omega$ along $\dot{\tilde{\gamma}}$ w.r.t. the adjoint Left-Cartan connection $\nabla^*$ on $T^*(SE(d))$
& (\ref{CartanCVF}), and Appendix \ref{app:remcohappy}
\\ \hline
$\nabla^* {\rm d}\tilde{U}$
& Coordinate-free definition of the Hessian
& (\ref{CoordinateFreeDefinitionHessian}), (\ref{CoordinateFreeHessianDerivation}), and Appendix \ref{app:remcohappy}
\\ \hline
\end{tabular}
\end{table*}

\end{document}